%% file: main.tex
\documentclass[12pt,letterpaper]{article}
\usepackage[utf8]{inputenc}
\usepackage[english]{babel}
\usepackage{amsfonts}
\usepackage{amssymb}
\usepackage[pdftex]{graphicx}
\usepackage[font=small,labelfont=bf]{caption}
\usepackage{epstopdf}
\usepackage[T1]{fontenc}
\usepackage{etoolbox}
\usepackage{multirow}
\usepackage[top=1.3in, bottom=1.3in, left=0.8in, right=0.8in]{geometry}
\usepackage{tikz}
\usetikzlibrary{patterns}
\usepackage{subcaption}
\usepackage{tikz-3dplot}
\usepackage{array}
\usepackage{rotating}
\usepackage{ctable}
\usepackage{geometry}
\usepackage{color,soul}
\usepackage{float}
\usepackage{wrapfig}
\usepackage{listings}
\usepackage{comment}
\usepackage{mathtools}
\usepackage{tabularx}
\usepackage{url}
\usepackage{enumitem}
\usepackage{longtable}
\usepackage{algorithm}
\usepackage[noend]{algpseudocode}
\usepackage{hyperref}
\usepackage{cleveref}
\usepackage{amsmath}

\usepackage{color, colortbl}
\definecolor{mygreen}{RGB}{28,172,0} 
\definecolor{mylilas}{RGB}{170,55,241}
\definecolor{YellowHighlight}{rgb}{1,0.929,0.525}
\definecolor{RedHighlight}{rgb}{0.9412,0.502,0.502}
\definecolor{GreenHighlight}{rgb}{0.5961,0.9843,0.5961}
\definecolor{niall}{RGB}{46,139,87}
\definecolor{sasha}{RGB}{68,119,170}

\newcolumntype{M}[1]{>{\arraybackslash}m{#1}}
\newcolumntype{N}{@{}m{0pt}@{}}

\newcommand{\mathbfgl}[1]{\boldsymbol{\mathbf{#1}}}

\allowdisplaybreaks

\makeatletter
\def\BState{\State\hskip-\ALG@thistlm}
\makeatother

\usetikzlibrary{arrows.meta}

\usepackage{calc}
\usepackage{soul}

\begin{document}

\newlength\Dwidth
\settowidth{\Dwidth}{\textbf{Sensors/Ladders}}
\newcolumntype{C}{>{\centering}p{(\Dwidth-4\tabcolsep)/3}}

\lstset{language=Matlab,%
    breaklines=true,%
    morekeywords={matlab2tikz},
    keywordstyle=\color{blue},%
    morekeywords=[2]{1}, keywordstyle=[2]{\color{black}},
    identifierstyle=\color{black},%
    showstringspaces=false,
    numbers=left,%
    numberstyle={\tiny \color{black}},
    numbersep=9pt, 
    basicstyle=\linespread{0.9}\small\ttfamily,
}

\newcommand{\ipnp}{\ensuremath{_{i+1}^{n+1}}}
\newcommand{\imnp}{\ensuremath{_{i-1}^{n+1}}}
\newcommand{\inp}{\ensuremath{_{i}^{n+1}}}

\newcommand{\ipn}{\ensuremath{_{i+1}^{n}}}
\newcommand{\imn}{\ensuremath{_{i-1}^{n}}}
\newcommand{\inn}{\ensuremath{_{i}^{n}}}

\newcommand{\np}{\ensuremath{^{n+1}}}
\newcommand{\n}{\ensuremath{^{n}}}

\newcommand{\erf}{\ensuremath{\text{erf}}}

\newcommand{\etaltri}{\emph{et al.}$\;$} 

\newcommand{\red}{\textcolor{red}}
\newcommand{\blue}{\textcolor{blue}}

\newcommand{\ud}{\mathop{}\!\mathrm{d}}

\title{Model selection of chaotic systems from data with hidden variables using sparse data assimilation}

\author{H. Ribera, S. Shirman, A. V. Nguyen, N. M. Mangan}

\date{\today}


\maketitle


\begin{abstract}

Many natural systems exhibit chaotic behaviour such as the weather, hydrology, neuroscience and population dynamics. Although many chaotic systems can be described by relatively simple dynamical equations, characterizing these systems can be challenging, due to sensitivity to initial conditions and difficulties in differentiating chaotic behavior from noise. Ideally, one wishes to find a parsimonious set of equations that describe a dynamical system. However, model selection is more challenging when only a subset of the variables are experimentally accessible. Manifold learning methods using time-delay embeddings can successfully reconstruct the underlying structure of the system from data with hidden variables, but not the equations. Recent work in sparse-optimization based model selection has enabled model discovery given a library of possible terms, but regression-based methods require measurements of all state variables. We present a method combining variational annealing -- a technique previously used for parameter estimation in chaotic systems with hidden variables --  with sparse optimization methods to perform model identification for chaotic systems with unmeasured variables. We applied the method to experimental data from an electrical circuit with Lorenz-system like behavior to successfully recover the circuit equations with two measured and one hidden variable. We discuss the robustness of our method to varying noise and manifold sampling using ground-truth time-series simulated from the classic Lorenz system.

\end{abstract}

\clearpage

\begin{center}
    \textbf{Significance statement}
\end{center}

Chaos represents a challenge for studying the dynamic behavior of many physical and biological systems. Since the 80s we have known that time-series measurements from one variable of a chaotic system contain information about the underlying structure of the full multi-dimensional system. However, recovery of the full system from data with hidden variables has remained elusive. This work develops a novel data-assimilation technique to identify governing equations of chaotic systems from data with hidden variables. This method identifies fairly simple, low-dimensional, and deterministic models from seemingly incomplete data. Discovery of such equations can enable rich mathematical study and physical insight for problems across nearly every discipline including climate science, hydrology, neuroscience, ecology, medicine and engineering.


\clearpage





\addtocontents{toc}{\protect\setcounter{tocdepth}{0}}

\section{Introduction}

Hypothesis generation through data-driven model identification has the potential to revolutionise science. Uncovering the interactions, structure, and mechanisms that determine the behaviour of chaotic systems in particular could improve scientific understanding in almost every discipline with dynamical systems \cite{Gardini2020} including climate science \cite{Slingo2011}, hydrology \cite{Sivakumar2000}, population dynamics \cite{Hassell1991}, and neuroscience \cite{Rabinovich1998}. Many chaotic systems can be informatively described by relatively simple dynamical equations. However, characterization and control of these systems can be challenging \cite{Boccaletti2000}, due to sensitivity to initial conditions and difficulties in differentiating chaotic behavior from noise \cite{Sugihara1990}. Characterization through statistical, geometric, or model-based means becomes more challenging when only a subset of the variables are experimentally accessible. Our goal is to identify a parsimonious set of equations to describe a chaotic system from measurements with hidden variables.

Much data-analysis for chaotic systems has focused on learning the attracting manifold structure from time-series. In the early 80s, Takens’s theorem \cite{Takens1981} describes the conditions under which one can use the time-delay embedding from a single variable to construct a manifold that preserves the topological properties of the full system. Takens's result formalized the idea that the information of the manifold structure, and therefore chaotic dynamics, could be recovered from the time-history of a single state variable. Manifold reconstruction methods \cite{Packard1980,Fraser1986,Kennel1992} based on partial information provide insight into the system structure, dimensionality, and statistics of chaotic systems. By constructing manifolds from time-delays, Sugihara \etaltri developed methods discriminating chaos from noise \cite{Sugihara1990} and detecting causality between measured variables \cite{Sugihara2012}.
Methods including reservoir computing \cite{Tang2020, Fan2020},  other deep learning frameworks \cite{Yeo2017}, data assimilation combined with neural networks \cite{Brajard2020}, support vector machine \cite{Mukherjee1997}, and nearest neighbours \cite{Altman1992} can accurately predict the dynamics of chaotic systems using a data-trained model with no specific physical knowledge of the system. For a review of predictive methods see \cite{Amil2019}. Assuming a reasonable model structure is known, data-assimilation methods \cite{Bannister2017, Bezruchko2006} including variational annealing \cite{Ye2015} can estimate model parameters for chaotic systems from incomplete, indirect, and noisy measurements. Although these methods are designed to assimilate information from data-streams with hidden variables and learn about chaotic systems, they are not designed for the purpose of hypothesizing parsimonious models or identifying model structure.

Data-driven discovery of parsimonious dynamical systems models to describe chaotic systems is by no means new. Early on, least-squares fitting of combinatorial sets of polynomial basis functions to time-series data followed by information-theory based selection produced models that reproduced manifold structure and statistics of the system \cite{Crutchfield1987}. Symbolic regression demonstrated successful recovery of the widely accepted equations for the chaotic double-pendulum system \cite{Schmidt2009}. More recently sparse regression \cite{Tibshirani1996,Candes2008,Hesterberg2008}, motivated model selection techniques such as SINDy \cite{Brunton2016}, which recover the ground-truth equations for chaotic systems from a relatively large library of functions, without needing a computationally intensive combinatorial search. Other sparsity-promoting frameworks have improved upon robustness for chaotic systems equation recovery through integral formulations \cite{Schaeffer2017,Pantazis2018,Messenger2020}, data assimilation methods \cite{Bocquet2019}, Bayesian frameworks \cite{Bocquet2020}, and entropic regression \cite{AlMomani2020}. However, all these methods require measurements of all state-variables that significantly impact the desired dynamic. Notably, Champion \etaltri recently used an autoencoder framework for automatic discovery system coordinates and equations, but required input time-series of a higher dimension than intrinsic dimension of the system \cite{Champion2019}. 

Model selection with hidden variables require different methodology. By `hidden variables' we mean that the number of measured variables is smaller than the intrinsic dimension of the system. Measured variables are not considered hidden if they are corrupted by noise or indirectly sampled through a measurement function. A few methods address the problem of model selection with hidden variables, but they have not been demonstrated for chaotic systems. For example, Daniels \etaltri \cite{Daniels2015,Daniels2019} combinatorially fit each model in a predefined model space using data assimilation and subsequently use Bayesian inference to select the best model. Successful recovery of mass-action kinetic systems for chemical reactions was demonstrated with hidden variables using a neural network approach \cite{Ji2020}.
A recent method uses LASSO to select predictive models for chaotic systems from a library with higher order derivatives given a single state variable \cite{Somacal2020}. This method effectively finds a higher-order ODE representation of the Lorenz and R\"{o}ssler systems, but it is unclear how the recovered structures relate to the ground truth models.

In this paper we present a new method to perform model selection in dynamical systems with hidden variables. This method combines the data assimilation technique variational annealing, which has been used to estimate parameters when the structure of the system is known, with sparse model selection via hard thresholding. We call this method Data Assimilation for Hidden, Sparse Inference (DAHSI). To demonstrate that our method could identify interpretable models for chaotic systems, we followed the philosophy of earlier works \cite{Schmidt2009, Brunton2016} and demonstrated recovery of accepted parsimonious models from experimental data and simulated time-series where the ground truth is known. In the Results section DAHSI successfully selected a set of models for a circuit that has Lorenz-like behaviour from experimental data of two state variables (one hidden). One of the identified models has the same structure as the Lorenz system. The other identified models with high AIC/BIC support exhibit nearly indistinguishable dynamics and suggest novel terms which may better represent the experimental circuit system. Moreover, we used ground truth simulations of the canonical Lorenz system to study how our method performs with varying data size and noise. In the Materials and Methods section we describe the DAHSI algorithm for model selection with hidden variables.  


\section{Results: Model selection for chaotic systems}
\label{sec-results}

\subsection{Identification of models for the Lorenz circuit from experimental data}
\label{sec-circuit}

\begin{figure}[H]
    \centering
    \includegraphics[width=0.91\textwidth]{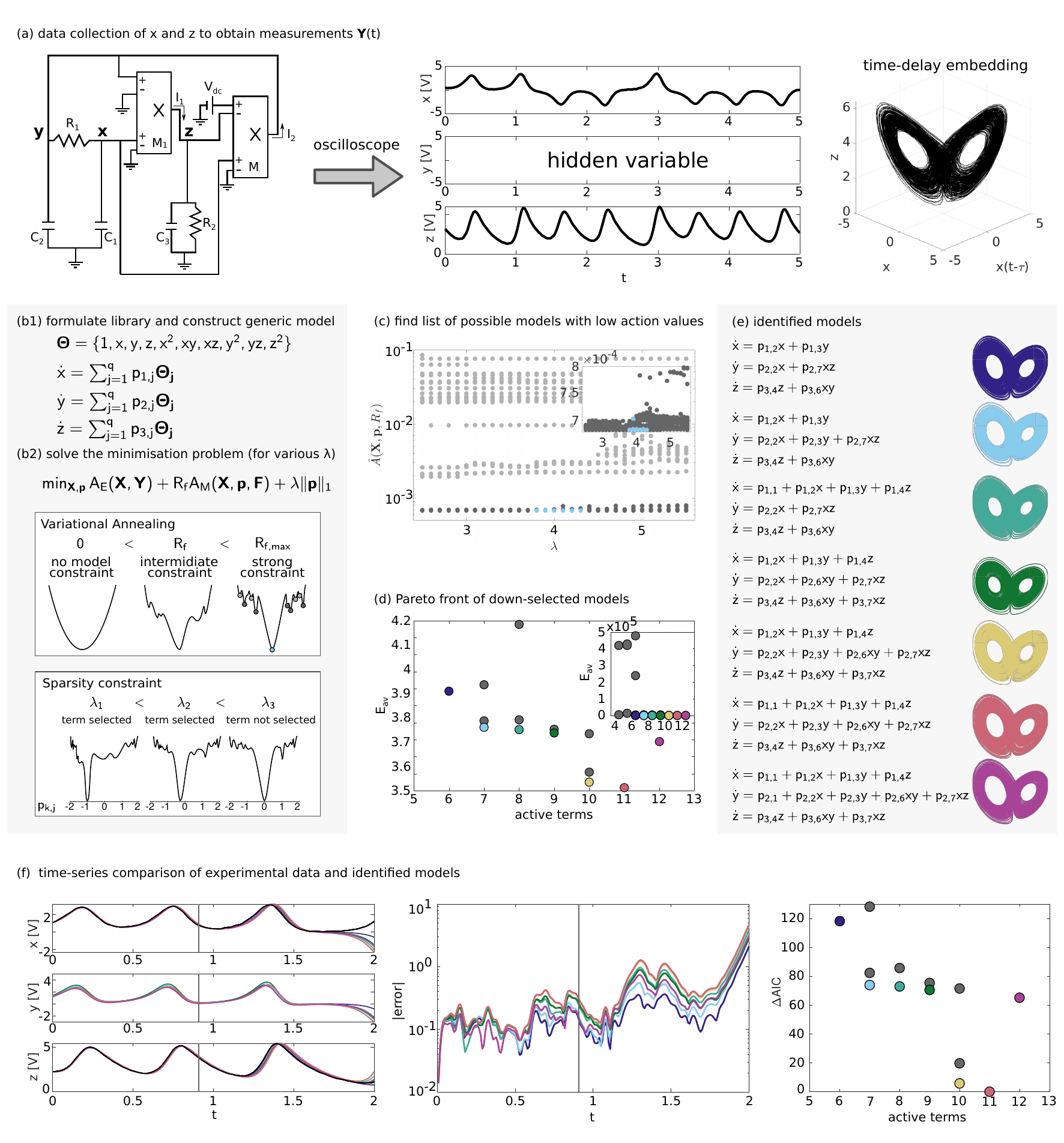}
    \caption{DAHSI model selection for the Lorenz-like system. 
    (a) Electrical circuit from \cite{Blakely2007}, training data of measured variables $x$ and $y$, and time-delay embedding of test data ($\tau = 0.02$).
    (b1) Model library and generic governing equation for each variable.
    (b2) Cost function as model error weight and the sparsity constraint vary.
    (c) Local minima with high cost (light grey), low cost (dark grey) and Lorenz structure (blue) as function of $\lambda$.
    (d) 25 low cost models are down-selected.
    (e) Model structure identified near the Pareto front.
    (f) Time series, error, and relative AIC for identified models (coloured), and higher error models (grey).
    }
    \label{figure_Scheme}
\end{figure}

The Lorenz system \cite{Lorenz1963} was originally developed to forecast the weather and has become a canonical example when developing new methods to characterize chaotic systems. To demonstrate model selection on experimental data with hidden variables, we considered high-quality data from the electrical circuit in Blakely \etaltri \cite{Blakely2007} (Fig. \ref{figure_Scheme}(a)). This system exhibits similar structure and behavior to the highly studied Lorenz system and is well described by relatively simple circuit equations 
\begin{align}
    \label{rescaledx}
    \frac{\ud x}{\ud t} &= \hat{\sigma} (y-x),\\
    \frac{\ud y}{\ud t} &= \hat{\rho}x - \hat{\gamma}y - \hat{\varepsilon} xz, \\
    \label{rescaledz}
    \frac{\ud z}{\ud t} &= -\hat{\beta} z + \hat{\eta} xy.
\end{align}
The structure of this system is similar to the Lorenz system, but in the standard Lorenz formulation $\hat{\varepsilon} = \hat{\eta} = \hat{\gamma}$. Here, $\mathbf{X} = (x,y,z)$ denote the voltages across the capacitors $C_1$, $C_2$ and $C_3$ in the circuit (Fig. \ref{figure_Scheme}(a)). The measured variables are $x$ and $z$, and $y$ is unmeasured or hidden. We denote the noisy measurements of $x$ and $z$ by $x_e$ and $z_e$, respectively, and the measurement function $\mathbf{h}((x,z)) = (x_e,z_e) = \mathbf{Y}$.
The experimental sampling rate is $\Delta t_e = 80$ ns resulting in $55,000$ time points. A low-pass filter was applied to remove high-frequency measurement error \cite{Blakely2007}. 
We re-scaled the experimental time by $\Delta t = \frac{\Delta t_e}{ 1.6 \times 10^{-5} \textrm{ns}} = 0.005$ so that the right hand side terms of \eqref{rescaledx}-\eqref{rescaledz} are around $\mathcal{O}(1)$.
We trained our method with $N=501$ time points (Fig. \ref{figure_Scheme}(a)), at a sampling rate $2 \Delta t = 0.01$ (re-scaled). The attractor is reasonably well sampled with 501 points (SI Appendix, Fig. 2), and we retain the remaining data for validation.

We demonstrated model identification with hidden variables of the Lorenz-like system (\eqref{rescaledx}-\eqref{rescaledz}) using DAHSI (Fig. \ref{figure_Scheme}).
First, we constructed a model library based on domain-knowledge. In this case we used monomials up to degree two in three variables, representing $10^9$ possible models composed of subsets of possible terms. From this library we generated a generic governing equation for each variable via the linear combination of all the candidate functions (Fig. \ref{figure_Scheme}(b1)). Our goal was to find a small subset of candidate functions within the library which describe the dynamics in the data. We did not assume that we knew the ``correct'' model complexity {\sl a priori}, and searched for the set of models which balance error and simplicity.

To perform model selection, we minimised a cost function composed of the measurement error, $A_E$, model error, $A_M$, and sparse penalty, $\lambda \Vert \mathbf{p} \Vert_1$ as a function of the parameters, $\mathbf{p}$ and library of functions $\mathbfgl{\Theta}$ (Fig. \ref{figure_Scheme} (b), and Materials and Methods section). The model error contains the coupling between variables, taking advantage of the information about hidden variables in the time-history of the measured variables. The measurement error only depends on the measurements and measured variables estimated from the model. Model selection is enabled through the sparse penalty which determines the number of parameters, $p_{k,j}$, that will be active in the model or zero.

To minimize the cost function, we used variational annealing (VA) \cite{Ye2015}, a data-assimilation technique for non-convex parameter estimation in nonlinear, chaotic systems. The problem is highly non-convex, with many local minima, due to the incoherence between the data and the model \cite{Pecora1990,Abarbanel2013}. Decreasing the information or measurements \cite{Kadakia2017} and increasing the number of terms in the library will both increase the number of local minima (SI Appendix, Fig. 1). VA works by varying  $R_f$ which sets the balance between model error and measurement error (Fig. \ref{figure_Scheme}(b2)). When $R_f = 0$, only measurement error contributes leading to a convex cost function with an easy to find global minima. As the model is enforced by gradually increasing $R_f$, the landscape increases in complexity and many local minima appear. By initialising the search near the minima for the previous $R_f$ the solution remains near the first minima found. Varying $\lambda$ leads to different model structures or candidate models. As the penalty strength, $\lambda$, increases, the global minima moves to 0 in a larger number of parameters (Fig. \ref{figure_Scheme}(b2)). Because there are many local minima, we need to choose $N_I = 500$ random initial guesses to fully explore the landscape. 

The sparse-variational annealing process generates 169 candidate models, which must be further down-selected and validated to complete the model-selection process. We down-selected to the 25 models (SI Appendix) with a cost function value less than $10^{-3}$  (Fig. \ref{figure_Scheme}(c)). 
In our system there is a clear gap in cost-function value at this value, but the criteria and gap size will be system dependent. To ensure we have the best parameter fit for each down-selected model we performed parameter estimation via VA without sparsity constraint. 

To validate the models, we needed to estimate an initial condition for the hidden variable $y$, for which there is no experimental data. We used an 8th order finite difference approximation of the time derivative of $x$ for each model structure and solve the resulting algebraic equation for $y_{0}$ (SI Appendix). We used the dynamic equation for $x$ since all down selected models contain $y$ but not any higher order $y$ terms. Estimation of the initial condition for hidden variables is only possible after the candidate models are found and must be done for the initial condition of each segment of validation data. This procedure takes advantage of Takens's theorem that the information in $y$ is available in the time-delay of $x$.

Validation within the Lyapunov time ensures that the time-series do not diverge due to the inherent sensitivity to differences in initial conditions introduced by measurement and numerical error. 
All down-selected models have a similar Lyapunov time around 0.9 time units. 
We considered $S = 1083$ segments of the experimental data (excluding the training set), each of length 1/4 of a Lyapunov time to calculate the sum of the average error for each model (Fig. \ref{figure_Scheme}(d)). We discarded the first four points of each time segment as these points were used to predict the initial condition for $y$. The average error for the $s$-th time segment of the $m$-th model is defined as $E^s_{av,m} = \frac{1}{2M} \sum_{i=1}^M (x_{i,e}^s - x_i^s)^2 + (z_{i,e}^s - z_i^s)^2$, where $x_{i,e}$ and $z_{i,e}$ are the $x$ and $z$ components of the experimental data, respectively, and $i$ is the time index. The sum of all average errors over the time segments $S$ of the $m$-th model is $E_{av,m} = \sum_{s=1}^S E^s_{av,m}$. 

The candidate models on the Pareto front (Fig. \ref{figure_Scheme}(d), and SI Appendix, Table S1) best balance model complexity and error (Fig. \ref{figure_Scheme}(e)). We successfully recovered the Lorenz-like structure derived by Blakely \etaltri \cite{Blakely2007}, which has the lowest average error of recovered models with 7 active terms. For $\lambda = 3.9$ the system presented in \cite{Blakely2007} is selected for 10.6\% of the $N_I = 500$ randomly chosen initialisation. However, we have no guarantee that this model is the "true model" for the circuit system. 
All models have a similar manifold structure (Fig. \ref{figure_Scheme}(e)) and low error within a Lyapunov time (Fig. \ref{figure_Scheme}(f)). We believe the main limitation of our prediction window is the uncertainty introduced by the hidden variable into the parameter estimation during VA. This uncertainty then propagates into the $y_0$ estimate required for each validation data set and magnifies noise (SI Appendix, Figs. S4 and S5).
Given the difficulty in selecting between proposed chaotic models that exhibit such similar behaviour \cite{Aguirre1994,Aguirre200}, the primary goal of DAHSI as a model identification method is to generate possible models. However, we were able to consistently identify a unique model (salmon with 11 terms, Fig. \ref{figure_Scheme}(f)) with the most support using Akaike information criteria as done in \cite{Mangan2017} and Bayesian information criteria (SI Appendix, Fig. 6), as well as identifying a weakly supported model (gold with 10 terms, Fig. \ref{figure_Scheme}(f)). 
By generating multiple models that lie near the Pareto front DAHSI has effectively generated hypothesis for additional terms, which could be tested with further experimentation.

While DAHSI identified the same equation terms as Lorenz and the circuit formulation from \cite{Blakely2007}, the parameters fit through the final step of VA are not the same. We compare the ability to predict the experimental data with the classical Lorenz system, the circuit formulation from \cite{Blakely2007} and the DAHSI-recovered models, each of which have a different number of free parameters (Table \ref{table_params_comparison}).
We perform parameter estimation via VA for each model and use the validation data-set described above to calculate $E_{av}$, $\Delta$AIC, and $\Delta$BIC. Although the average error is similar for the VA-estimated circuit formulation and all DAHSI models, the DAHSI recovered models with 10 and 11 terms have substantially more $\Delta$AIC and $\Delta$BIC support. The classical Lorenz parameter structure, which only has 4 free parameters, is unable to capture the dynamics of the system. The parameters estimated via VA for the circuit model with 6 free parameters perform much better than those estimated from first principles \cite{Blakely2007}. Notably, the parameters estimated for the 7-term DAHSI model are very close to the parameters estimated for the original circuit model. Further experimentation is needed to determine if the coefficients in the $\dot{x}$ equation should be equal, $p_{1,2} = p_{1,3}$, and if the coefficient on the $y$ term in the $\dot{y}$ equation, $p_{2,3}$ should be positive, negative, or zero (Fig. \ref{figure_Scheme}(e)). The additional terms suggested by the 10 and 11 term DAHSI recovered models are strongly supported by the AIC/BIC calculations, but would require further experimentation to conclusively validate. They may represent parasitic resistances or other physical effects which have a small but real impact on the circuit dynamics and were neglected during the original derivation by Blakely {\sl et al.}  \cite{Blakely2007}. Recovery of the Lorenz-like model and identification of other models with AIC/BIC support demonstrates that DAHSI can successfully identify parsimonious models for chaotic systems.

\begin{table}[H]
\centering
\caption{Parameter estimation for the classical Lorenz formulation (4 free parameters, $p_{1,2} = p_{1,3}$,  $p_{2,3 } = p_{2,7} = -p_{3,6}$); the circuit formulation in \cite{Blakely2007} (6 free parameters, $p_{1,2} = p_{1,3}$); and the DAHSI-recovered models.}
\label{table_params_comparison}
    \resizebox{0.98\textwidth}{!}{%
    \begin{tabular}{l|r|c|r|rr|rrr}
    & & & & \multicolumn{2}{c|}{circuit formulation} & \multicolumn{3}{c}{DAHSI-recovered} \\
    \hline
     & Term & Parameter & classical & as in \cite{Blakely2007} & estimated & 7-terms & 10-terms & 11-terms \\ \hline
    
    \multirow{4}{*}{eq. $\dot{x}$} & $1$ & $p_{1,1}$ & -- & -- & -- & -- & -- & $-0.2514$ \\
    & $x$ & $p_{1,2}$ & $-29.7560$ & $-12.9032$ & $-16.5369$ & $-16.9554$ & $-17.0172$ & $-17.0582$ \\
    & $y$ & $p_{1,3}$ & $29.7560$ & $12.9032$ & $16.5369$ & $18.7853$ & $19.9884$ & $19.9840$ \\  
    & $z$ & $p_{1,4}$ & -- & -- & -- & -- & $0.1596$ & $0.1833$ \\  
    \hline

    \multirow{4}{*}{eq. $\dot{y}$} & $x$ & $p_{2,2}$ & $68.5427$ & $54.2903$ & $28.0876$ & $24.3535$ & $22.6028$ & $22.6017$ \\    
    & $y$ & $p_{2,3}$ & $-12.8815$ & $-1.2903$ & $-0.0763$ & $0.2580$ & $0.3346$ & $0.3567$ \\ 
    & $xy$ & $p_{2,6}$ & -- & -- & -- & -- & $-0.0906$ & $-0.0843$ \\ 
    & $xz$ & $p_{2,7}$ & $-12.8815$ & $-14.2857$ & $-7.6252$ & $-6.7054$ & $-6.2507$ & $-6.2561$ \\     
    \hline
    
    \multirow{2}{*}{eq. $\dot{z}$} & $z$ & $p_{3,4}$ & $-3.4168$ & $-3.8259$ & $-3.6547$ & $-3.6835$ & $-3.6954$ & $-3.6966$ \\ 
    & $xy$ & $p_{3,6}$ & $12.8815$ & $3.4843$ & $4.3315$ & $4.8273$ & $5.1412$ & $5.1292$ \\ 
    & $xz$ & $p_{3,7}$ & -- & -- & -- & -- & $0.0903$ & $0.0791$ \\     
    \hline
    $E_{av}$ & -- & -- & 2165 & 319 & 10.37 & 9.7441 & 9.0995 & 9.0345 \\ 
    \hline
    $\Delta$AIC & -- & -- & 5920 & 3852 & 139.315 & 73.887 & 5.758 & 0 \\ 
    \hline
    $\Delta$BIC & -- & -- & 5885 & 3827 & 114.377 & 53.937 & 0.77 & 0
    \end{tabular}
    }
\end{table}

\subsection{Robustness study on the simulated Lorenz system}

To study the robustness of our method to varying noise and manifold sampling we used ground-truth time series simulated from the classic Lorenz system,
\begin{align}
    \label{eq_for_x}
    \frac{\ud x}{\ud t} &= \sigma(y-x),\\
    \frac{\ud y}{\ud t} &= x(\rho-z) - y,\\
    \label{eq_for_z}
    \frac{\ud z}{\ud t} &= -\beta z + xy,
\end{align}
where $\sigma = 10$, $\rho = 28$, and $\beta = 8/3$. We numerically simulated the system using Runge-Kutta 4th order and a time step of $\Delta t = 0.01$ and $N=501$, producing time-series similar to the experimental data set. As in the experimental data set, we considered $y$ to be the hidden variable. We studied the recovery rates of DAHSI as a function of the VA tuning parameter, $\alpha$, and found trends similar to previous work \cite{Rozdeba2018}, (SI Appendix, Table S3).

First, we studied the robustness of our method to measurement error modeled as additive Gaussian noise of mean zero and varying standard deviation, $\omega$. Therefore, the measurement function is $\mathbf{h}(\mathbf{X}) = \mathbf{X} + \mathcal{N}(0,\omega)$. We expect that different noise instances, controlled by the random number generator seed, will change our recovery rate due to random corruption of essential parts of the data or overall poor manifold sampling.

\begin{figure}[H]
    \centering
    \includegraphics{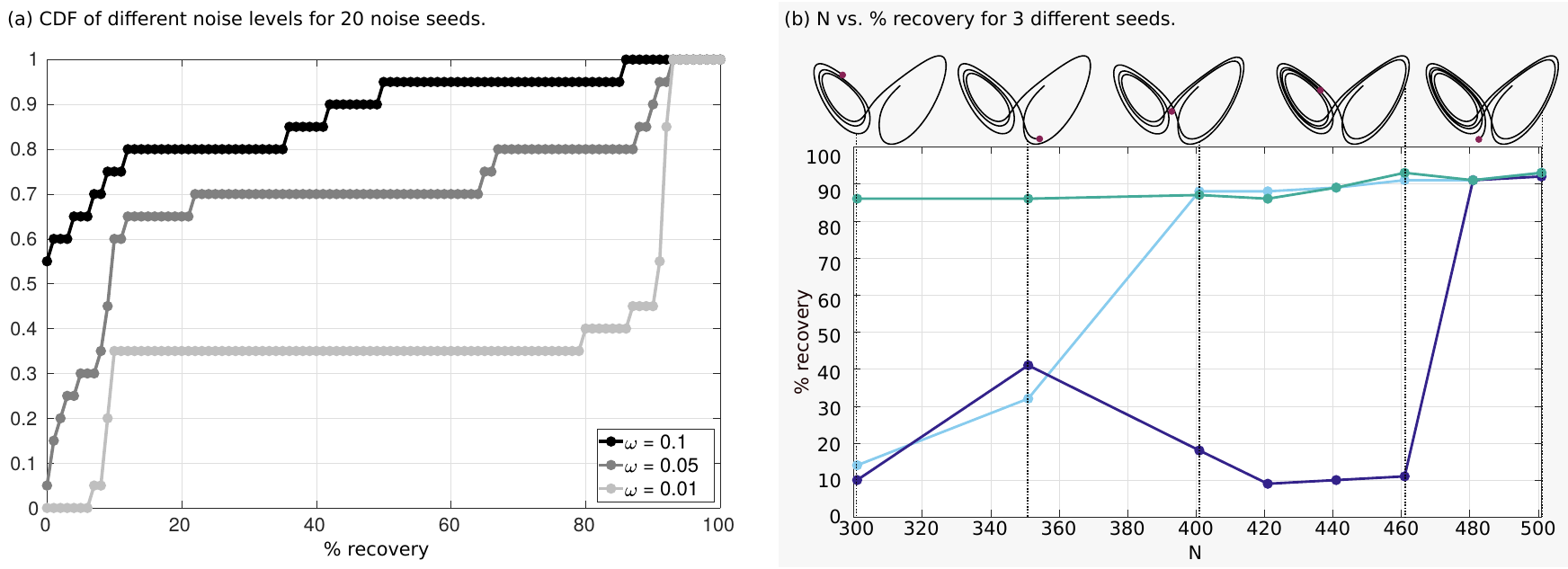}
    \caption{
    Robustness to noise and manifold sampling.
    (a) Cumulative distribution function of the recovery rate for three noise levels and 20 noise seeds. $\omega = 0.01$ (light grey); $\omega = 0.05$ (dark grey); $\omega = 0.1$ (black).
    (b) Recovery rate on eight different manifold sizes, for three different noise seeds (colors).
    }
    \label{figure_Robustness}
\end{figure}

We calculated recovery for 3 different standard deviations of noise with 20 noise seeds each and calculate the cumulative distribution function of the recovery rate (Fig. \ref{figure_Robustness}(a)). The random noise seeds produced wide variation in recovery rate between 10-90\% for the lowest noise, indicating that the minimal data set used here is not very robust.
As the noise strength increased, the cumulative distributions shifted left as more seeds have lower recover. Setting $\omega = 0.01$ produced a binomial distribution, with either a high recovery rate (> 80\%) (the majority of simulations), or a low recovery rate (< 15\%). 
For $\omega = 0.05$ there were some seeds with intermediate recovery rates, more low recovery rates, and a few seeds that with a very high recovery rate. 
The noise level dramatically affected the recovery rate for $\omega = 0.1$. The vast majority of simulations led to less than 10\% recovery. More than half had 0\% recovery, and only one had higher than 85\% recovery.

Next, we investigated how manifold sampling affected the recovery rate of our system. We chose 3 different noise seeds, and varied the number number of time points $N$ by increasing the length of the time-series (Fig. \ref{figure_Robustness}(b)). Varying the length of the time-series changed the sampling of the manifold, demonstrating that sampling lobe transitions is crucial for accurate model recovery.
For one seed (light blue line) the recovery was high for $N = 501$ through $N = 401$. 
There were sharp drops in recovery of $\approx 60$\% and $\approx 15$\% when the data-set lost a lobe crossing in the attractor, as happens at $N = 351$ and $N = 301$, respectively.
Sharp drops in another seed (dark blue line) also occurred when the sampling of the crossing between lobes is reduced at $N=460$ and $301$. Decreased sampling of each lobe did not appear to have as dramatic an effect ($N = 401$ to $460$). The increase of recovery rate for the dark-blue noise instance at $N = 351$ suggests that optimal sampling requires some nontrivial balance of different dynamic regions.
The specific corruption of noise instance had a big impact on how many crossings are needed to get a high recovery as the recovery was consistently high for one seed (cyan). These results suggest that optimal manifold sampling to counter noise corruption would vastly improve DAHSI performance on data sets with high noise.

\section{Discussion}
\label{sec-discussion}

In this paper we have presented DAHSI, a method to identify non-linear dynamical systems from data with hidden variables. DAHSI combines variational annealing, a data assimilation technique, with sparse thresholding.
We applied DAHSI to an experimental data set from a circuit that exhibits Lorenz-like dynamics \cite{Blakely2007}. The outcome is a set of candidate models, including a model with the same Lorenz-like structure derived by Blakely \etaltri from circuit equations \cite{Blakely2007}. Two additional parsimonious models with strong support based on AIC/BIC-based validation were also identified. The unanticipated terms suggested by these models may represent real physical processes in the circuit, such parasitic resistances or other factors not included in the idealized model derivation. Through this example, we demonstrated that DAHSI works as an effective tool for generating models and functional hypothesis from data.

To analyze recovery and the effects of noise and manifold sampling in a system where we know the ground truth, we studied the performance of DAHSI applied to simulated time-series from the classical Lorenz system.
Notably, we successfully selected the ground truth model as most likely from those generated by DAHSI using information-criteria based validation techniques (SI Appendix. Fig. 3).
Our noise studies showed recovery rates of 80\% for $\sim \mathcal{N}(0,0.01)$ and  10\% or lower for  $\sim \mathcal{N}(0,0.1)$. Therefore we anticipate that the current formulation of DAHSI will have reasonable recovery rates for noise levels $< 10\%$ of the signal value.
Further robustness to noise could be achieved through integral formulations similar to those used for sparse regression, rather than the discretized mapping between time-points used here \cite{Schaeffer2017,Pantazis2018,Messenger2020}.
Manifold sampling impacts recovery and we conclude that recovery is especially sensitive to sampling at the saddle point transition between the lobes. 
Moreover, the noise seed used to generate the synthetic data impacts the recovery and we suspect this is due to random corruption of measurements from different regions of the manifold.  For chaotic systems, increasing the time of experiment will eventually ensure robust sampling of the manifold. However, the computational time of DAHSI scales with the length of the input time-series \cite{Gondzio2012}.
Therefore, we anticipate that short bursts of time-series designed to optimally sample the manifold would provide optimal sampling and computational efficiency. Further metrics for analyzing the information content of our data and minimal data-requirements for recovering models \cite{Ho2020} would lead to optimal manifold sampling. 

One of the main benefits of a sparse model selection framework is that we identify likely model structures while avoiding combinatorial testing and validation of all potential models. For example, the number possible models described three variables with monomials up to degree two is approximately $10^9$. Doing parameter estimation on each of these models and validating would be computationally intractable, taking at least $10^5$ processor-days with our setup. For comparison, our entire model selection and validation process took just over a day of computational time. Running one initialisation of the problem and sweeping through $\lambda = 2.5:0.1:5.5$ with $N=501$ (as done in Example \ref{sec-circuit}) took 4 hours. We parallelized simulations using Northwestern's High Performance Computing Cluster Quest, running about 100 of simulations at a time, leading to a total computational time of roughly 20 hours. Performing parameter estimation without thresholding on a single model takes between 15 second and 15 minutes, depending on model structure. Parameter estimation on 25 down-selected models took 5 hours with our set up. Times estimates are for a Intel(R) Xeon(R) CPU E5-2680 v4 @ 2.40GHz processor. In order to understand the impact of library size on a call to IPOPT, the optimiser used in DAHSI, we tested model libraries with 7, 10, 13, 16, 19, and 30 terms (SI Appendix, Fig. 7). The computational time does not scale monotonically with library size. Instead, we find that a library with 10 terms can take 100 times longer to run than the library of 30 monomials. We suspect that the variation in optimization time depends on correlations between library functions \cite{Mangan2019}, model symmetries, and other structural features. 

In addition to the chaotic systems presented in the results, we have applied DAHSI on two non-chaotic systems: on time-series data from a Lotka-Volterra-like system with no hidden variables and on simulated time-series for a mass action kinetics systems with hidden variables. Although DAHSI recovered reasonable models for both systems, there are several caviats. Recover of Lotka-Volterra required an iterative formulation (SI Appendix, Fig. 15). We also compared DAHSI to SINDy \cite{Brunton2016} for the Lotka-Volterra system and found that SINDy was far superior in speed when all variables are observed. Recall that a comparison between DAHSI and SINDy is not possible for Lorenz-like circuit system, as SINDy requires access to the unmeasured $y$ variable. 
The mass action kinetic system modeled a semiconductor with two trap levels differing by one electronic unit of charge (SI Appendix). The recovery rate for the ground truth model was low, around 3\%. Unlike chaotic systems, which are highly non-convex, the mass-action kinetic system has a very flat cost function due to structural parameter identifiability issues (SI Appendix, Figs. S11-S14), 
\cite{Apgar2010,Gabor2017,Meshkat2009,Eisenberg2014}. Stochastic gradient decent algorithms such as IPOPT are known to perform poorly for flat cost functions so switching to an optimiser designed for such systems \cite{Kantabutra2002} may improve recovery. Other data-assimilation methods for parameter estimation with hidden variables such as 3D-Var, 4D-Var, Kalman filtering, and hybrid methods \cite{Bannister2017} may be more cost-effective if VA is unnecessary to navigate to the global minimum of a highly non-convex function.

The formulation of cost function and sparsity constraint also likely impacts recovery. Different methods for sparse model-selection include stepwise and all-subsets regression, ridge regression \cite{Hoerl1970}, LASSO \cite{Tibshirani1996}, least angle regression \cite{Efron2004}, and SR3 \cite{Zheng2018}. SR3 accelerates convergence and has been shown to outperform other methods and improves performance but has an extra tuning parameter.
The parameter path for the first four methods is shown to be different in \cite{Hesterberg2008} and therefore, we expect that different regularisation methods will lead to different model identification. Comparison between different sparsity-enforcement mechanisms within DAHSI framework could improve recovery but may be somewhat system dependent.

We anticipate many future applications and extensions of DAHSI. The framework for DAHSI does not have any intrinsic restrictions about the functional form of the equations, in particular the function library need not be linear in the unknown parameters. Variational annealing is designed to handle stochasticity through the model error. In addition, data assimilation is commonly used for PDE systems, including PDE discovery \cite{Chang2019}. Therefore, we anticipate we can apply or extend our framework to broader applications, without reformulation as was needed in sparse-regression based frameworks for rational functions \cite{Mangan2016}, stochastic systems \cite{Boninsegna2018}, and PDEs \cite{Rudy2017,Kang2019}. Modifications to the optimization methodology and further investigation of optimal data-sampling strategies could improve the computational efficiency of DAHSI, opening up higher dimensional problems to model selection with hidden variables.


\section{Methods: Mathematical formulation of cost function and algorithm}
\label{sec-mat-and-met}

The dynamics of many physical systems can be described by models with only a few terms. Our goal is to retrieve the sparse system representation of these type of systems given the measurements of some, but not all, of the state variables. We consider a dynamical system with unknown governing equations
\begin{equation}
\label{dynamics-2}
    \frac{\ud \mathbf{X}}{\ud t} = \mathbf{F}(\mathbf{X}(t),\mathbf{p}),
\end{equation}
where $\mathbf{X} = (x_1,x_2,\dots,x_D) \in \mathbb{R}^D$ are the state variables, $\mathbf{F} = (F_1,\,F_2,\dots,F_D)$ are the unknown functions that govern the dynamics of the system and $\mathbf{p}$ is a set of unknown parameters.

For a system with hidden variables, the measurements $\mathbf{Y} = (y_1,y_2,\dots,y_L) \in \mathbb{R}^L$ are lower dimensional $L \leq D$ than the underlying variables. The measurement function $\mathbf{h}(\mathbf{X}) = \mathbf{Y}$ is a known transformation of a subset of the state variables in \eqref{dynamics-2}. 
In principle, the measurement function could map some combination of state variables to a lower-dimension, as in $\mathbf{h}(\mathbf{X})= x_1 + x_2$. In this work we assume $\mathbf{h}$ captures Gaussian experimental noise such that, $\mathbf{Y} = \mathbf{X} + \mathcal{N}(0,\omega)$. The measurements are taken at $N$ equally spaced point in time between $[t_1,\,t_N]$.

The function capturing the nonlinear dynamics of each state variable, $F_k$, is assumed to be sparse in function space as has been done previously \cite{Hastie2009, Brunton2016}. Given a library of possible functions $\mathbfgl{\Theta} = (\theta_1,\theta_2,\dots,\theta_q)$, we can write a candidate function $\hat{F}_k$ as
\begin{equation}
    \label{generic-model}
    \hat{F}_k \coloneqq \hat{F}_k(\mathbf{X},\mathbf{p}) = p_{k,1} \theta_1(\mathbf{X}) + p_{k,2} \theta_2(\mathbf{X}) + \cdots + p_{k,q} \theta_q(\mathbf{X}),
\end{equation}
for $k=1,2,\dots,D$. There is no inherent restriction that the functions be linearly additive. The set of $p_{k,j}$ defines the vector $\mathbf{p} \in \mathbb{R}^{P}$, where $P = Dq$ is the total number of unknown parameters. 

We want to estimate the unknown parameters $p_{k,j}$ and all state variables $\mathbf{X}$ using only the measurements $\mathbf{Y}$ with the constraint that $\mathbf{p}$ is sparse. This is equivalent to minimising the negative log likelihood
\begin{equation}
    \label{costfunk}
    \begin{split}
        &A(\mathbf{X},\mathbf{p}) = \frac{1}{N}\sum_{i=1}^N \Vert \mathbf{X}(t_{i}) - \mathbf{Y}(t_{i}) \Vert^2 \\
        & + \frac{1}{N}\sum_{i=1}^{N-1} R_f \left\{ \Vert \mathbf{X}(t_{i+1}) - \mathbf{f}(\mathbf{X}(t_i),\mathbf{p},\mathbf{\hat{F}}) \Vert^2 \right\} + \lambda \Vert \mathbf{p} \Vert_1.        
    \end{split}
\end{equation}

Here, $\mathbf{f}(\mathbf{X}(t_i),\mathbf{p},\mathbf{\hat{F}}) = \mathbf{X}(t_{i+1})$ defines the discrete time model dynamics and is obtained by discretising \eqref{dynamics-2} using a Hermite-Simpson collocation. We note that if $\lambda = 0$ in \ref{costfunk} we obtain the cost function used in VA.
Following the statistical derivation in \cite{Talagrand1987,Evensen2009,Abarbanel2013}, the experimental error, $A_E(\mathbf{X},\mathbf{Y}) = \frac{1}{N}\sum_{i=1}^N \Vert \mathbf{X}(t_{i}) - \mathbf{Y}(t_{i}) \Vert^2$ assumes Gaussian noise and the model error, $A_M(\mathbf{X},\mathbf{p},\mathbf{\hat{F}}) = \frac{1}{N}\sum_{i=1}^{N-1} \left\{ \Vert \mathbf{X}(t_{i+1}) - \mathbf{f}(\mathbf{X}(t_i),\mathbf{p},\mathbf{\hat{F}}) \Vert^2 \right\}$ assumes a relaxed delta function. We assume that the state at the $t_{i+1}$ depends only on the state at $t_i$. We assume that each element in $\mathbf{p}$ follows a Laplace distribution with \emph{mean} $0$ (SI Appendix). The details and necessary background to minimise \eqref{costfunk} are presented in the following sections.

\subsection{DAHSI: Data Assimilation for Hidden Sparse Inference}
\label{sec-DAHSI}

\begin{figure}[H]
    \centering
    \includegraphics[width=0.68\textwidth]{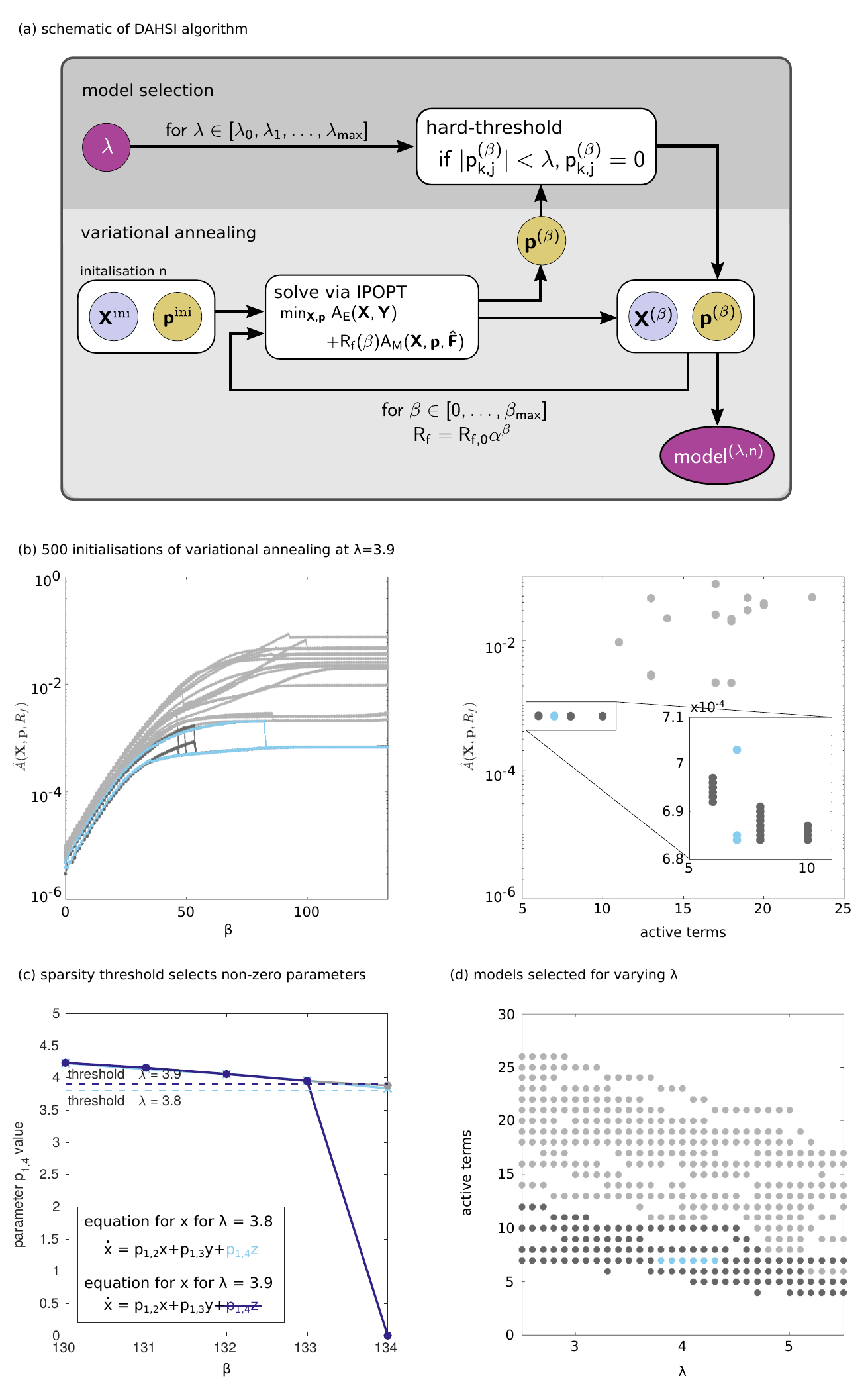}
    \caption{DAHSI Algorithm. 
    (a) Schematic of Algorithm \ref{algorithm}. 
    (b) Action paths as function of $\beta$ for $N_I=500$ and$\lambda=3.9$ (left). Final action values (right) for high (light grey) and low (dark grey) action values; and the Lorenz-like structure (blue).
    (c) Parameter $p_{1,4}$ in the last steps of VA for $\lambda = 3.8$ and $3.9$ 
    (d) Model complexity as function of $\lambda$.
    }
    \label{figure_MethodLambdamodels}
\end{figure}

Our algorithm, Data Assimilation for Hidden Sparse Inference (DAHSI), performs model identification for chaotic systems from data with hidden variables. It combines the data assimilation technique VA with sparse thresholding (Fig. \ref{figure_MethodLambdamodels}(a)). The code base for DAHSI can be found at \cite{Ribera2021GitHub}.

As the desired model complexity is unknown ahead of time, DAHSI sweeps through different hard-threshold values, $\lambda$. For each $\lambda$, the cost function \eqref{costfunk}, is minimized by iterating between VA \cite{Ye2015,Ye2015b} and hard-thresholding of the parameters. We chose the iterative framework over direct incorporation of the $\ell_1$ penalty into the minimized cost function, based on the results that show that least square with thresholding converges locally, often outperforming convex variants \cite{Zheng2018,Chartrand2008}, and recent demonstrations that LASSO makes mistakes early in the recovery pathway \cite{Su2017}. 

At each VA step, we minimize $A_E + R_f A_M$, which is 4DVar in its "weak" formulation \cite{Talagrand1987,Evensen2009}, over $\mathbf{X}$ and $\mathbf{p}$ given $R_f$ using IPOPT, an optimisation package that uses a gradient descent method \cite{Wachter2002}. The state variables $\mathbf{X}^{\text{ini}}$ are initialized as $\mathbf{Y}$ for the measured states and random values from a uniform distribution within specified bounds for the unmeasured states. Since we expect the parameter vector $\mathbf{p}$ to be sparse, it is initialized as $\mathbf{p^{\text{ini}}} = 0$. 

Initially $R_f$ takes some small value $R_{f,0}= \epsilon$, as $R_f = 0$ would lead to an unconstrained solution on the unmeasured states and $\mathbf{p}$.
At each step $\beta = 0,1,2,\dots,\beta_{\max}$ of VA, $R_f$ is updated to $R_f = R_{f,0} \alpha^{\beta}$, for $\alpha>1$. After each step $\beta$ of VA, we enforce sparsity by applying a hard threshold, $\lambda$, to $\mathbf{p}^{(\beta)}$.
The solution, $\{\mathbf{X}^{(\beta)},\mathbf{p}^{(\beta)}\}$, at each step of the VA process is used as the initialization for the next step. 
We choose $\beta_{\max}$ so that the cost function plateaus, Fig. \ref{figure_MethodLambdamodels}(b), and our final solution is $\{\mathbf{X}^{\text{fin}},\mathbf{p}^{\text{fin}}\}$. 
Because there are many local minima, we run $N_I$ different initial guesses to fully explore the landscape of $A_E + R_f A_M$.
It is important to note that the same $\lambda$ yields multiple models due to the $N_I$ different initializations of the unmeasured states. For example, if we consider $N_I = 500$ with a fixed $\lambda = 3.9$ in our Example \ref{sec-circuit}, we find a total of 20 models (Fig. \ref{figure_MethodLambdamodels}(b)). 

To produce candidate models with varying sparsity, the entire $\beta$ sweep with VA and thresholding is repeated for each $\lambda$.
As with other model identification methods, different $\lambda$ will yield different models (for the same initialisation of unmeasured states). For one particular initialisation in Example \ref{sec-circuit}, with $\lambda = 3.8$ the term $z$ is selected in the first equation of the system. With larger $\lambda = 3.9$, the term $z$ is no longer selected (Fig. \ref{figure_MethodLambdamodels}(c)). Although the same $\lambda$ yields multiple models due to the difference of the initial choice of unmeasured states, as we would expect, higher values of $\lambda$ produce models with fewer active terms (Fig. \ref{figure_MethodLambdamodels}(d)). 

\begin{algorithm}
\caption{DAHSI Algorithm.}
\label{algorithm}    
    \begin{algorithmic}[1]
    
    \Procedure{DAHSI}{}
        \State Input: measurements $\mathbf{Y}$, generic model library $\mathbf{\Theta}$, $\lambda_{max}$, $\beta_{max}$, $\alpha$
        \State Calculate discrete function $\mathbf{\hat{F}}$ from $\mathbf{\Theta}$
        \For{$l = 1:L$}
            \State $x_l = y_l$  \Comment{Fit measurements to data}
        \EndFor
        \State Randomly initialise unobserved variables $\{x_{l+1},\dots,x_D\}$
        \State $\mathbf{X}^{\text{ini}} = \{x_1,x_2,\dots,x_l,x_{l+1},\dots,x_D \}$
        \State Initialise $\mathbf{p}^{\text{ini}} = 0$  \Comment{Force sparsity}
        \State Assemble pair $\{\mathbf{X}^{\text{ini}},\mathbf{p}^{\text{ini}}\}$
        \State $R_{f,0} = \epsilon$
        
        \While{$\lambda < \lambda_{\max}$}
            \For{$\beta = 0:\beta_{\max}$}  \Comment{Variational Annealing}
                \State $R_f = R_{f,0} \alpha^{\beta}$ 
                \State $\{\mathbf{X}^{(\beta)},\mathbf{p}^{(\beta)}\}$ = $\min_{\mathbf{X},\mathbf{p}} A_E(\mathbf{X},\mathbf{Y}) + R_f A_M(\mathbf{X},\mathbf{p},\mathbf{\hat{F}})$
                \Comment{Minimize via IPOPT}
                
                \If{$p^{(\beta)}_{k,j} < \lambda$}  \Comment{Hard-threshold $\mathbf{p}$}
                    \State $p^{(\beta)}_{k,j} = 0$
                \EndIf
            \EndFor
            
            \State model$^{(\lambda)} \leftarrow \mathbf{p}^{(\beta)}$
            \Comment{Store models}
            
            \State $\lambda = 2\lambda$ \Comment{Increase $\lambda$}
        \EndWhile
    \EndProcedure

    \end{algorithmic}
\end{algorithm}


\clearpage

\bibliography{bibliography} 
\bibliographystyle{siam}

\clearpage

\include{SupplementaryInformation/main_SI_V2}


\end{document}

%% file: SupplementaryInformation/main_SI_V2.tex

\addtocontents{toc}{\protect\setcounter{tocdepth}{3}}

\begin{center}
    {\Large Supplementary Information for \\
    
    \vskip5mm

    Model selection of chaotic systems from data with hidden variables using sparse data assimilation}
    
    \vskip10mm

    H. Ribera, S. Shirman, A. V. Nguyen and N. M. Mangan
\end{center}

\tableofcontents

\clearpage

\setcounter{figure}{0}   
\setcounter{table}{0} 
\setcounter{equation}{0} 

\newcounter{myc}[part] 
\renewcommand{\thesection}{S\arabic{myc}} 
\let\osection\section 
\renewenvironment{section}{\stepcounter{myc}\osection} 

\section{Cost function analysis}
\label{SI-costfunk-analysis}

Our aim is now to explore the landscape of the cost function as to understand the problem that we are solving and why it is very challenging. For illustrative purposes, in the following discussion we are only considering two dimensions of the cost function $\hat{A} = A_E+R_f A_M$.  We use the classical Lorenz system and take all parameters in the structure fixed and we add two extra parameters (highlighted in red),
\begin{align}
    \label{lorenzx2}
    \dot{x} &= \sigma (y - x) + {\color{red}p_{1,1}},\\
    \dot{y} &= x (\rho - z) - y + {\color{red}p_{2,1}},\\
    \label{lorenzz2}
    \dot{z} &= xy - \beta z.
\end{align}

\begin{figure}[H]
    \centering
    \includegraphics[width=0.32\textwidth]{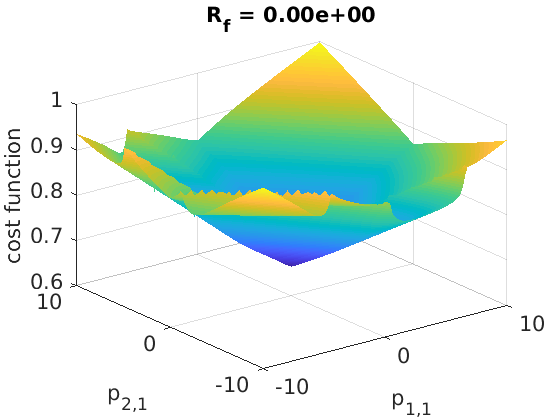}
    \includegraphics[width=0.32\textwidth]{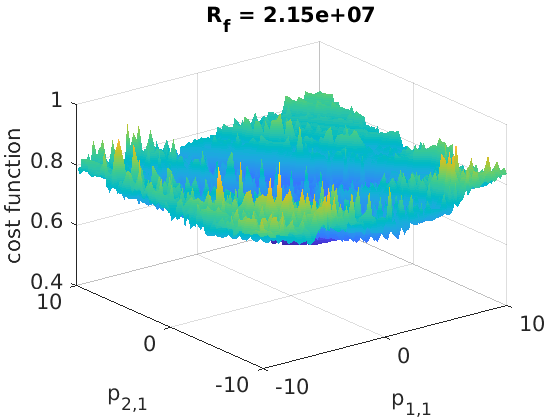}
    \includegraphics[width=0.32\textwidth]{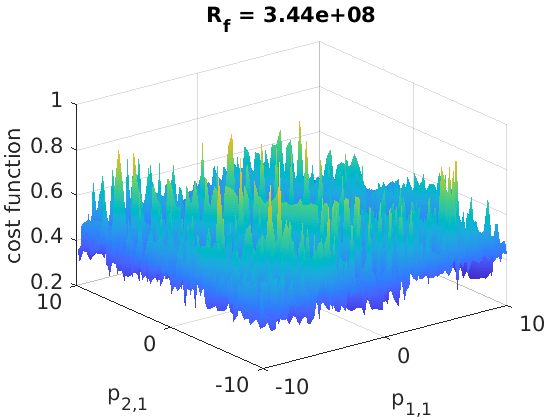}
    \caption{Varying $p_{1,1}$ and $p_{2,1}$ for three different $R_f$.}
    \label{figure_surfaces}
\end{figure}

We then vary these two parameters and plot what the cost function looks like, for three different values of $R_f$. The cost function that we want to minimise, is the one that has a large $R_f$ value (Fig. \ref{figure_surfaces}, right).

The cost function $\hat{A}$ is highly non-convex and the task of finding its global minima a priori is a difficult task.

\clearpage

\section{Time-delay embedding of used training data}
\label{SI-timeembedding}

\begin{figure}[H]
    \centering
    \includegraphics[width=0.75\textwidth]{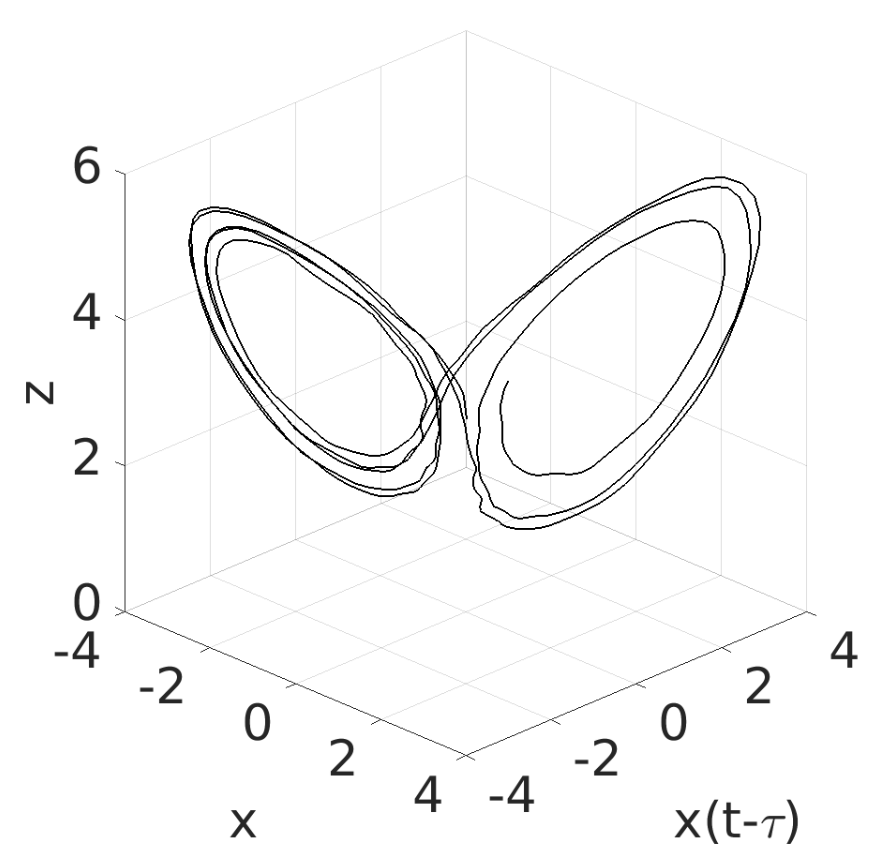}
    \caption{Time-delay embedding of training data ($\tau = 0.02$).}
    \label{figure_SI_timeembedding}
\end{figure}

\clearpage


\section{AIC calculation for synthetic data}
\label{SI-AICsynthetic}

We want to find which model is the one that best represents the data synthetic data generated (in which we added some noise $\sim \mathcal{N}(0,0.01)$). Since we are working with chaotic systems, we only expect prediction up to the Lyapunov time of the system. We consider 1/4 of the shortest Lyapunov time out of all the down-selected models for the synthetic data, $t_{M} \approx 0.3$. We use $S = 300$ time series of length $t_M$ as our validation set, but discard the first four points as they will be used to predict the initial condition for $y_0$ (as shown in the following section). To calculate the AIC score, we define the residual sum of squares of the $m$-th model as
\begin{equation}
    \text{RSS}_m = \sum_{s=1}^S E^s_{av,m}(\mathbf{Y}_s,\mathbf{F}_m, \mathbf{p}_m),
\end{equation}
where $\mathbf{Y}_s=[x_e,z_e]_s$ is the synthetic data of the time-series $s$, $\mathbf{F}_m$ the governing equations of the $m$-th model, and $\mathbf{p}_m$ denotes the parameters found via parameter estimation for the $m$-th model. $E^s_{av,m}$ is the average absolute error over the time-series $s$ and is defined as
\begin{equation}
    \label{Eavm}
    E^s_{av,m}(\mathbf{Y}_s,\mathbf{F}_m, \mathbf{p}_m) = \frac{1}{2M} \sum_{i=1}^M (x_{i,e}^s - x_i^s)^2 + (z_{i,e}^s - z_i^s)^2,
\end{equation}
where $x^s$ and $z^s$ denote the $x$ and $z$ component, respectively, of the solution of the $m$-th model in the $s$ time series, found via RK4 with $\Delta t = 0.01$. $M$ denotes 1/4 of a Lyapunov time, excluding the first four points as we have mentioned before.

Finally, we can define the AIC of the $m$-th model as 
\begin{equation}
    \label{aic_def}
    \text{AIC}_m = S \log \left( \frac{\sum_{s=1}^S E^s_{av,m}(\mathbf{Y}_s,\mathbf{F}_m, \mathbf{p}_m)}{S} \right) + 2N_{p,m}, 
\end{equation}
where $N_{p,m}$ is the number of free parameters in the $m$-th model.

We finally re-scale by the minimum AIC value, denoted by AIC$_{\text{min}}$, and so
$\Delta \text{AIC}_m = \text{AIC}_m - \text{AIC}_{\text{min}}$. 

\subsection{Initial condition choice for unmeasured y}
\label{SI-y0choice}

We need an initial condition for each time series to be able to simulate each model. We have an initial condition for both $x$ and $z$ given by the experimental data, but we do not have any information for the $y$ component. 
We cannot use the VA to estimate $y_0$ and parameters simultaneously (which would lead to better prediction windows see next section) because our validation data will then have been used for training.
Let us consider the 8th order finite difference approximation of the time derivative of $x$
\begin{equation}
    \label{discretisation}
    \frac{\ud x(t)}{\ud t} \approx \frac{
      \splitfrac{
        3x(t+4\Delta t)-32x(t+3\Delta t)+168x(t+2\Delta t)-672x(t+\Delta t)
      }{
        +672x(t+\Delta t)-168x(t+2\Delta t)+32x(t+3\Delta t)-3x(t+4\Delta t)
      }
      }{
        840 \Delta t
      }
\end{equation}

For each model, we have that 
\begin{equation}
    \label{eqforx_generic}
    \frac{\ud x(t)}{\ud t} = F_{1,m}(x(t),y(t),z(t),\mathbf{p}_m).
\end{equation}
Putting \eqref{discretisation} and \eqref{eqforx_generic} together we have
\begin{equation}
    \label{approximation_xeq}
    \frac{-x(t+2\Delta t) + 8x(t+\Delta t) - 8x(t-\Delta t) + x(t-2\Delta t)}{12 \Delta t} \approx F_{1,m}(x(t),y(t),z(t),\mathbf{p}_m).
\end{equation}
We need to solve for $y(0)$. We note that for the down-selected models in Example A in our manuscript the only terms with $y$ in the first equation in all the models is just the first order term, so for this case this is a particularly simple equation to solve.

\begin{figure}[htb!]
    \centering
    \includegraphics[width=1\textwidth]{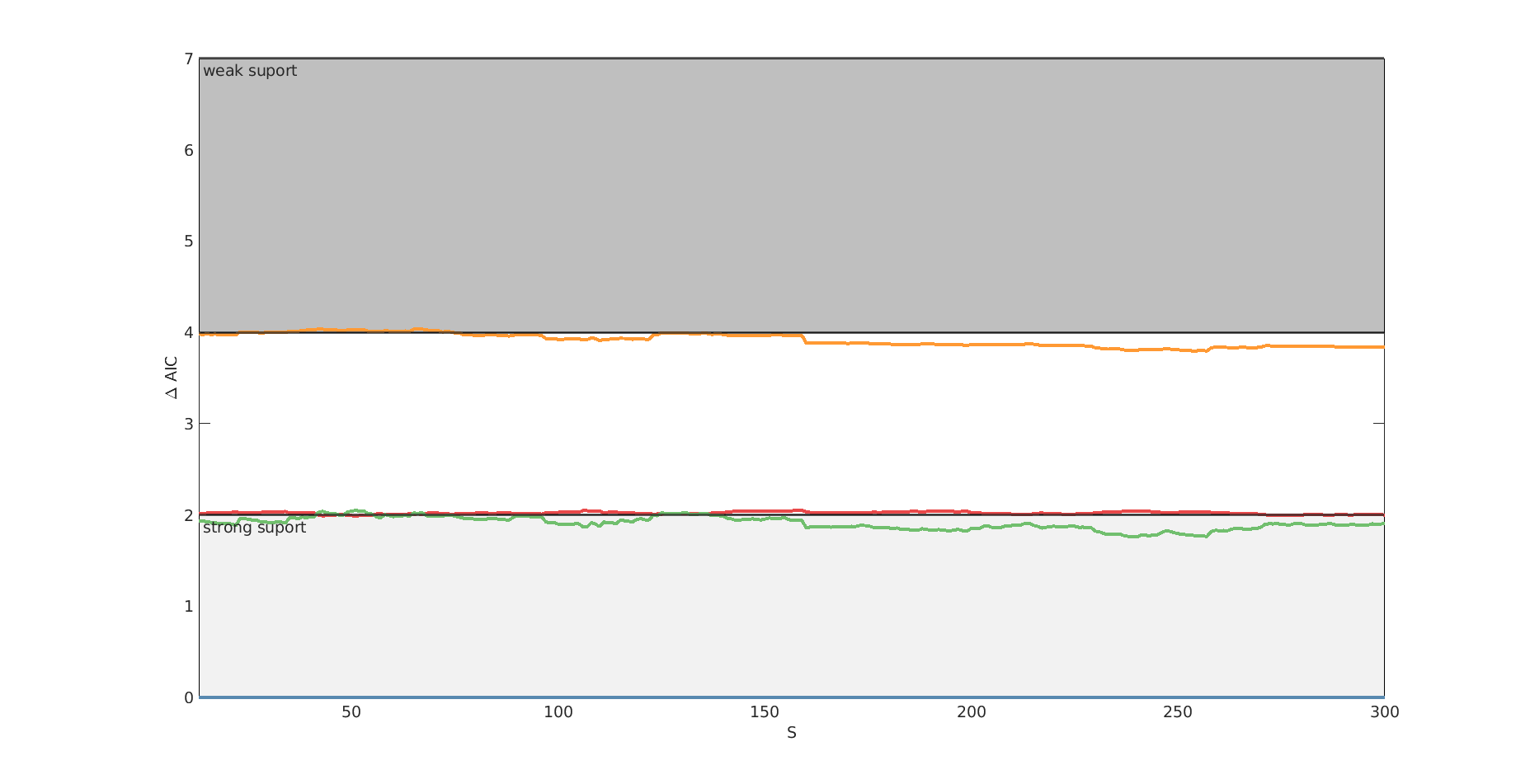}
    \caption{$\Delta$AIC from the different models DAHSI found using the synthetic data.}
    \label{figure_AIC_synthetic}
\end{figure}

The results in the synthetic data indicate that there are only four candidate models that best represent the data. Even though the $\Delta$AIC from incorrect models (Fig. \ref{figure_AIC_synthetic}, red, green and yellow lines) does not increase as we add more time series $S$ in the calculation of AIC, we consistently pick the correct model structure (blue line) as the one with lowest $\Delta AIC$.

\subsection{Prediction window}
\label{SI-predictionwindow}

\begin{figure}[htb!]
    \centering
    \includegraphics[width=\textwidth]{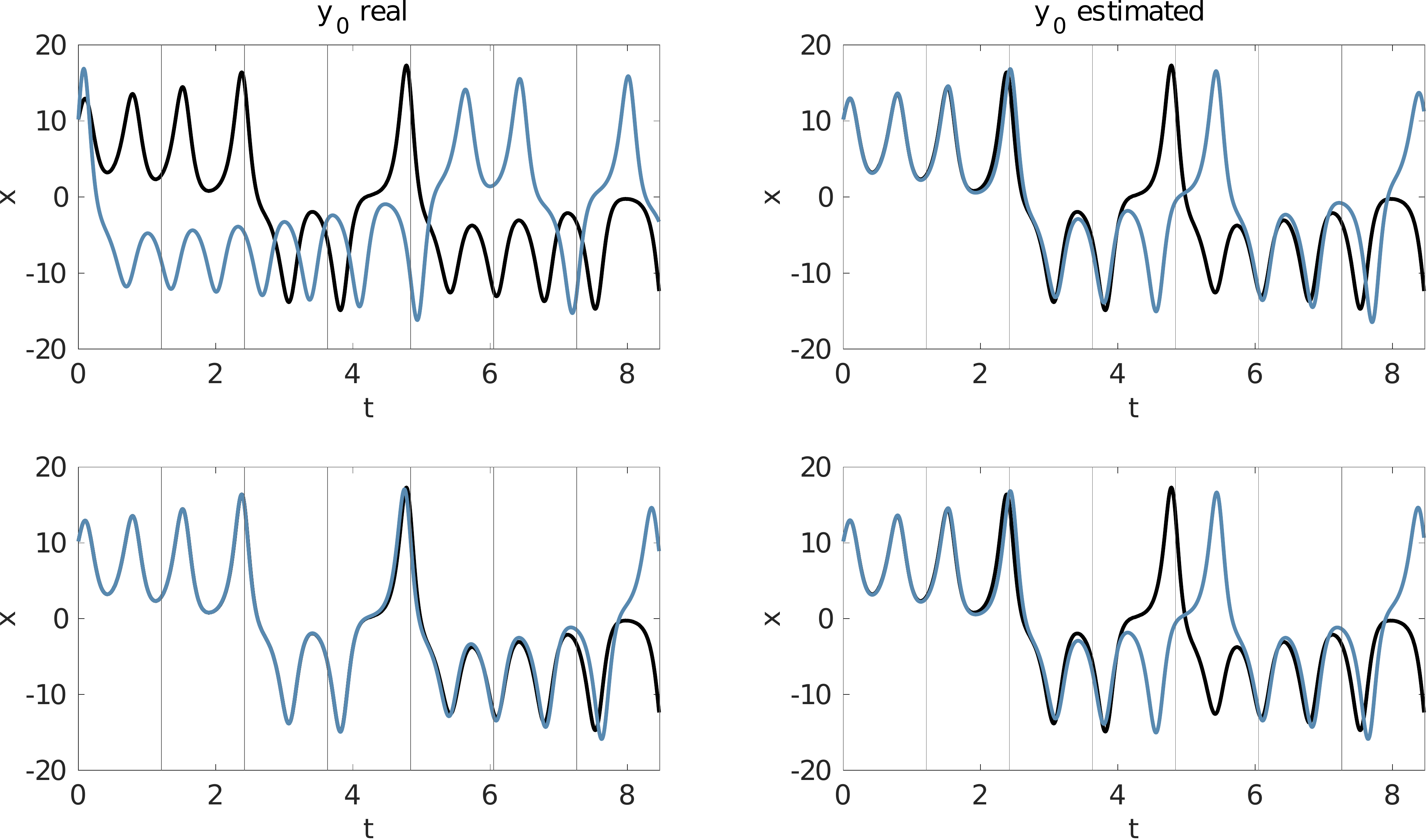}
    \caption{Top: Prediction of the model with parameters estimated using 2 observed variables ($x$ and $z$), when using the real $y_0$ and the estimated $y_0$ (as shown in \S \ref{SI-y0choice}). Bottom: Prediction of the model with parameters estimated using 3 observed variables, when using the real $y_0$ and the estimated $y_0$. The noise added to the synthetic data is $\mathcal{N}(0,\omega)$, $\omega = 0.01$.}
    \label{figure_predictionlownoise}
\end{figure}

\begin{figure}[htb!]
    \centering
    \includegraphics[width=\textwidth]{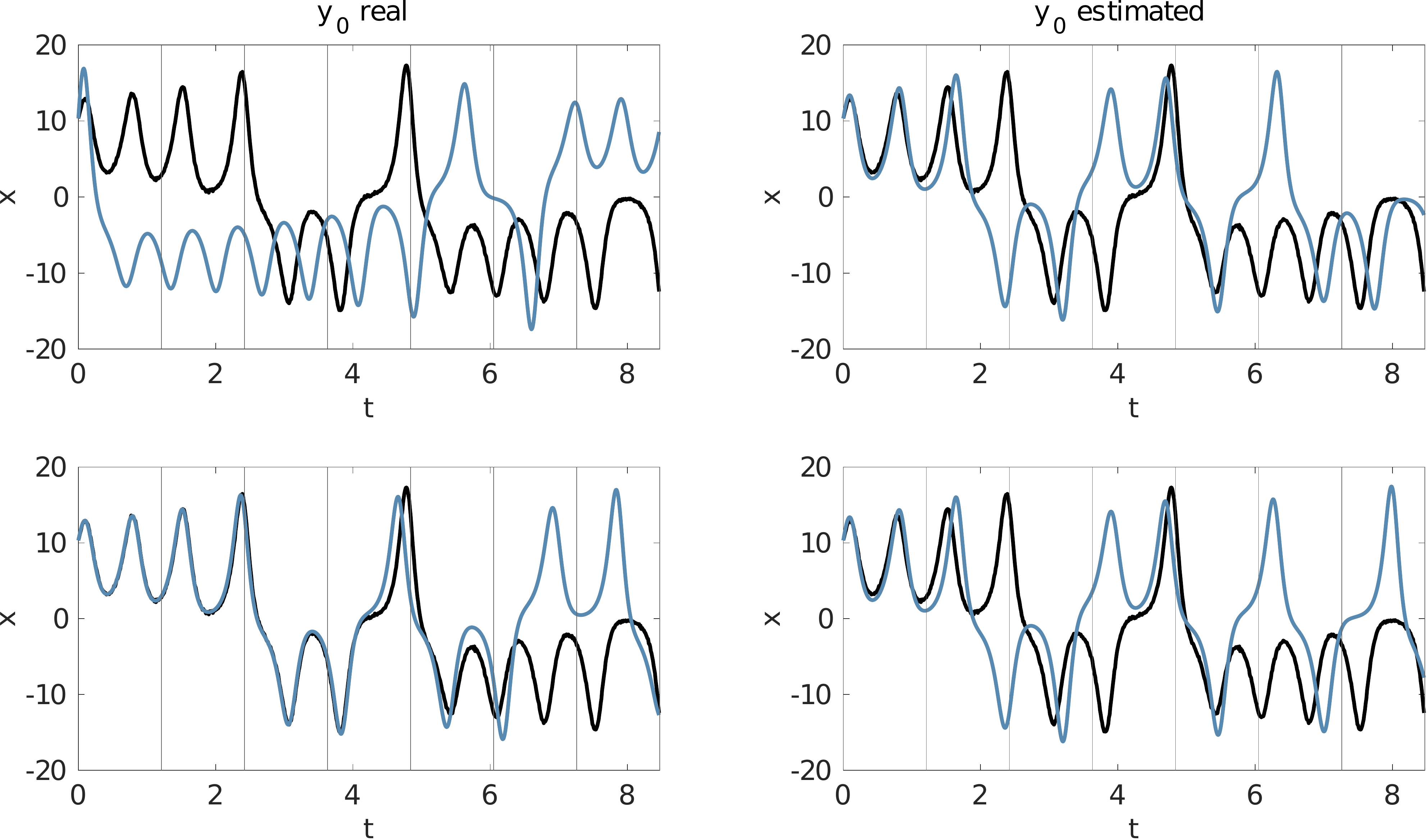}
    \caption{Top: Prediction of the model with parameters estimated using 2 observed variables ($x$ and $z$), when using the real $y_0$ and the estimated $y_0$ (as shown in \S \ref{SI-y0choice}). Bottom: Prediction of the model with parameters estimated using 3 observed variables, when using the real $y_0$ and the estimated $y_0$. The noise added to the synthetic data is $\mathcal{N}(0,\omega)$, $\omega = 0.1$.}
    \label{figure_predictionhighnoise}
\end{figure}

We compare how the prediction window changes from having two observed variables to having three observed variables. For noise $\omega = 0.01$, having one hidden variable (Fig. \ref{figure_predictionlownoise}, top row) and using the real value of $y_0$, leads to no prediction at all. However, using the estimated $y_0$ calculated as in the previous section leads to a prediction window of about 3.5 Lyapunov times. This shows that the parameter estimates and $y_0$ estimate are compensating for each other. For the case of all variables observed (Fig. \ref{figure_predictionlownoise}, bottom row), we see that using the real value of $y_0$ leads to a prediction window of about 6 Lyapunov times. If we estimate $y_0$ the prediction window reduces to about 3.5 Lyapunov times.
For a higher noise $\omega = 0.1$, having one hidden variable (Fig. \ref{figure_predictionhighnoise}, top row) and using the real value of $y_0$, again leads to no prediction at all. Moreover, using the estimated $y_0$ calculated as in the previous section leads to a shorter prediction window than for lower noise, about 1 Lyapunov time. This shows that with increased noise amplified by hidden variables, the $y_0$ estimate cannot compensate for the parameter estimate that well. For the case of all variables observed with $\omega = 0.1$ (Fig. \ref{figure_predictionhighnoise}, bottom row), we see that using the real value of $y_0$ leads to a prediction window of about 3.5 Lyapunov times. If we estimate $y_0$ the prediction window reduces to about 1 Lyapunov time.

\clearpage

\section{Down-selected models}
\label{SI-downselected}

We present the structure of the 25 down-selected models, but we do not provide the parameter estimation (the parameter values can be found in the code package \cite{Ribera2021GitHub}).

\begin{align}
	\frac{\ud x}{\ud t} &= p_{1,2}x+p_{1,3}y,\\
	\frac{\ud y}{\ud t} &= p_{2,2}x+p_{2,7}xz,\\
	\frac{\ud z}{\ud t} &= p_{3,4}z+p_{3,6}xy.
\end{align}

\begin{align}
	\frac{\ud x}{\ud t} &= p_{1,2}x+p_{1,3}y,\\
	\frac{\ud y}{\ud t} &= p_{2,2}x+p_{2,3}y+p_{2,7}xz,\\
	\frac{\ud z}{\ud t} &= p_{3,4}z+p_{3,6}xy.
\end{align}

\begin{align}
	\frac{\ud x}{\ud t} &= p_{1,2}x+p_{1,3}y+p_{1,4}z,\\
	\frac{\ud y}{\ud t} &= p_{2,2}x+p_{2,3}y+p_{2,7}xz,\\
	\frac{\ud z}{\ud t} &= p_{3,4}z+p_{3,6}xy.
\end{align}

\begin{align}
	\frac{\ud x}{\ud t} &= p_{1,1}+p_{1,2}x+p_{1,3}y,\\
	\frac{\ud y}{\ud t} &= p_{2,2}x+p_{2,7}xz,\\
	\frac{\ud z}{\ud t} &= p_{3,4}z+p_{3,6}xy.
\end{align}

\begin{align}
	\frac{\ud x}{\ud t} &= p_{1,1}+p_{1,2}x+p_{1,3}y+p_{1,4}z,\\
	\frac{\ud y}{\ud t} &= p_{2,2}x+p_{2,7}xz,\\
	\frac{\ud z}{\ud t} &= p_{3,4}z+p_{3,6}xy.
\end{align}

\begin{align}
	\frac{\ud x}{\ud t} &= p_{1,1}+p_{1,2}x+p_{1,3}y+p_{1,4}z,\\
	\frac{\ud y}{\ud t} &= p_{2,1}+p_{2,2}x+p_{2,3}y+p_{2,7}xz,\\
	\frac{\ud z}{\ud t} &= p_{3,4}z+p_{3,6}xy.
\end{align}

\begin{align}
	\frac{\ud x}{\ud t} &= p_{1,2}x+p_{1,3}y,\\
	\frac{\ud y}{\ud t} &= p_{2,2}x+p_{2,7}xz,\\
	\frac{\ud z}{\ud t} &= 0.
\end{align}

\begin{align}
	\frac{\ud x}{\ud t} &= p_{1,2}x+p_{1,3}y,\\
	\frac{\ud y}{\ud t} &= p_{2,2}x+p_{2,7}xz,\\
	\frac{\ud z}{\ud t} &= p_{3,6}xy.
\end{align}

\begin{align}
	\frac{\ud x}{\ud t} &= p_{1,2}x+p_{1,3}y,\\
	\frac{\ud y}{\ud t} &= p_{2,2}x+p_{2,3}y+p_{2,7}xz,\\
	\frac{\ud z}{\ud t} &= 0.
\end{align}

\begin{align}
	\frac{\ud x}{\ud t} &= p_{1,1}+p_{1,2}x+p_{1,3}y,\\
	\frac{\ud y}{\ud t} &= p_{2,2}x+p_{2,7}xz+,\\
	\frac{\ud z}{\ud t} &= 0.
\end{align}

\begin{align}
	\frac{\ud x}{\ud t} &= p_{1,2}x+p_{1,3}y,\\
	\frac{\ud y}{\ud t} &= p_{2,2}x+p_{2,3}y+p_{2,7}xz,\\
	\frac{\ud z}{\ud t} &= p_{3,6}xy.
\end{align}

\begin{align}
	\frac{\ud x}{\ud t} &= p_{1,2}x+p_{1,3}y+p_{1,4}z,\\
	\frac{\ud y}{\ud t} &= p_{2,2}x+p_{2,3}y+p_{2,7}xz,\\
	\frac{\ud z}{\ud t} &= 0.
\end{align}

\begin{align}
	\frac{\ud x}{\ud t} &= p_{1,1}+p_{1,2}x+p_{1,3}y,\\
	\frac{\ud y}{\ud t} &= p_{2,2}x+p_{2,7}xz,\\
	\frac{\ud z}{\ud t} &= p_{3,6}xy.
\end{align}

\begin{align}
	\frac{\ud x}{\ud t} &= p_{1,2}x+p_{1,3}y+p_{1,4}z,\\
	\frac{\ud y}{\ud t} &= p_{2,2}x+p_{2,7}xz,\\
	\frac{\ud z}{\ud t} &= p_{3,4}z+p_{3,6}xy.
\end{align}

\begin{align}
	\frac{\ud x}{\ud t} &= p_{1,2}x+p_{1,3}y+p_{1,4}z,\\
	\frac{\ud y}{\ud t} &= p_{2,2}x+p_{2,3}y+p_{2,7}xz,\\
	\frac{\ud z}{\ud t} &= p_{3,6}xy.
\end{align}

\begin{align}
	\frac{\ud x}{\ud t} &= p_{1,1}+p_{1,2}x+p_{1,3}y+p_{1,4}z,\\
	\frac{\ud y}{\ud t} &= p_{2,2}x+p_{2,7}xz,\\
	\frac{\ud z}{\ud t} &= p_{3,6}xy.
\end{align}

\begin{align}
	\frac{\ud x}{\ud t} &= p_{1,2}x+p_{1,3}y+p_{1,4}z,\\
	\frac{\ud y}{\ud t} &= p_{2,2}x+p_{2,7}xz,\\
	\frac{\ud z}{\ud t} &= p_{3,4}z+p_{3,6}xy+p_{3,7}xz.
\end{align}

\begin{align}
	\frac{\ud x}{\ud t} &= p_{1,1}+p_{1,2}x+p_{1,3}y+p_{1,4}z,\\
	\frac{\ud y}{\ud t} &= p_{2,2}x+p_{2,3}y+p_{2,7}xz,\\
	\frac{\ud z}{\ud t} &= p_{3,6}xy.
\end{align}

\begin{align}
	\frac{\ud x}{\ud t} &= p_{1,2}x+p_{1,3}y+p_{1,4}z,\\
	\frac{\ud y}{\ud t} &= p_{2,2}x+p_{2,6}xy+p_{2,7}xz,\\
	\frac{\ud z}{\ud t} &= p_{3,4}z+p_{3,6}xy+p_{3,7}xz.
\end{align}

\begin{align}
	\frac{\ud x}{\ud t} &= p_{1,1}+p_{1,2}x+p_{1,3}y+p_{1,4}z,\\
	\frac{\ud y}{\ud t} &= p_{2,2}x+p_{2,7}xz,\\
	\frac{\ud z}{\ud t} &= p_{3,4}z+p_{3,6}xy+p_{3,7}xz.
\end{align}

\begin{align}
	\frac{\ud x}{\ud t} &= p_{1,1}+p_{1,2}x+p_{1,3}y+p_{1,4}z,\\
	\frac{\ud y}{\ud t} &= p_{2,1}+p_{2,2}x+p_{2,3}y+p_{2,7}xz,\\
	\frac{\ud z}{\ud t} &= p_{3,6}xy.
\end{align}

\begin{align}
	\frac{\ud x}{\ud t} &= p_{1,2}x+p_{1,3}y+p_{1,4}z,\\
	\frac{\ud y}{\ud t} &= p_{2,2}x+p_{2,3}y+p_{2,6}xy+p_{2,7}xz,\\
	\frac{\ud z}{\ud t} &= p_{3,4}z+p_{3,6}xy+p_{3,7}xz.
\end{align}

\begin{align}
	\frac{\ud x}{\ud t} &= p_{1,1}+p_{1,2}x+p_{1,3}y+p_{1,4}z,\\
	\frac{\ud y}{\ud t} &= p_{2,2}x+p_{2,6}xy+p_{2,7}xz,\\
	\frac{\ud z}{\ud t} &= p_{3,4}z+p_{3,6}xy+p_{3,7}xz.
\end{align}

\begin{align}
	\frac{\ud x}{\ud t} &= p_{1,1}+p_{1,2}x+p_{1,3}y+p_{1,4}z,\\
	\frac{\ud y}{\ud t} &= p_{2,2}x+p_{2,3}y+p_{2,6}xy+p_{2,7}xz,\\
	\frac{\ud z}{\ud t} &= p_{3,4}z+p_{3,6}xy+p_{3,7}xz.
\end{align}

\begin{align}
	\frac{\ud x}{\ud t} &= p_{1,1}+p_{1,2}x+p_{1,3}y+p_{1,4}z,\\
	\frac{\ud y}{\ud t} &= p_{2,1}+p_{2,2}x+p_{2,3}y+p_{2,6}xy+p_{2,7}xz,\\
	\frac{\ud z}{\ud t} &= p_{3,4}z+p_{3,6}xy+p_{3,7}xz.
\end{align}

\subsection{Models identified in the Pareto front edge}
\label{SI-Paretomodels}

\begin{table}[H]
\centering
\caption{Models identified in the Pareto front in Figure \ref{figure_Scheme}(d) in the main text.}
\label{table_params_comparison_pareto}
    \resizebox{0.78\textwidth}{!}{%
    \begin{tabular}{l|r|rrrrrrr}
     & & \multicolumn{7}{c}{number of active terms} \\
     & Term & 6 & 7 & 8 & 9 & 10 & 11 & 12 \\ \hline
    \multirow{10}{*}{eq. $\dot{x}$} & $1$ & 0 & 0 & -0.8112 & 0 & 0 & -0.2514 & -2.4053 \\
    & $x$ & -16.5556 & -16.9554 & -16.4666 & -16.5603 & -17.0172 & -17.0582 & -17.0627 \\  
    & $y$ & 19.8000 & 18.7853 & 19.8120 & 16.7514 & 19.9884 & 19.9840 & 19.9862 \\  
    & $z$ & 0 & 0 & 0.0276 & 0.1486 & 0.1596 & 0.1833 & 1.4595 \\  
    & $x^2$ & 0 & 0 & 0 & 0 & 0 & 0 & 0 \\  
    & $xy$ & 0 & 0 & 0 & 0 & 0 & 0 & 0 \\  
    & $xz$ & 0 & 0 & 0 & 0 & 0 & 0 & 0 \\  
    & $y^2$ & 0 & 0 & 0 & 0 & 0 & 0 & 0 \\  
    & $yz$ & 0 & 0 & 0 & 0 & 0 & 0 & 0 \\  
    & $z^2$ & 0 & 0 & 0 & 0 & 0 & 0 & 0 \\  \hline
    \multirow{10}{*}{eq. $\dot{y}$} & $1$ & 0 & 0 & 0 & 0 & 0 & 0 & 0.7892 \\
    & $x$ & 23.2613 & 24.3535 & 23.0763 & 27.3789 & 22.6028 & 22.6017 & 22.6061 \\  
    & $y$ & 0 & 0.2580 & 0 & 0 & 0.3346 & 0.3567 & 0.3298 \\  
    & $z$ & 0 & 0 & 0 & 0 & 0 & 0 & 0 \\  
    & $x^2$ & 0 & 0 & 0 & 0 & 0 & 0 & 0 \\  
    & $xy$ & 0 & 0 & 0 & -0.0922 & -0.0906 & -0.0843 & -0.3647 \\  
    & $xz$ & -6.3345 & -6.7054 & -6.2868 & -7.4621 & -6.2507 & -6.2561 & -6.2691 \\  
    & $y^2$ & 0 & 0 & 0 & 0 & 0 & 0 & 0 \\  
    & $yz$ & 0 & 0 & 0 & 0 & 0 & 0 & 0 \\  
    & $z^2$ & 0 & 0 & 0 & 0 & 0 & 0 & 0 \\  \hline
    \multirow{10}{*}{eq. $\dot{z}$} & $1$ & 0 & 0 & 0 & 0 & 0 & 0 & 0 \\
    & $x$ & 0 & 0 & 0 & 0 & 0 & 0 & 0 \\  
    & $y$ & 0 & 0 & 0 & 0 & 0 & 0 & 0 \\  
    & $z$ & -3.6646 & -3.6835 & -3.6736 & 4.3951 & -3.6954 & -3.6966 &  -3.6941 \\  
    & $x^2$ & 0 & 0 & 0 & 0 & 0 & 0 & 0 \\  
    & $xy$ & 5.1948 & 4.8273 & 5.2315 & -3.6660 & 5.1412 & 5.1292 & 5.1326  \\  
    & $xz$ & 0 & 0 & 0 & 0.0883 & 0.0903 & 0.0791 & 0.2900 \\  
    & $y^2$ & 0 & 0 & 0 & 0 & 0 & 0 & 0 \\  
    & $yz$ & 0 & 0 & 0 & 0 & 0 & 0 & 0 \\  
    & $z^2$ & 0 & 0 & 0 & 0 & 0 & 0 & 0 \\  \hline    
    $E_{av}$ &  & 10.1693 & 9.7441 & 9.7174 & 9.6778 & 9.0995 & 9.0345 & 9.5765 \\ 
    \end{tabular}
    }      
\end{table}   

\subsection{AIC and BIC on the 25 down-selected models}
\label{SI-aic-bic-experimental}
Bayesian information criteria (BIC) is defined as 
\begin{equation}
    \label{bic_def}
    \text{BIC}_m = S \log \left( \frac{\sum_{s=1}^S E^s_{av,m}(\mathbf{Y}_s,\mathbf{F}_m, \mathbf{p}_m)}{S} \right) + S \log(N_{p,m}). 
\end{equation}
In the same way when we defined $\Delta$AIC in a previous section, we re-scale by the minimum BIC value, denoted by BIC$_{\text{min}}$, and so
$\Delta \text{BIC}_m = \text{BIC}_m - \text{BIC}_{\text{min}}$.

We will now calculate how AIC (\eqref{aic_def}) and BIC (\eqref{bic_def}) change as we add more time series into the calculation. For each $S$ that we use to calculate both AIC and BIC ($S \leq 1083$, which is the total number of time segments we have available that are of length 1/4 of a Lyapunov time), we will pick $S$ random time-segments to ensure that the $S$ time-segments used in the calculation are independent samples.

\begin{figure}[htb!]
    \centering
    \includegraphics[width=1\textwidth]{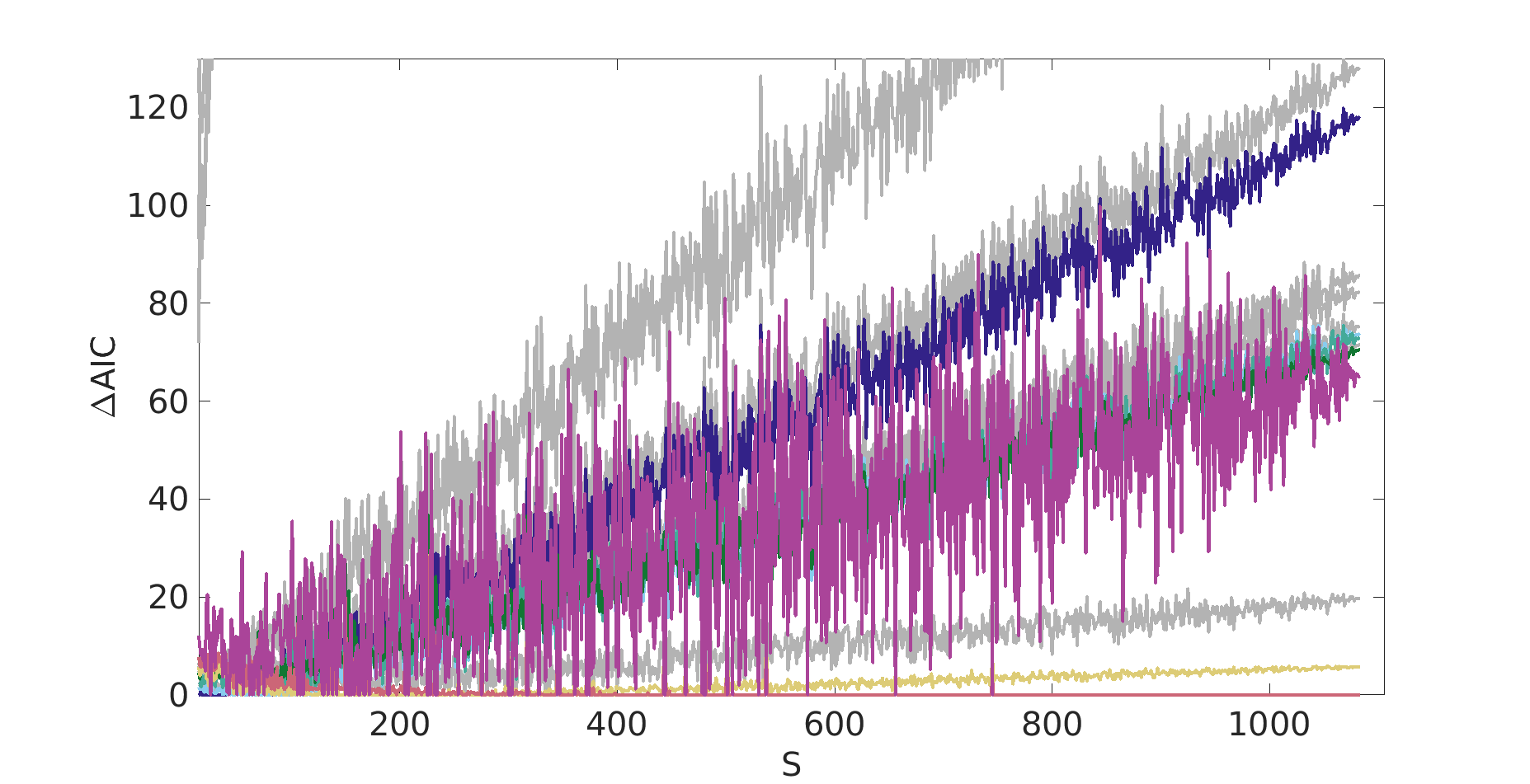}
    \includegraphics[width=1\textwidth]{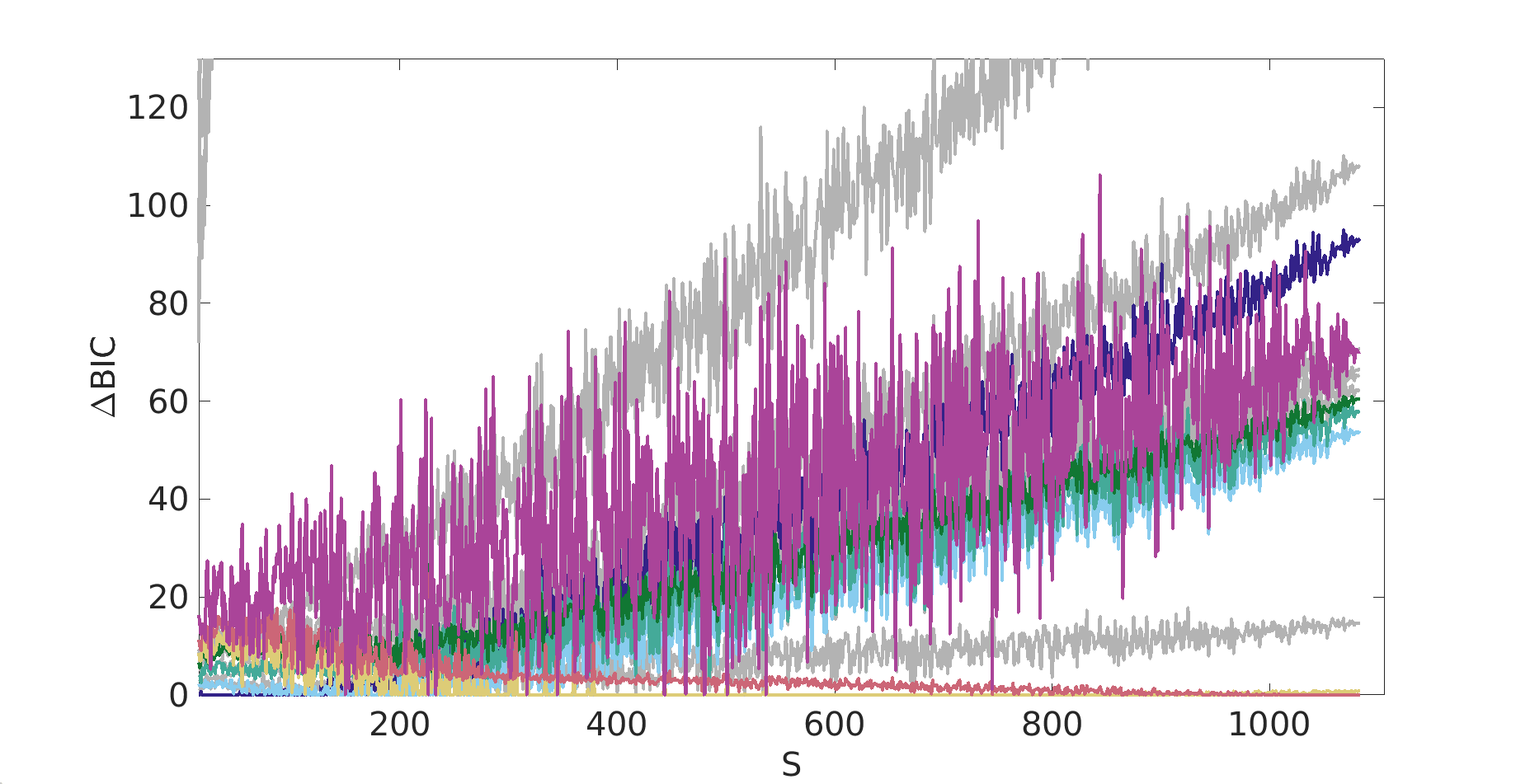}
    \caption{$\Delta \text{AIC}_m$ and $\Delta \text{BIC}_m$ from the different models DAHSI found using the experimental data in \cite{Blakely2007}.}
    \label{figure_AIC_BIC_experimental}
\end{figure}

For both $\Delta \text{AIC}_m$ and $\Delta \text{BIC}_m$ we are able to consistently identify a unique model (Fig. \ref{figure_AIC_BIC_experimental}). If we just look at the Pareto front (Fig. 1(d) in the main text), one might ask if the decrease between 9 and 10 terms is meaningful. Both AIC and BIC say that it is.

\clearpage


\section{Action derivation}
\label{SI-action-derivation}

We consider a dynamical system with unknown governing equations
\begin{equation}
\label{dynamics-2}
    \frac{\ud \mathbf{X}}{\ud t} = \mathbf{F}(\mathbf{X}(t),\mathbf{p}),
\end{equation}
where $\mathbf{X} = (x_1,x_2,\dots,x_D) \in \mathbb{R}^D$ are the state variables, $\mathbf{F} = (F_1,\,F_2,\dots,F_D)$ are the unknown functions that govern the dynamics of the system and $\mathbf{p}$ is a set of unknown parameters. Te measurements $\mathbf{Y} = (y_1,y_2,\dots,y_L) \in \mathbb{R}^L$ are lower dimensional $L \leq D$ than the underlying variables.

Our goal is to find $\mathbf{X}$ and $\mathbf{p}$ that maximise the probability $P(\mathbf{X}, \mathbf{p} \; | \; \mathbf{Y}, \mathbf{\hat{F}})$. We have that \cite{Abarbanel2013}
\begin{equation}
    P(\mathbf{X}, \mathbf{p} \; | \; \mathbf{Y}, \mathbf{\hat{F}}) = \int \exp\left[ -A_0(\mathbf{X},\mathbf{Y}) \right] \; \ud \mathbf{X}.
\end{equation}
Furthermore,
\begin{equation}
    \label{action}
    A_0(\mathbf{X},\mathbf{Y}) = - \sum_{i=1}^N \text{CMI}\left[\mathbf{X}(t_i),\mathbf{Y}(t_i) \; | \; \mathbf{Y}(t_0),\dots,\mathbf{Y}(t_{i-1})\right] - \sum_{i=1}^{N-1} \log\left[ P(\mathbf{X}(t_{i+1}),\mathbf{p} \; | \; \mathbf{X}(t_{i}), \mathbf{\hat{F}}) \right],
\end{equation}

We make the following assumptions:
\begin{enumerate}
    \item The measurements $\mathbf{Y}$ have uncorrelated Gaussian error and that there is no correlation between errors in measuring different quantities or at varying time points \cite{Abarbanel2013};
    \item The state at the next time point depends only on the state at the current time point, and that our model can have some error by widening the $\delta$ function it would follow otherwise using a Gaussian approximation of it \cite{Abarbanel2013};
    \item Each element in $\mathbf{p}$ follows a Laplace distribution with \emph{mean} $0$ and diversity $b$.
\end{enumerate}

With assumption 1 it can be shown that
\begin{equation}
    \text{CMI}\left[\mathbf{X}(t_i),\mathbf{Y}(t_i) \; | \; \mathbf{Y}(t_0),\dots,\mathbf{Y}(t_{i-1})\right] = \frac{1}{2 \sigma_m^2} \sum_{l=1}^L \left(x_l(t_i) - y_l(t_i)\right)^2.
\end{equation}

For the second term in the sum, we need to find an expression for $P(\mathbf{X}(t_{i+1}),\mathbf{p} \; | \; \mathbf{X}(t_{i}), \mathbf{\hat{F}})$. Let us now focus on the $k$-th component of $\mathbf{X}(t_{i+1})$, and so our goal is to find an expression for $P(x_k(t_{i+1}), \mathbf{p}_k \; | \; \mathbf{X}(t_{i}), F_k)$. 

We consider the library of $q$ possible functions and the generic expression for each equation of our model:
\begin{equation}
    \label{generic-model}
    \hat{F}_k \coloneqq \hat{F}_k(\mathbf{X},\mathbf{p}) = p_{k,1} \theta_1(\mathbf{X}) + p_{k,2} \theta_2(\mathbf{X}) + \cdots + p_{k,q} \theta_q(\mathbf{X}),
\end{equation}
for $k=1,2,\dots,D$.

We can rewrite the probability we are seeking as
\begin{equation}
    \label{main_eq_action}
    P(x_k(t_{i+1}),\mathbf{p}_k \; | \; \mathbf{X}(t_i), F_k) = P(\mathbf{p}_k \; | \; x_k(t_{i+1}), \mathbf{X}(t_i), F_k) P(x_k(t_{i+1}) \; | \; \mathbf{X}(t_i), F_k).
\end{equation}
Now each term in the right hand side can also be rewritten as
\begin{align}
    P(\mathbf{p}_k \; | \; x_k(t_{i+1}), \mathbf{X}(t_i), F_k) &= \frac{P(x_k(t_{i+1}),\mathbf{X}(t_i), F_k \; | \; \mathbf{p}_k) P(\mathbf{p}_k)}{P(x_k(t_{i+1}),\mathbf{X}(t_i), F_k)}, \\
    P(x_k(t_{i+1}) \; | \; \mathbf{X}(t_i), F_k) &= \frac{P(x_k(t_{i+1}),\mathbf{X}(t_i), F_k)}{P(\mathbf{X}(t_i),F_k)}.    
\end{align}

Thus, \eqref{main_eq_action} becomes
\begin{equation}
    \label{main_eq_action_2}
    P(x_k(t_{i+1}),\mathbf{p}_k \; | \; \mathbf{X}(t_i), F_k) = \frac{P(x_k(t_{i+1}),\mathbf{X}(t_i), F_k \; | \; \mathbf{p}_k) P(\mathbf{p}_k)}{P(\mathbf{X}(t_i),F_k)}.
\end{equation}

We can rewrite the first therm on the right hand side in \eqref{main_eq_action_2} as a likelihood,
\begin{equation}
    P(x_k(t_{i+1}),\mathbf{X}(t_i), F_k \; | \; \mathbf{p}_k) = \mathcal{L}(\mathbf{p}_k \; | \; x_k(t_{i+1}),\mathbf{X}(t_i), F_k).
\end{equation}
Assuming that our next state follows a normal distribution with mean $f_k$ and standard deviation $\sigma^2$,
\begin{equation}
    \mathcal{L}(\mathbf{p}_k \; | \; x_k(t_{i+1}),\mathbf{X}(t_i), F_k) = \frac{1}{\sigma \sqrt{2\pi}} \exp{
    \left(
    -\frac{
    \left[ 
    x_k(t_{i+1}) - f_k(\mathbf{X},\mathbf{p},F_k) \right]^2
    }{2\sigma^2}
    \right)
    }.
\end{equation}
With assumption 3, we know that each $p_{k,j}$ follows a Laplace distribution,
\begin{equation}
    p_{k,j} \sim \text{Laplace}(0,b) = \frac{1}{2b}\exp{\left(-\frac{|p_{k,j}|}{b}\right)},
\end{equation}
and so
\begin{equation}
    P(\mathbf{p}_k) = \prod_{j=1}^q \frac{1}{2b}\exp{\left(-\frac{|p_{k,j}|}{b}\right)}.
\end{equation}

With this we can write \eqref{main_eq_action_2} as
\begin{equation}
    \begin{split}
    \label{main_eq_action_3}
        P(&x_k(t_{i+1}),\mathbf{p}_k \; | \; \mathbf{X}(t_i), F_k) \propto \\
        &\propto \frac{1}{\sigma \sqrt{2\pi}} \exp{\left(-\frac{\left[ x_k(t_{i+1}) - f_k(\mathbf{X},\mathbf{p},F_k) \right]^2}{2\sigma^2}\right)} \prod_{j=1}^q \frac{1}{2b}\exp{\left(-\frac{|p_{k,j}|}{b}\right)}.        
    \end{split}
\end{equation}
Note that since we are going to be minimising the action $A_0$ (\eqref{action}) we forget about the constant term $P(\mathbf{X}(t_i), F_k)$ in the denominator and we just have a proportionality instead of an equality.

Note that because the $k$-th current state only depends upon the previous one,
\begin{equation}
    P\left(\mathbf{X}(t_{i+1}), \mathbf{p} \; | \; \mathbf{X}(t_{i}), \mathbf{\hat{F}}\right) = \prod_{k=1}^D P(x_k(t_{i+1}), \mathbf{p}_k \; | \; \mathbf{X}(t_i), F_k),
\end{equation}
and so, finally, we can write
\begin{equation}
    \begin{split}
    \label{main_eq_action_4}
        P(&\mathbf{X}(t_{i+1}),\mathbf{p} \; | \; \mathbf{X}(t_{i}), \mathbf{\hat{F}}) \propto \\ &\propto \prod_{k=1}^D \left\{ \frac{1}{\sigma \sqrt{2\pi}} \exp{\left(-\frac{\left[ x_k(t_{i+1}) - f_k(\mathbf{X},\mathbf{p},F_k) \right]^2}{2\sigma^2}\right)} \prod_{j=1}^q \frac{1}{2b}\exp{\left(-\frac{|p_{k,j}|}{b}\right)}\right\}.
    \end{split}
\end{equation}

Upon taking the logarithm to this expression above,
\begin{equation}
    \label{main_eq_action_5}
        \log(P(\mathbf{X}(t_{i+1}),\mathbf{p} \; | \; \mathbf{X}(t_{i}), \mathbf{\hat{F}})) \propto \sum_{k=1}^D \left\{ -\frac{\left[ x_k(t_{i+1}) - f_k(\mathbf{X},\mathbf{p},F_k) \right]^2}{2\sigma^2} - \lambda \Vert \mathbf{p}_k \Vert_1 \right\} + \frac{D}{\sigma \sqrt{2\pi}} + \frac{\lambda D}{2},
\end{equation}
where $\lambda = q/b$. 

We have seen that \eqref{action} becomes
\begin{equation}
    \label{costfunk-appendix}
    A(\mathbf{X},\mathbf{p}) = \frac{1}{N}\sum_{i=1}^N \Vert \mathbf{X}(t_{i}) - \mathbf{Y}(t_{i}) \Vert^2 + \frac{1}{N}\sum_{i=1}^{N-1} R_f \left\{ \Vert \mathbf{X}(t_{i+1}) - \mathbf{f}(\mathbf{X}(t_n),\mathbf{p},\mathbf{\hat{F}}) \Vert^2 \right\} + \lambda \Vert \mathbf{p} \Vert_1,
\end{equation}
which is what we wanted to show.

\clearpage


\section{Computational time}
\label{SI-computationaltime}

\begin{figure}[htb!]
    \centering
    \includegraphics[width=0.49\textwidth]{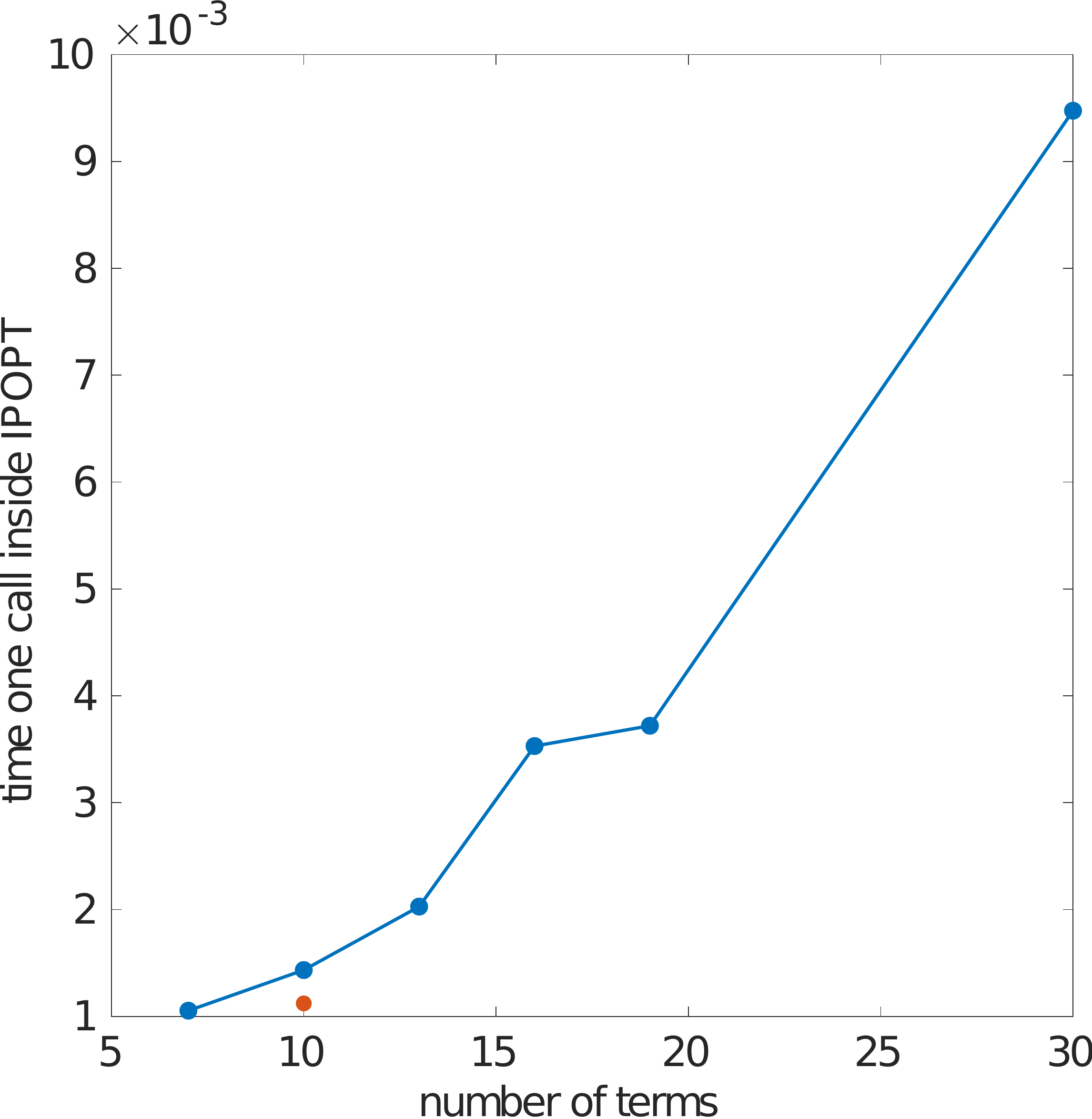}
    \includegraphics[width=0.49\textwidth]{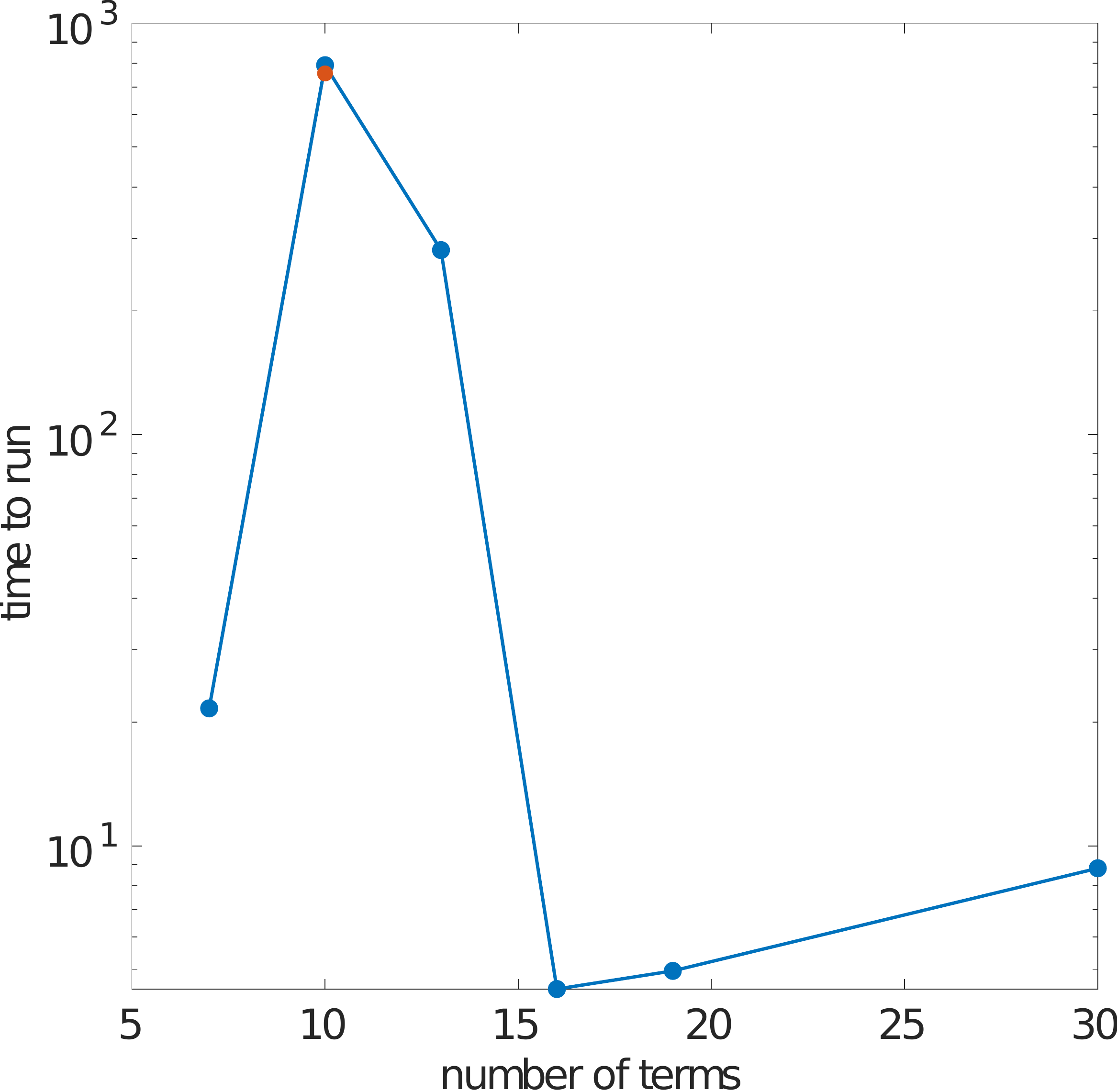}
    \caption{7 terms: parameter estimation. 10 terms: in blue $x^2$ in each equation; in red $1$ in each equation. 13 terms: $x^2$ and $y^2$ in each equation. 16 terms: $x^2$, $y^2$ and $z^2$ in each equation. 19 terms: $1$, $x^2$, $y^2$ and $z^2$ in each equation. 30 terms: model selection.}
    \label{figure_computational}
\end{figure}

We use the Lorenz system, with all variables observed, $N=1001$ time points, $\Delta t = 0.01$, and no noise. The more terms our library $\mathbfgl{\Theta}$, the more time it takes to evaluate the cost function associated, its Jacobian and its Hessian (Fig. \ref{figure_computational}(left)).
However, due to model symmetries and other structural features, the time to run our algorithm does not monotonically increase with increasing number of terms in our library. A library with 10 terms can take 100 times more to run than the full library of 30 monomials (Fig. \ref{figure_computational}(right)).

\clearpage


\section{Semiconductor}
\label{SI-semiconductor}

We consider this semiconductor model ($T$ trap levels with two possible states differing by one electronic unit
of charge),
\begin{align}
    \label{eq-n}
    \frac{\ud x}{\ud t} &= e_{n,01} y - R_{n,10} xz,\\
    \frac{\ud y}{\ud t} &= -e_{n,01} y + R_{n,10} xz,\\
    \label{eq-nd1}
    \frac{\ud z}{\ud t} &= e_{n,01} y - R_{n,10} xz. 
\end{align}
$x$ denotes the number of electrons in the conduction band, $y$ denotes the number of traps with 2 electrons, and $z$ denotes the number of traps with 1 electron. We chose $e_{n,01} = 0.5$ and $R_{n,10} = 0.25$.

\begin{figure}[H]
    \centering
    \includegraphics[width=.45\textwidth]{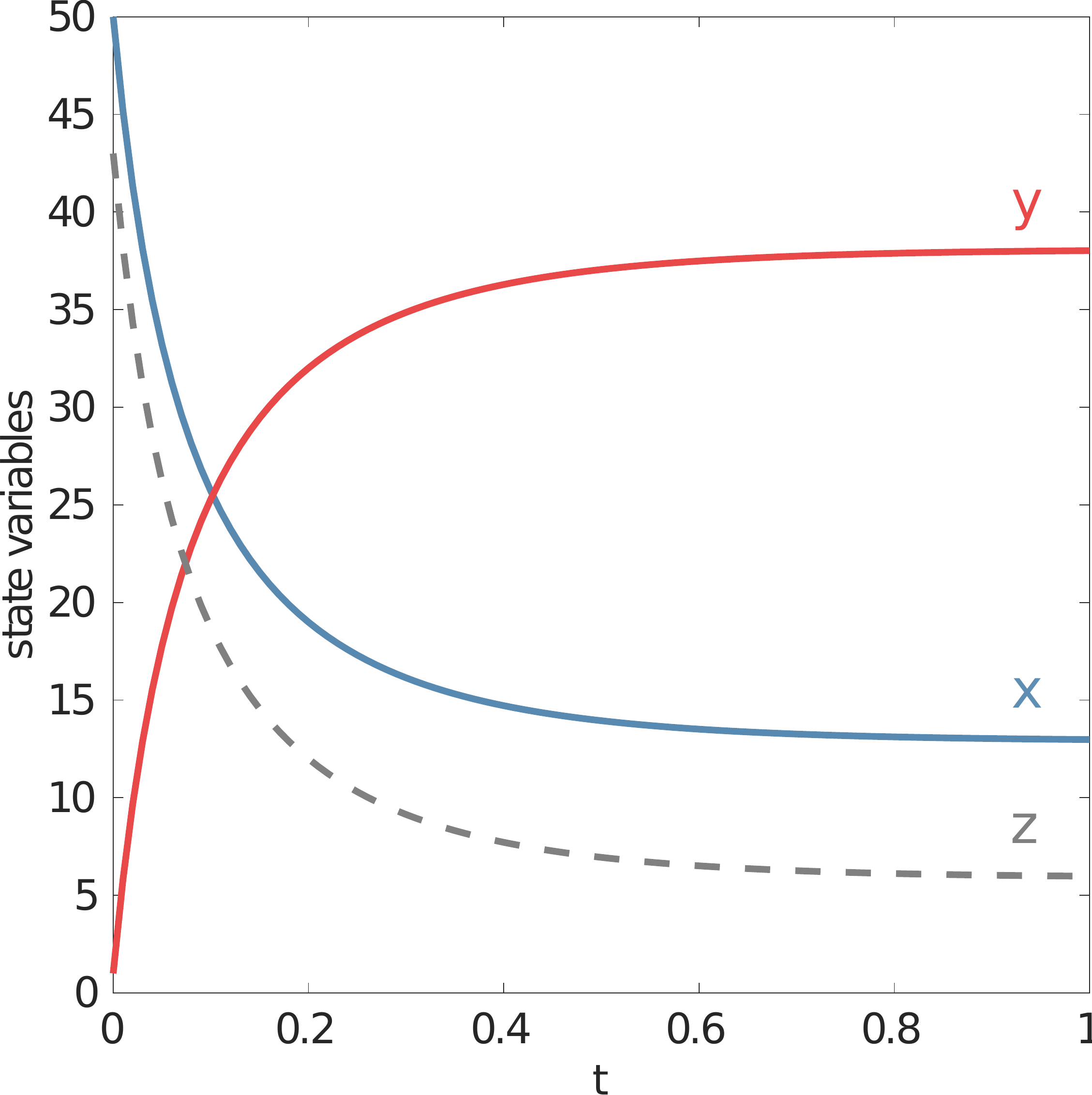}
    \caption{Dynamics from the original system \eqref{eq-n}-\eqref{eq-nd1}.}
    \label{figure_SemiCond_timeseries}
\end{figure}

Instead of using the library of all monomials in three variables up to degree two, we know that there are only a few terms make sense physically. Our generic model for this example is
\begin{align}
    \frac{\ud x}{\ud t} &= p_{1,1} + p_{1,2} x + p_{1,3} y + p_{1,4} z + p_{1,5} x^2 + p_{1,7} xz,\\
    \frac{\ud y}{\ud t} &= p_{2,1} + p_{2,2} x + p_{2,3} y + p_{2,4} z + p_{2,5} x^2 + p_{2,7} xz,\\
    \frac{\ud z}{\ud t} &= p_{3,1} + p_{3,2} x + p_{3,3} y + p_{3,4} z + p_{3,7} xz.   
\end{align}

\begin{figure}[H]
    \centering
    \includegraphics[width=0.45\textwidth]{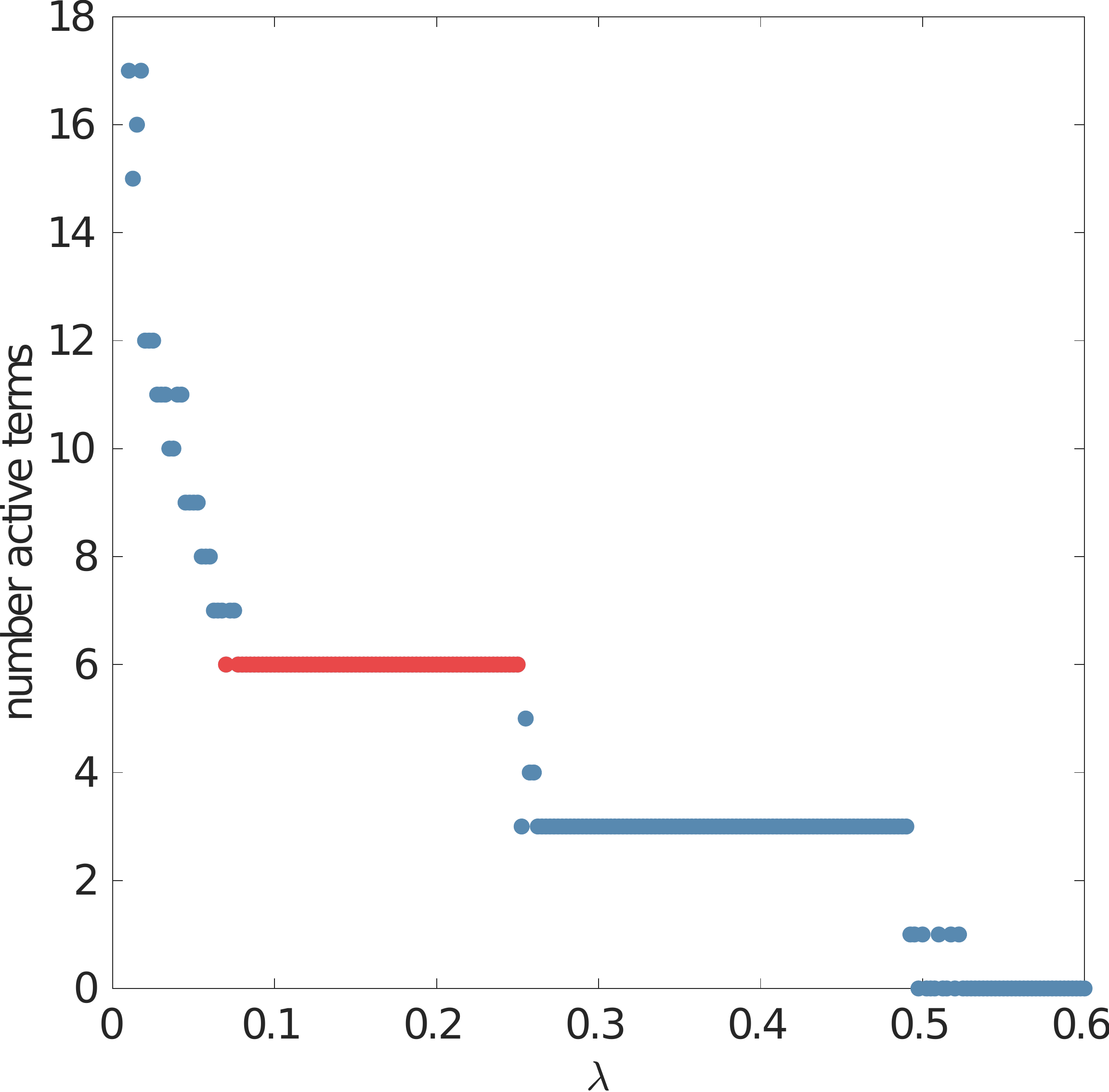}
    \hspace{10mm}
    \includegraphics[width=0.45\textwidth]{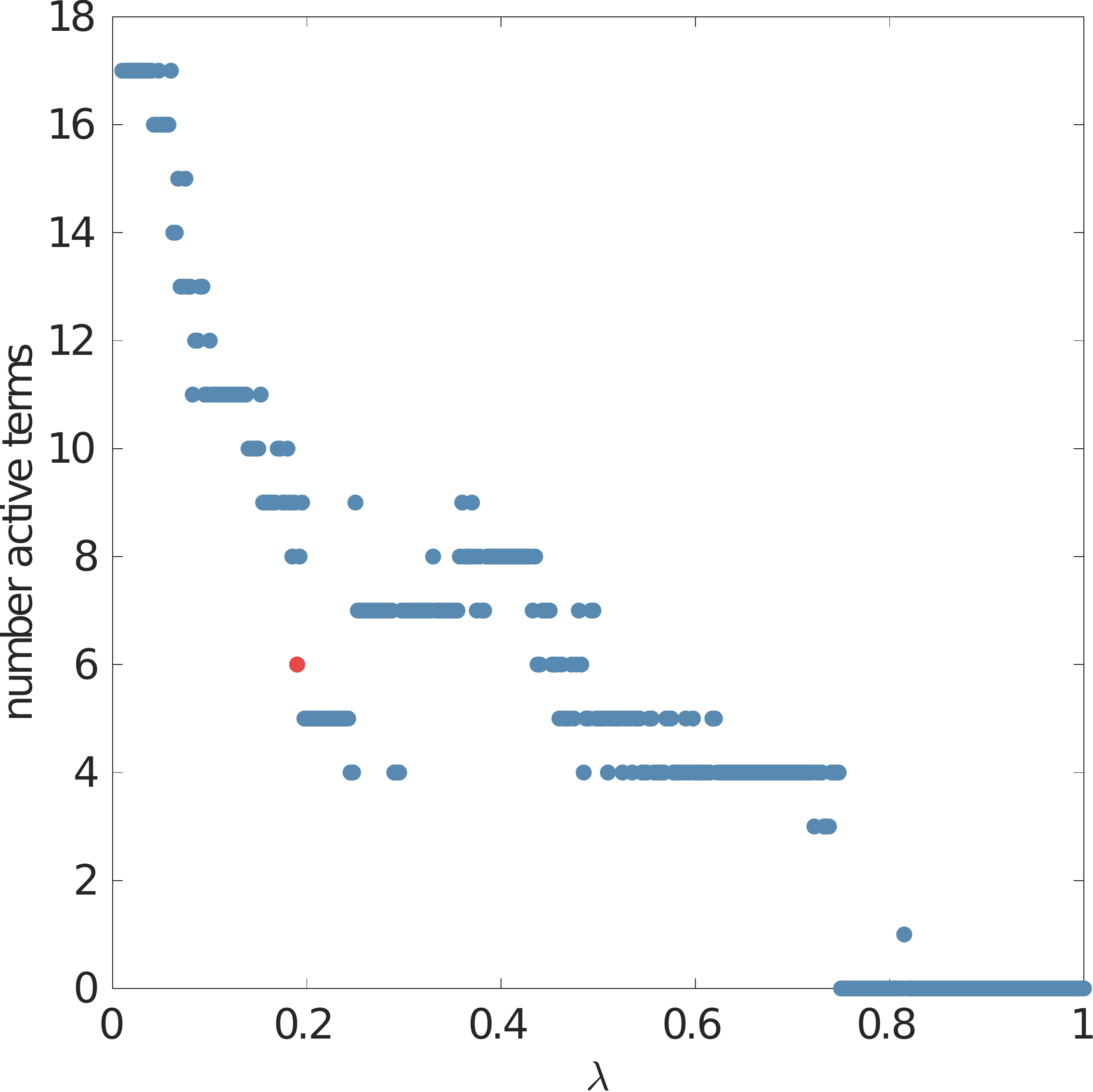}
    \caption{Left: all observed variables. Right: one hidden variable. Highlighted in red are the $\lambda$ that lead to model recovery.}
    \label{figure_lambdasweep_semicond}
\end{figure}

We first consider three observed variables, $D=L=3$. 
We consider a time series of $N = 101$ equally spaced time points, with $\Delta t = 0.01$. The $\lambda$ sweep results in a different amount of active terms for each value. See Fig. \ref{figure_lambdasweep_semicond} (left). Since we know the model from which our data comes from, we just want to see if the model that has the right number of terms (highlighted in red) corresponds to our original one, which it does.

\subsection{1 hidden variable}

\begin{figure}[htb!]
    \centering
    \includegraphics[width=0.45\textwidth]{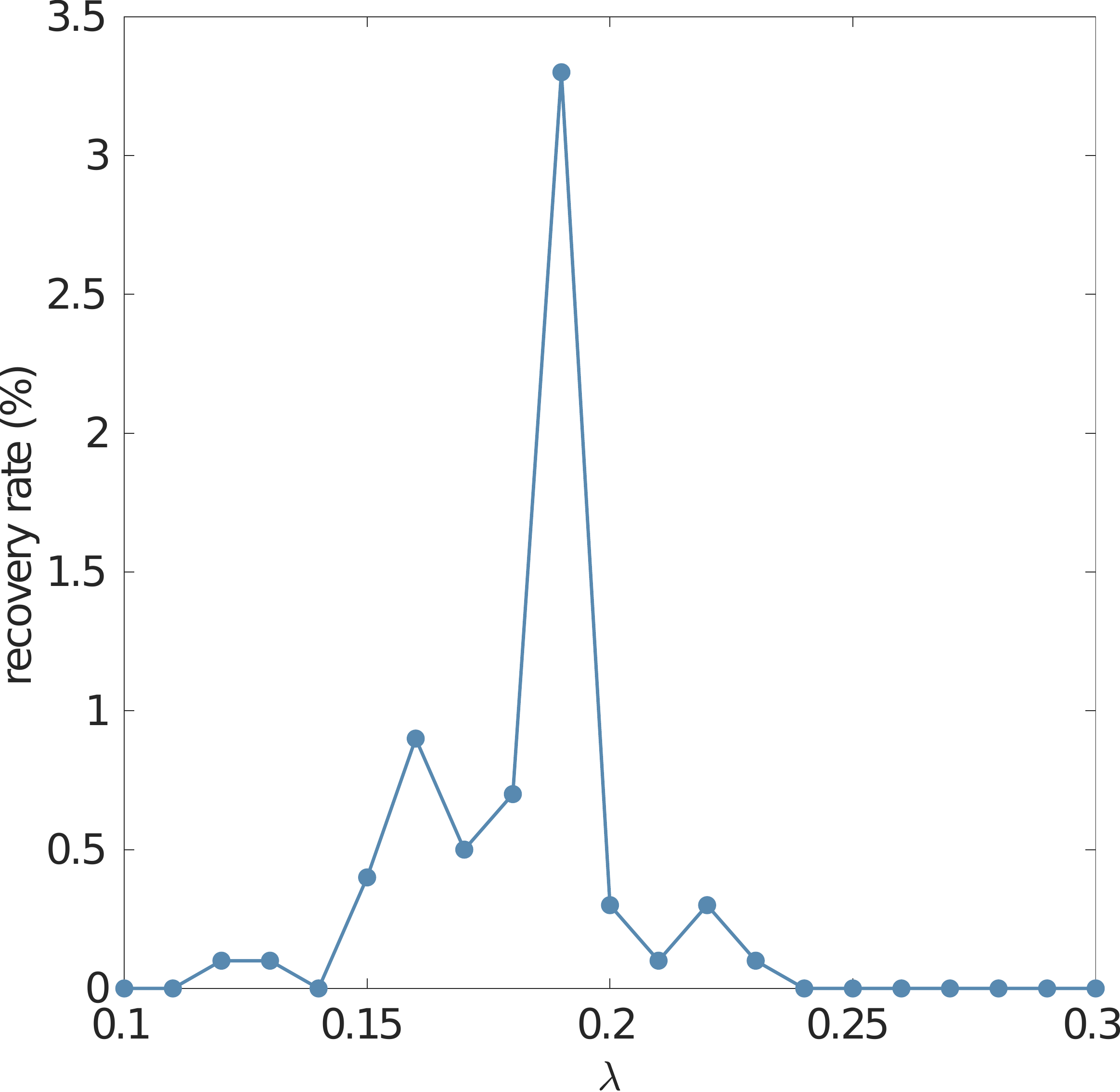}
    \caption{Percentage of recovery rate for different $\lambda$ values for 1,000 different initialisations. Initial guess for unmeasured variables is obtained through the derivatives of the measured variables.}
    \label{figure_lambdasweep_Grad}
\end{figure}

We consider two observed variables, $L=2$. We pick $x$ and $y$. We run $N_I = 1,000$ different initialisations. There is a question in this particular case on how the initial guess should be picked (see Algorithm \ref{Algorithm_IniGuess}.
We do a $\lambda$ sweep from $\lambda = 0.1$ through $\lambda = 0.3$. Out of all the 1,000 different initialisations, we recover the right sparsity pattern 68 times. The optimal $\lambda = 0.19$, for which we recover the right sparsity pattern 33 times (see Fig. \ref{figure_lambdasweep_Grad}).

\begin{table}[H]
    \centering
    \begin{tabular}{ccccccc}
    \hline
    observed & hidden & N & $\Delta$ t & $\beta_{\max}$ & $\lambda$ & recovery \\ \specialrule{1.3pt}{1pt}{1pt}
    2 & 1 ($z$) & 101 & 0.01 & 30 & 0.19 & 3.3\% \\ \hline
    \end{tabular}
    \caption{Recovery of the semiconductor system with one hidden variable.}
\end{table}

\begin{algorithm}[H]
\caption{Algorithm for picking an initial guess for unobserved variables in the semiconductor case}
\label{Algorithm_IniGuess}    
    \begin{algorithmic}[1]
    
        \For{$d = 1:(D-L)$}   \Comment{loop in unmeasured variables}
        
            \State Pick at random one of the observed variables.
            
            \State $dX_d \leftarrow$ Calculate gradient vector from its time series.
            
            \While{$Z_d$ out of bounds} \Comment{Make sure unmeasured variable is within bounds}  
                \State $Z_d(t_1) \leftarrow$ Random initial condition for unobserved variable within bounds.
                
                \For{$i = 1:N-1$} 
                    \State $Z_d(t_{i+1}) = \Delta t \times dX_d(t_i) + Z_d(t_i)$ 
                \EndFor
            \EndWhile
        \EndFor            

    \end{algorithmic}
\end{algorithm}

\subsection{Parameter identifiability}
\label{SI-parameter-identif}

There are two main reasons of why a parameter might not be identifiable: said parameter does not influence the model output; there is a interdependence among different parameters, that is, one can compensate the change of one parameter (that would influence the model output) by changing other parameter(s) and have the output be the same. In this section, we focus on the latter.

One way to detect pairwise interplay is by plotting contours of the cost function versus pairs of parameters. Largely eccentric contours or \emph{valleys} show that the cost function is almost unchanged in one direction, and the two parameters are highly correlated. 
The main drawback for our case in particular is that we will be limited to find relationships only between pairs of parameters instead of higher dimensional interactions. 

\begin{figure}
    \centering
    \includegraphics[height=0.195\textheight]{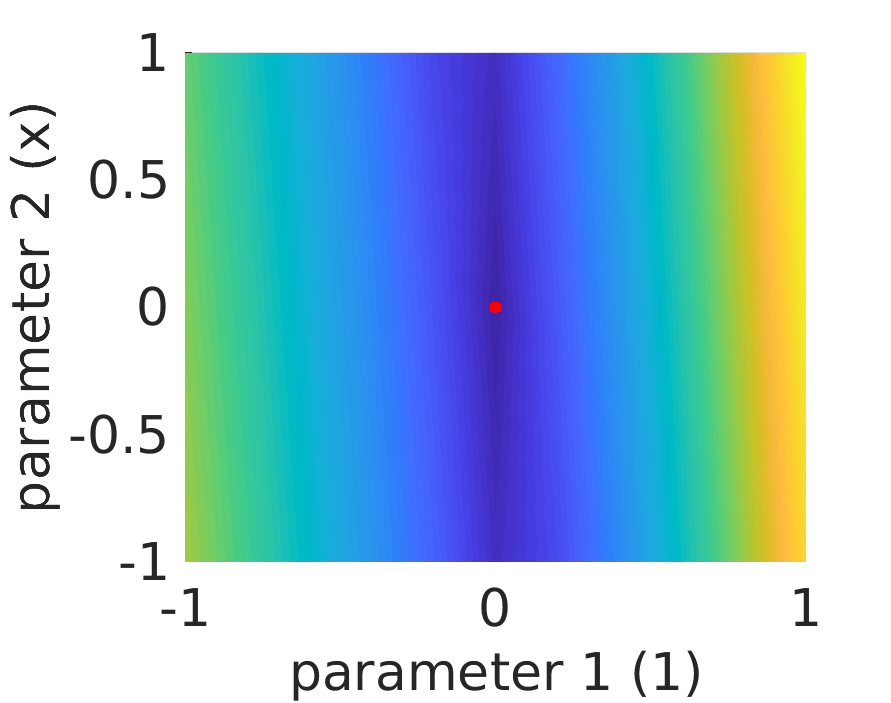}
    \includegraphics[height=0.195\textheight]{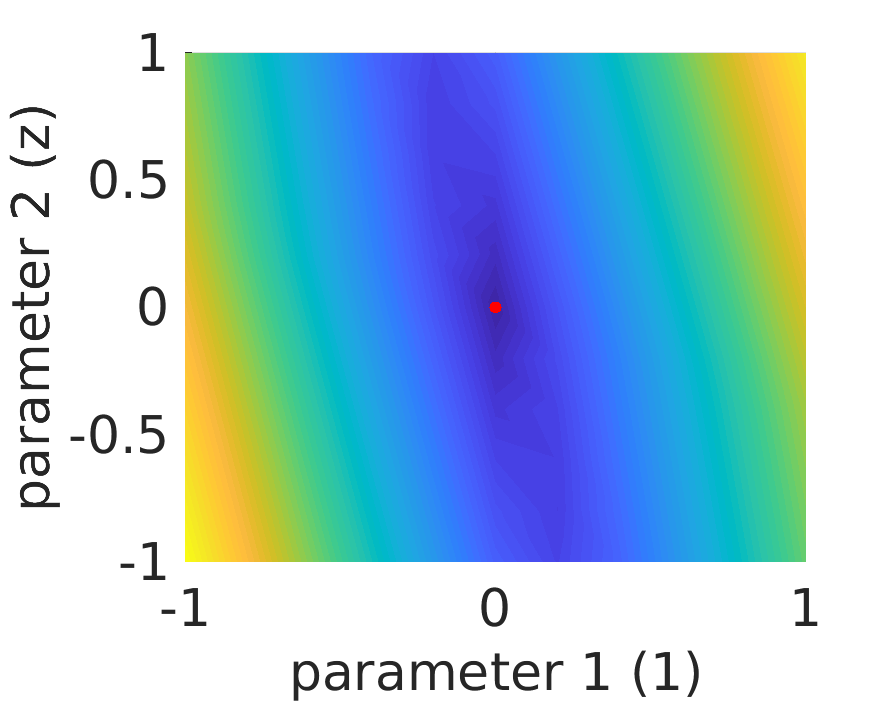}\\
    \includegraphics[height=0.195\textheight]{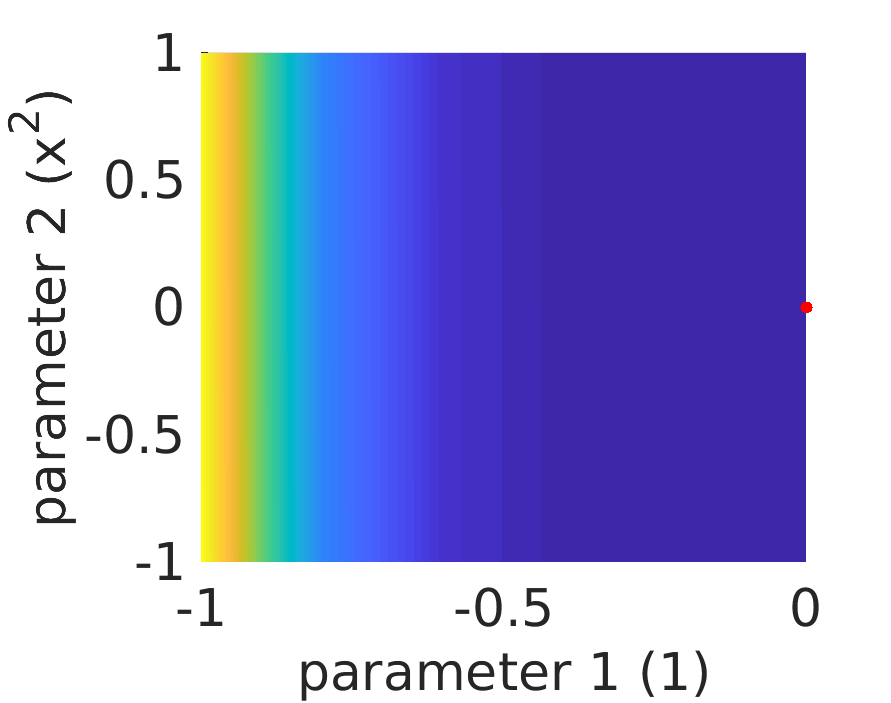}
    \includegraphics[height=0.195\textheight]{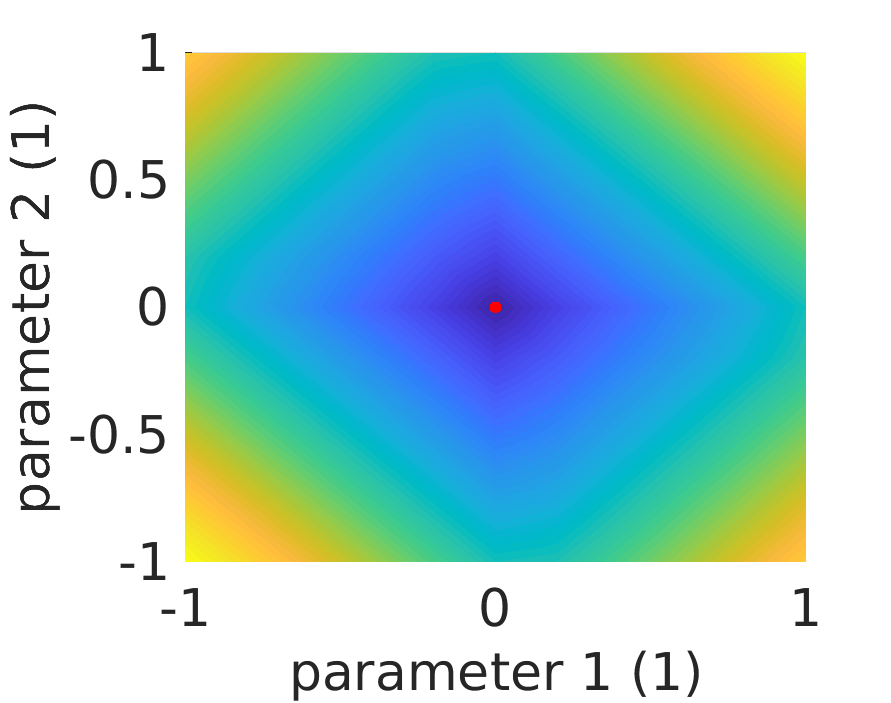}\\
    \includegraphics[height=0.195\textheight]{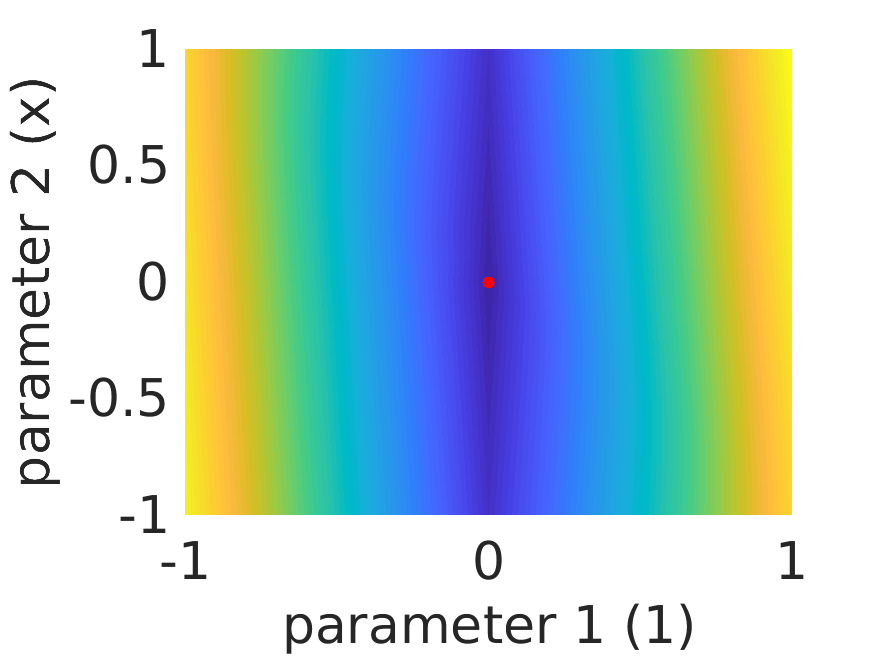}
    \includegraphics[height=0.195\textheight]{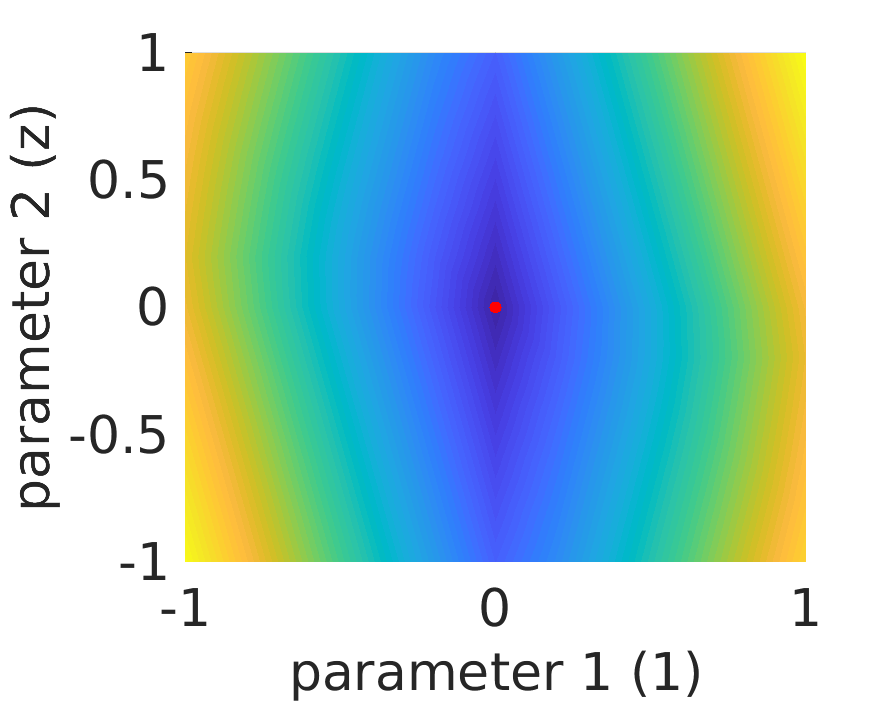}\\
    \includegraphics[height=0.195\textheight]{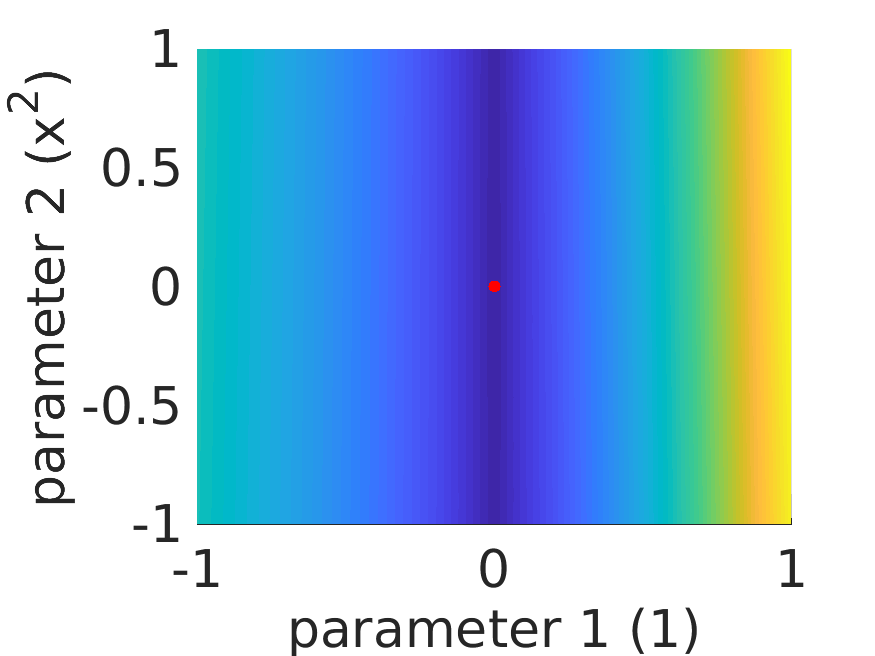}
    \includegraphics[height=0.195\textheight]{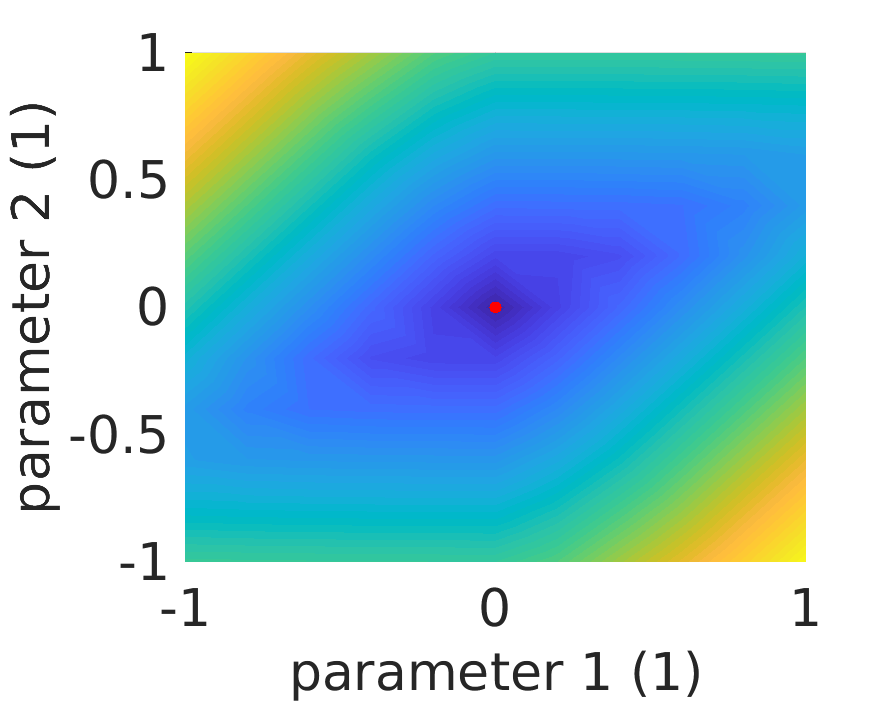}\\
    \includegraphics[height=0.195\textheight]{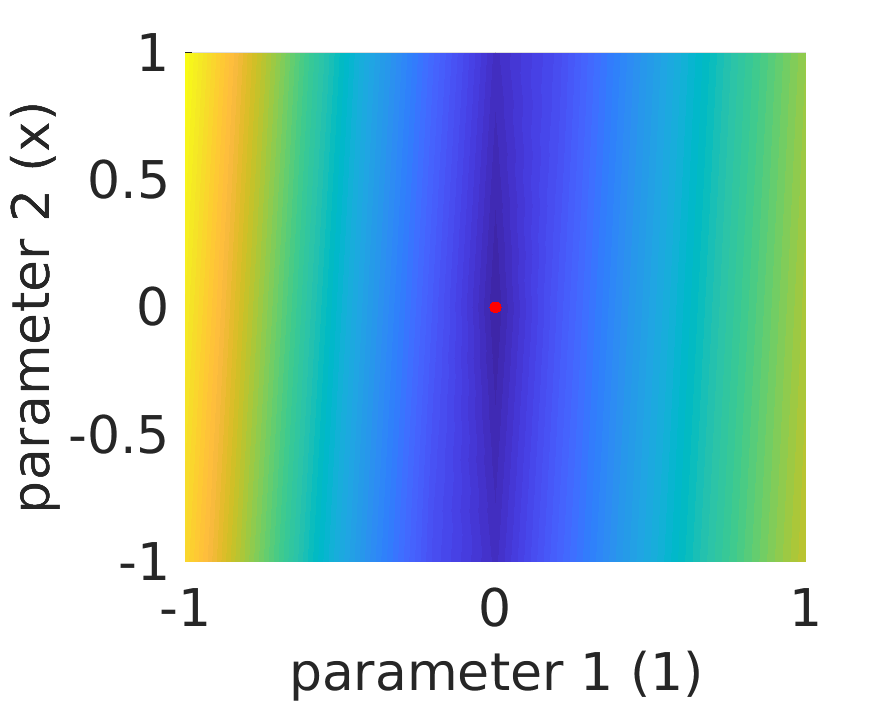}
    \includegraphics[height=0.195\textheight]{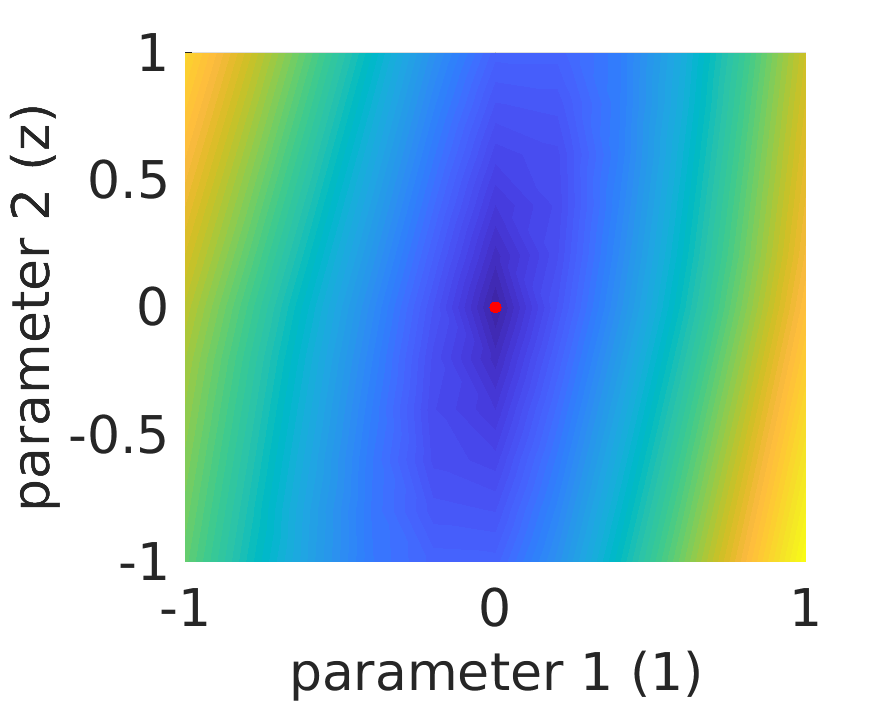}
    \caption{parameter 1 is $1$ on the $\ud x/\ud t$ equation.}
    \label{figure_identi_1}
\end{figure}

\begin{figure}
    \centering
    \includegraphics[height=0.195\textheight]{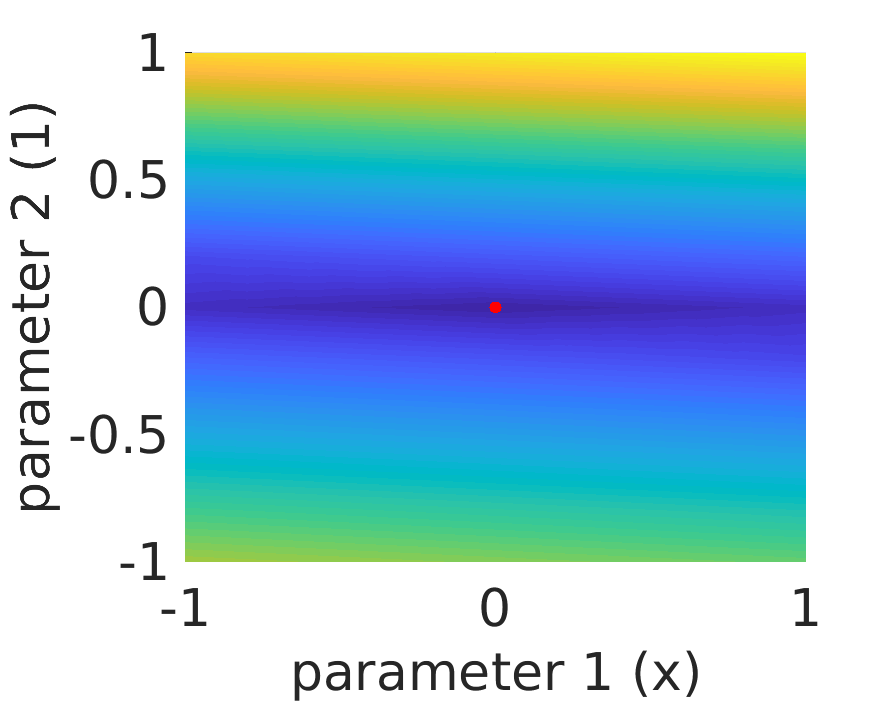}
    \includegraphics[height=0.195\textheight]{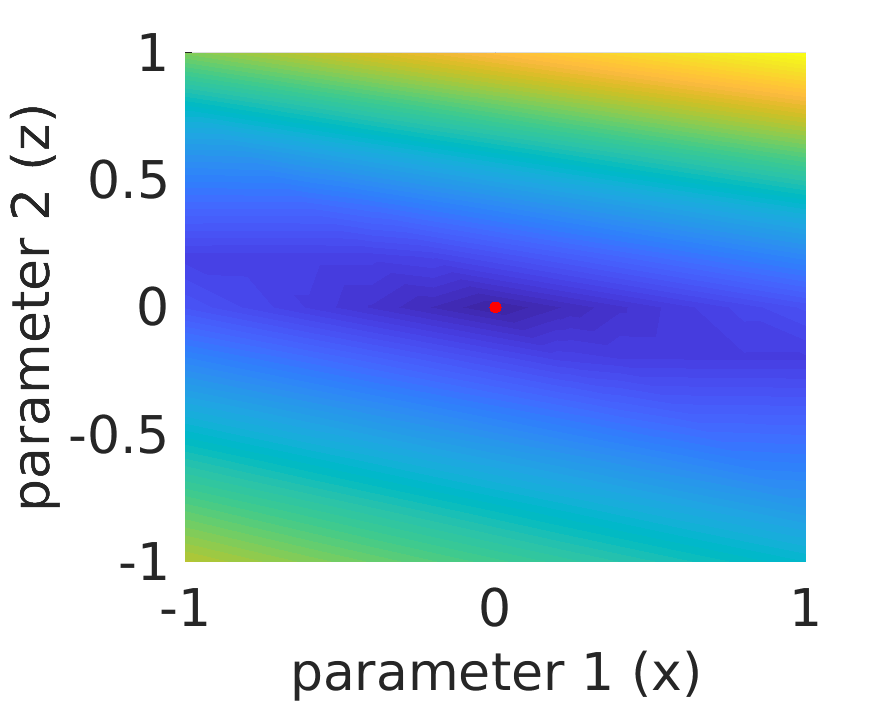}\\
    \includegraphics[height=0.195\textheight]{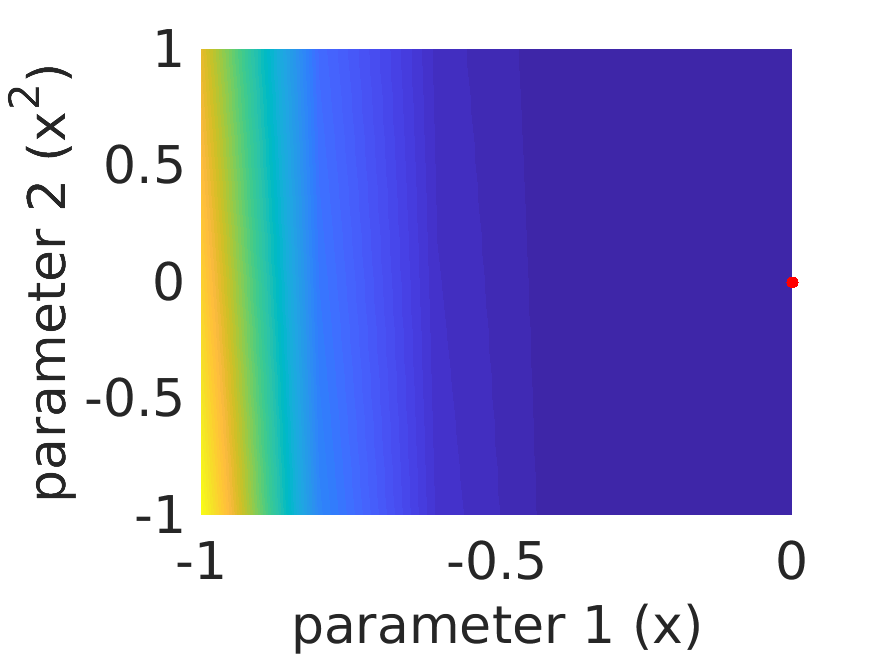}
    \includegraphics[height=0.195\textheight]{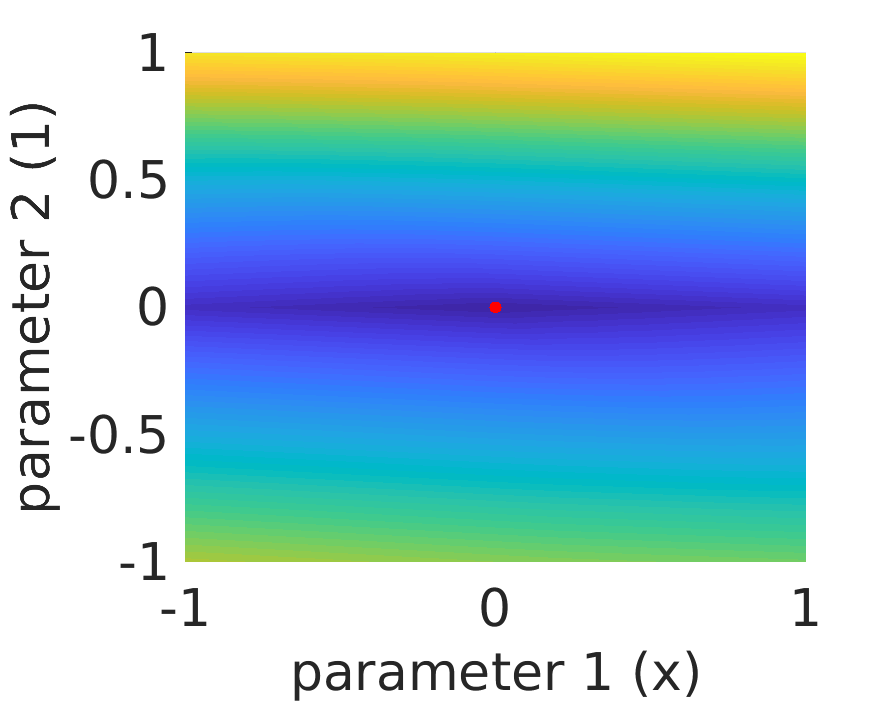}\\
    \includegraphics[height=0.195\textheight]{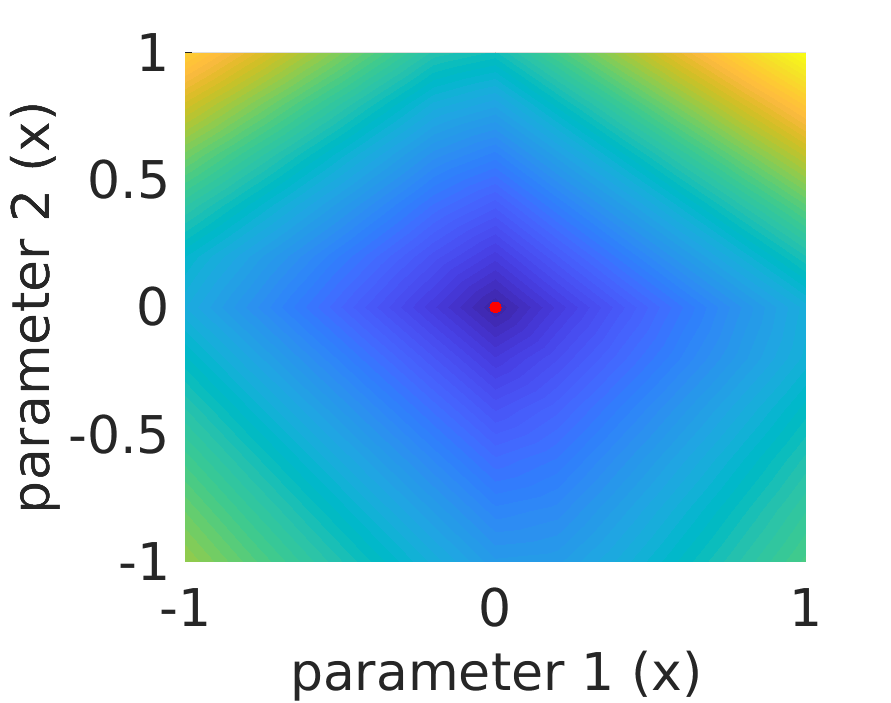}
    \includegraphics[height=0.195\textheight]{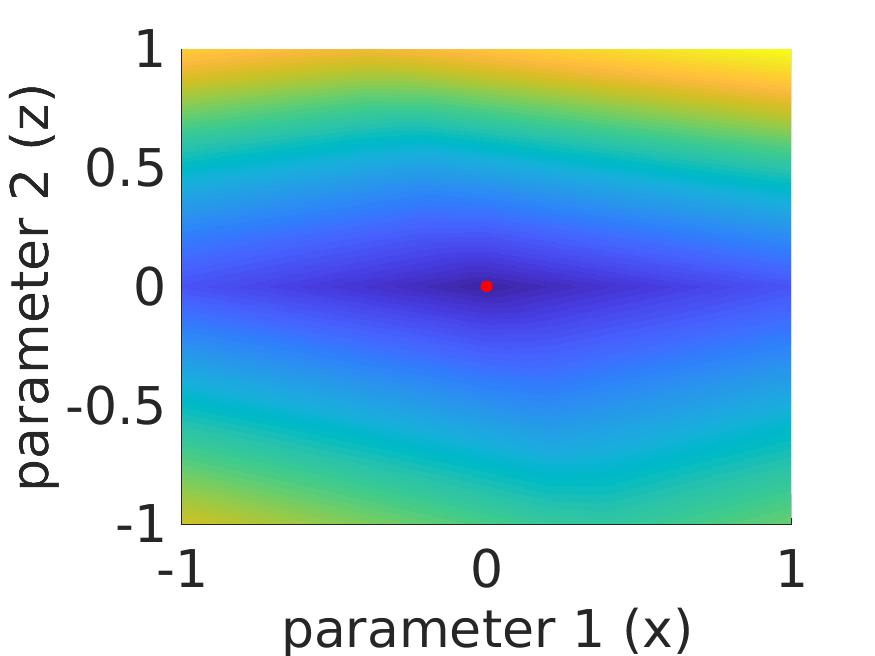}\\
    \includegraphics[height=0.195\textheight]{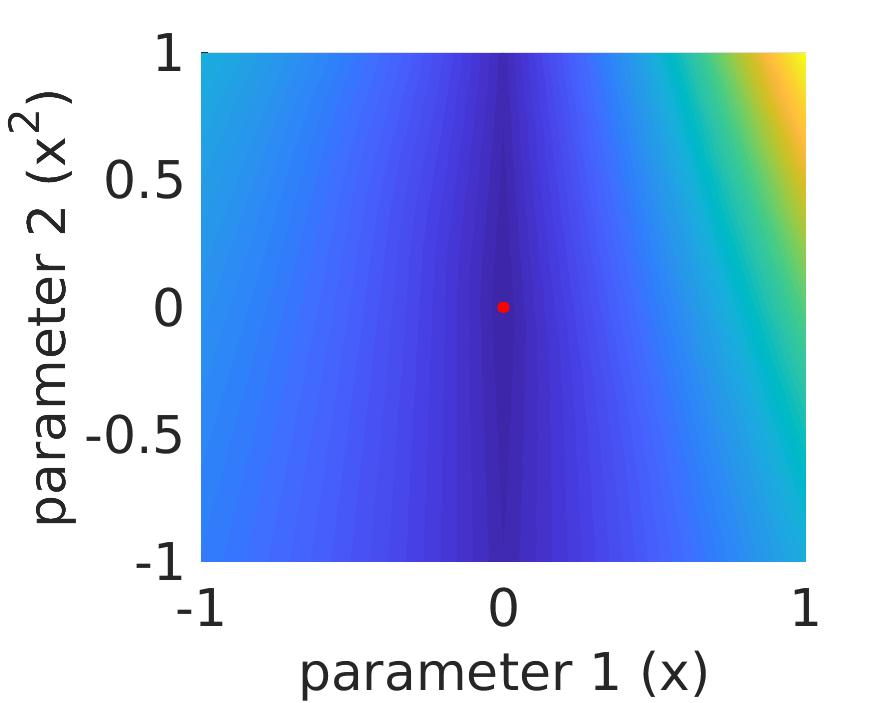}
    \includegraphics[height=0.195\textheight]{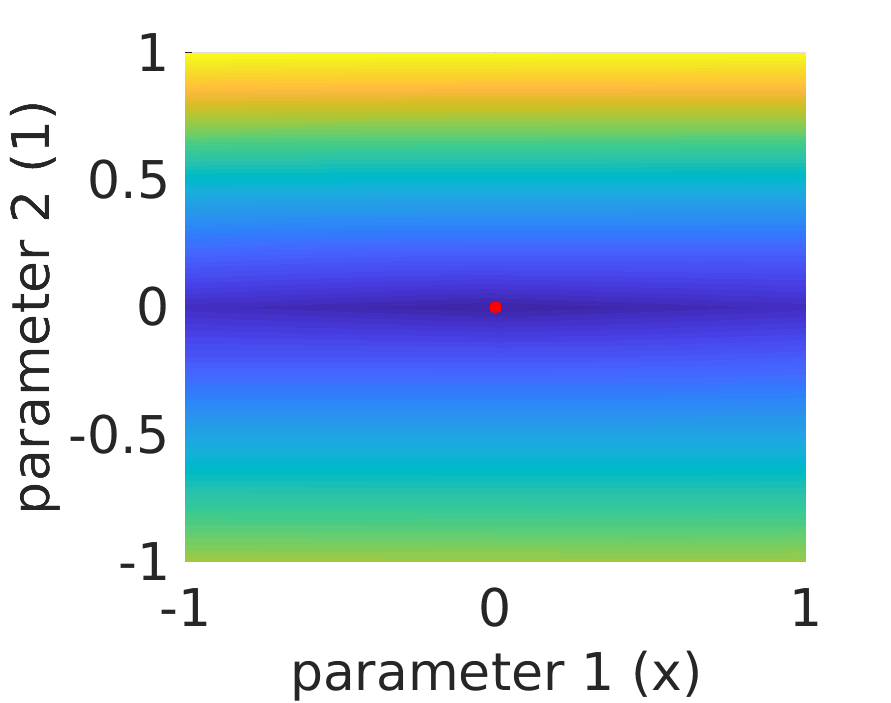}\\
    \includegraphics[height=0.195\textheight]{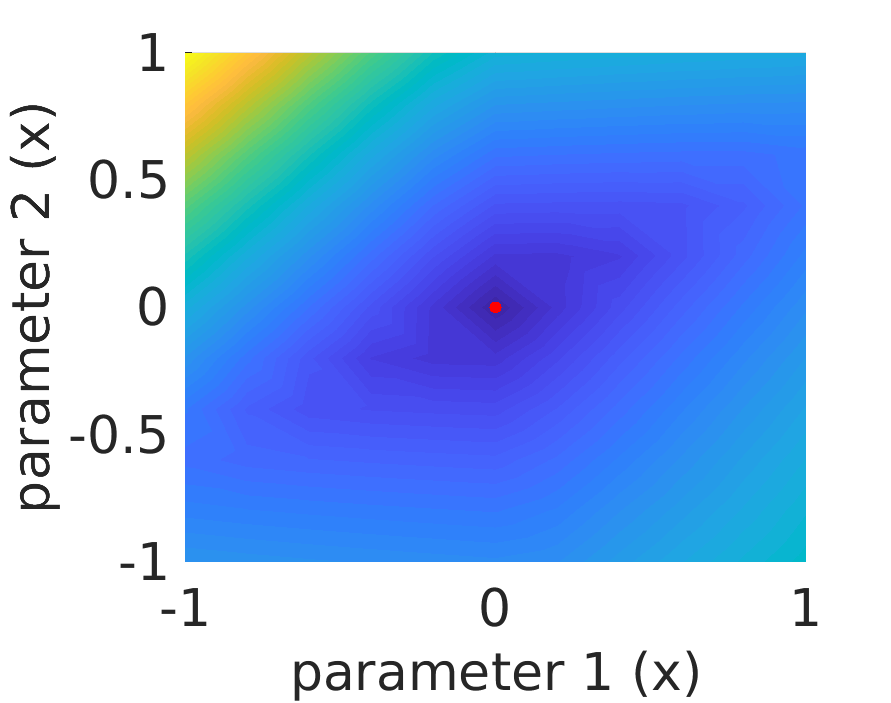}
    \includegraphics[height=0.195\textheight]{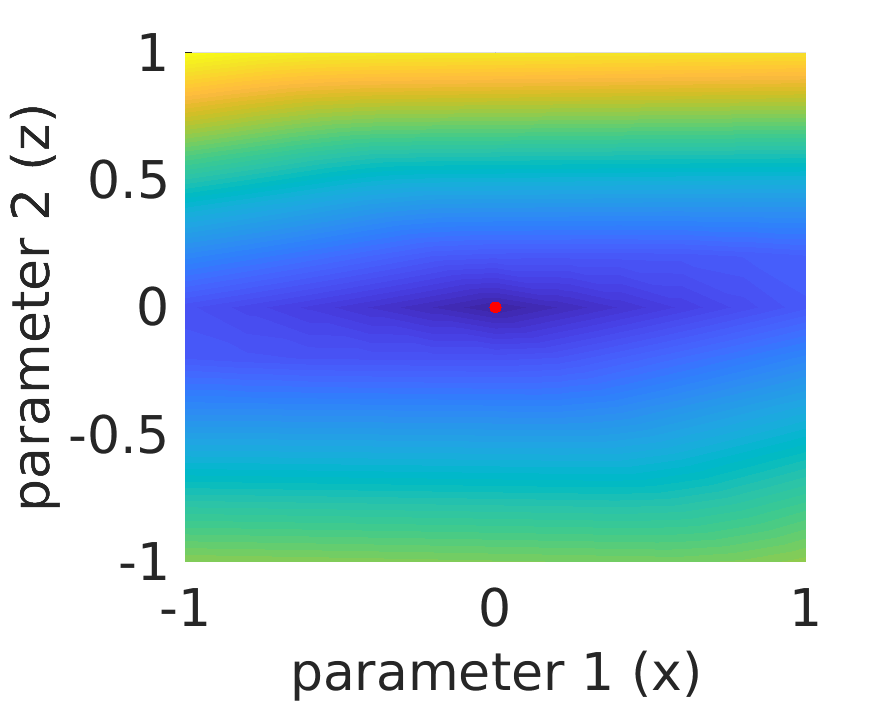}
    \caption{parameter 1 is $x$ on the $\ud x/\ud t$ equation.}
    \label{figure_identi_2}
\end{figure}

\begin{figure}
    \centering
    \includegraphics[height=0.195\textheight]{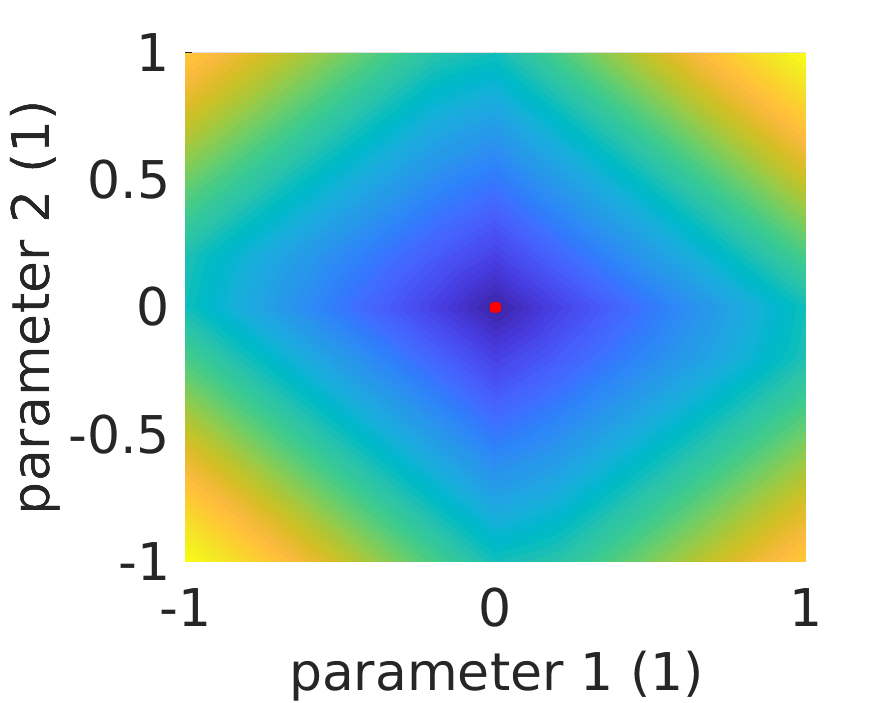}
    \includegraphics[height=0.195\textheight]{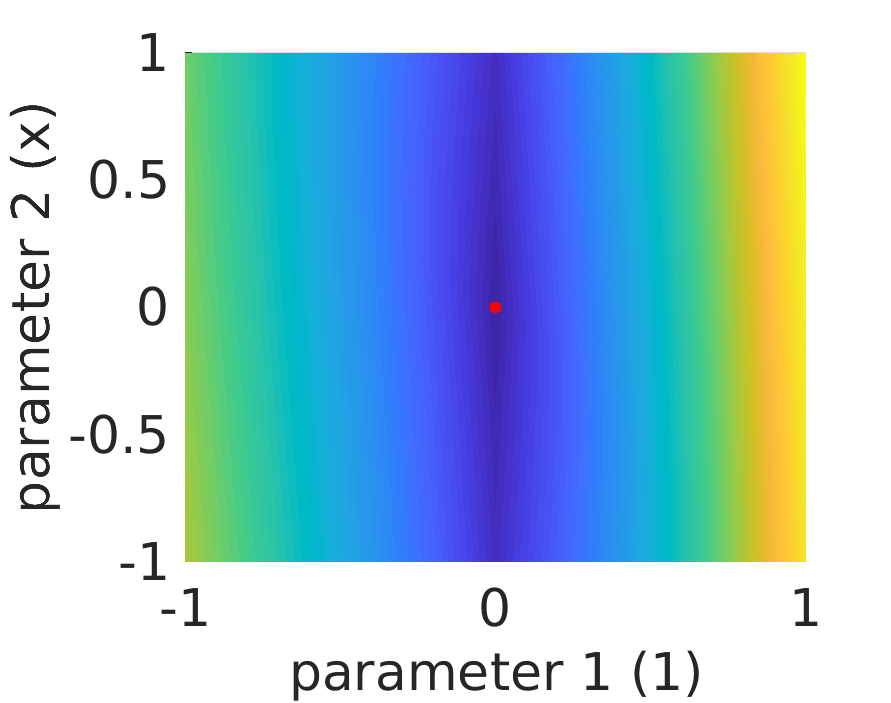}\\
    \includegraphics[height=0.195\textheight]{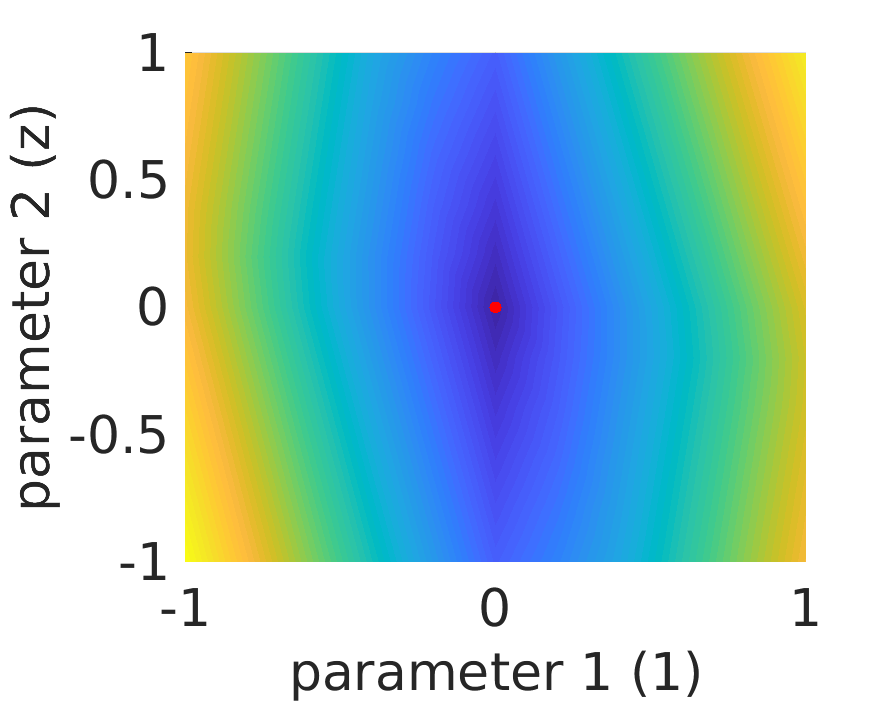}
    \includegraphics[height=0.195\textheight]{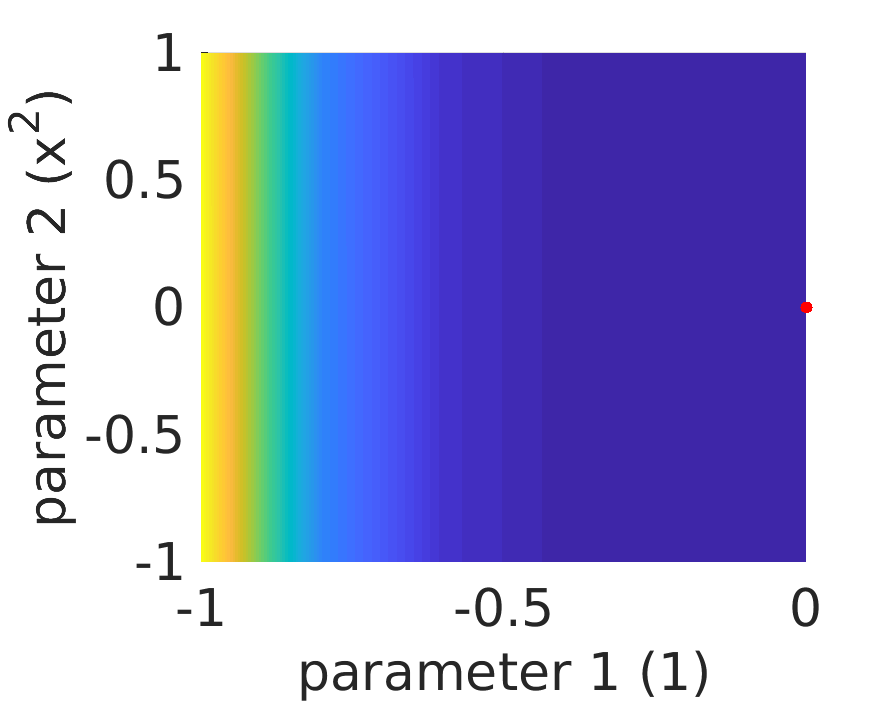}\\
    \includegraphics[height=0.195\textheight]{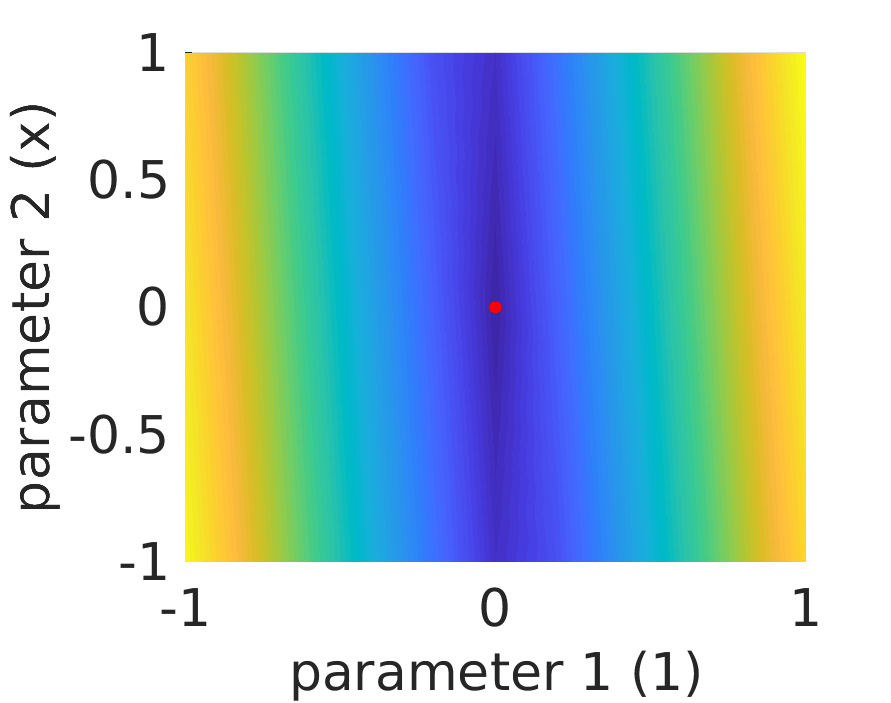}
    \includegraphics[height=0.195\textheight]{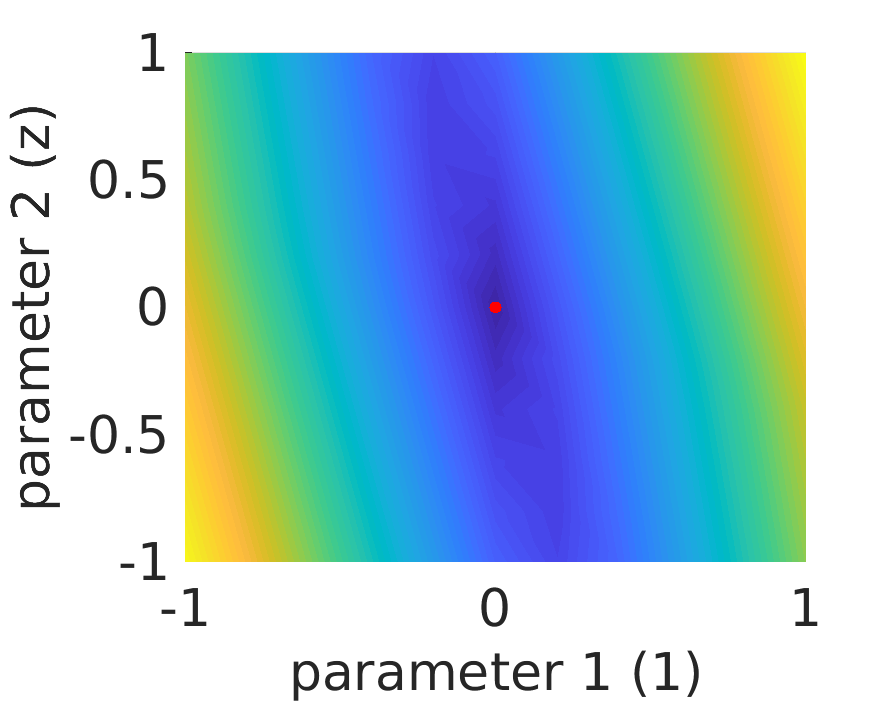}\\
    \includegraphics[height=0.195\textheight]{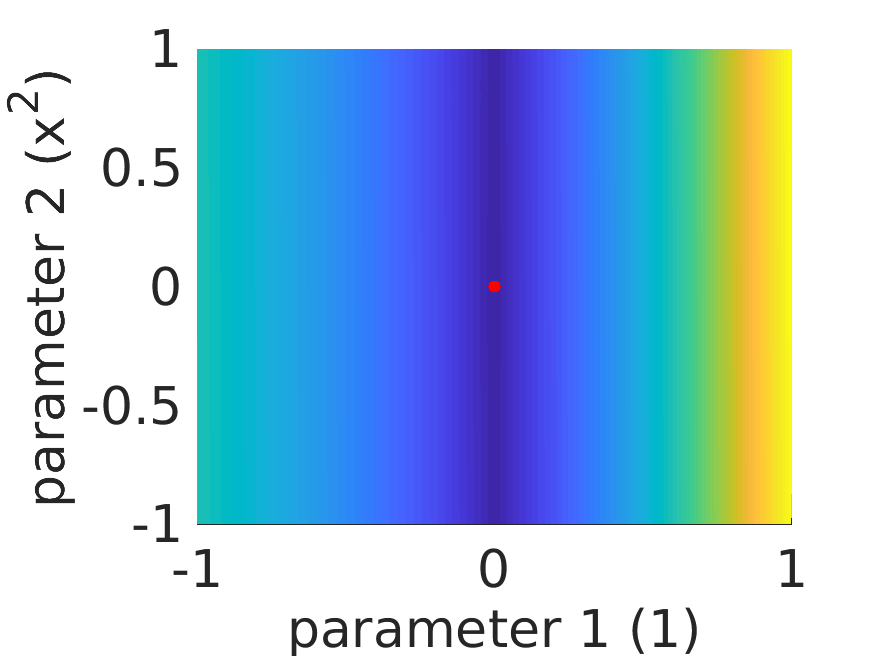}
    \includegraphics[height=0.195\textheight]{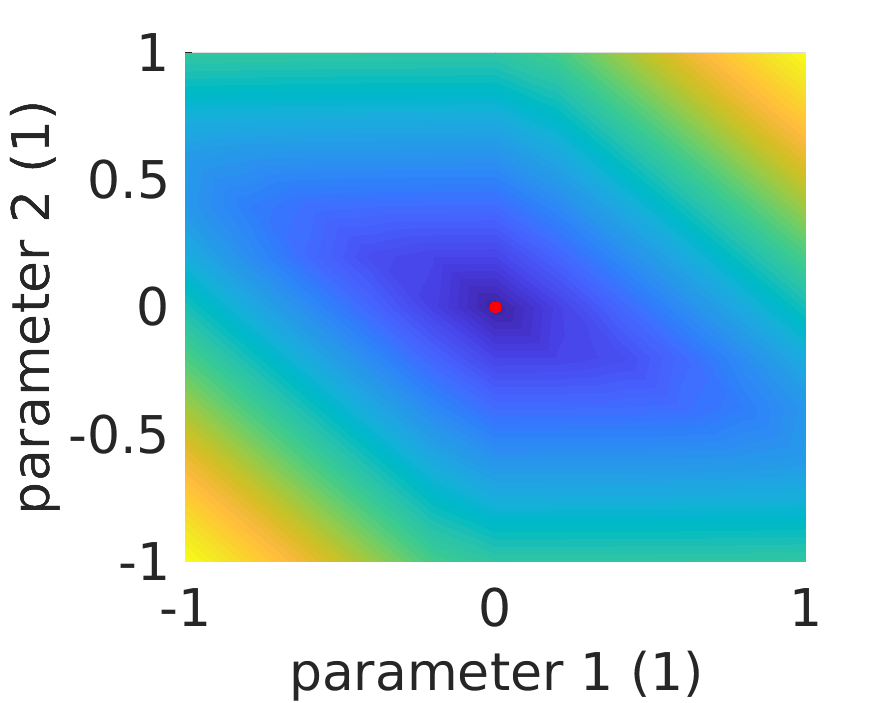}\\
    \includegraphics[height=0.195\textheight]{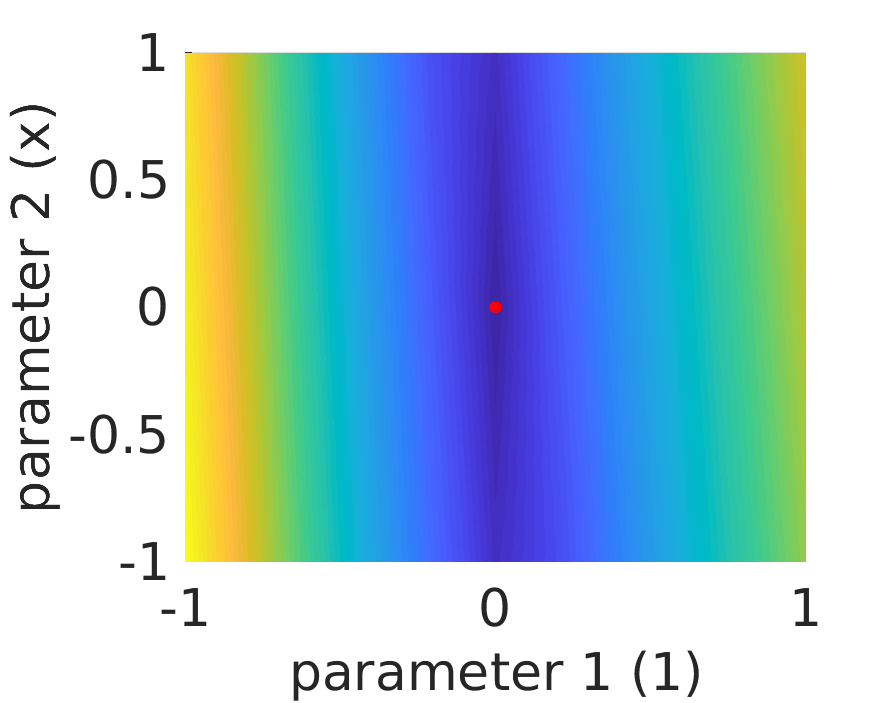}
    \includegraphics[height=0.195\textheight]{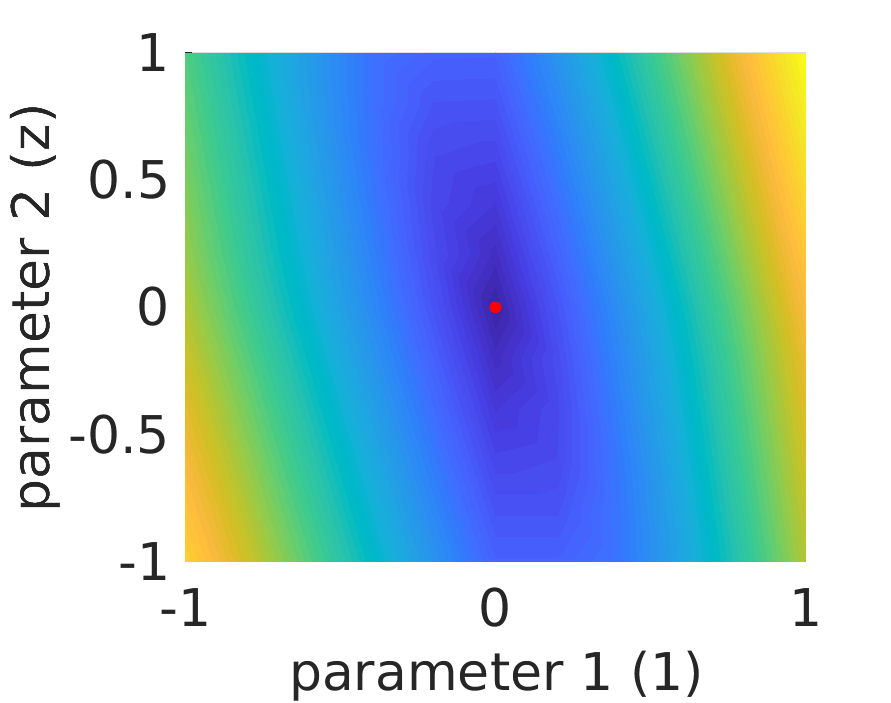}
    \caption{parameter 1 is $1$ on the $\ud y/\ud t$ equation.}
    \label{figure_identi_3}
\end{figure}

\begin{figure}
    \centering
    \includegraphics[height=0.195\textheight]{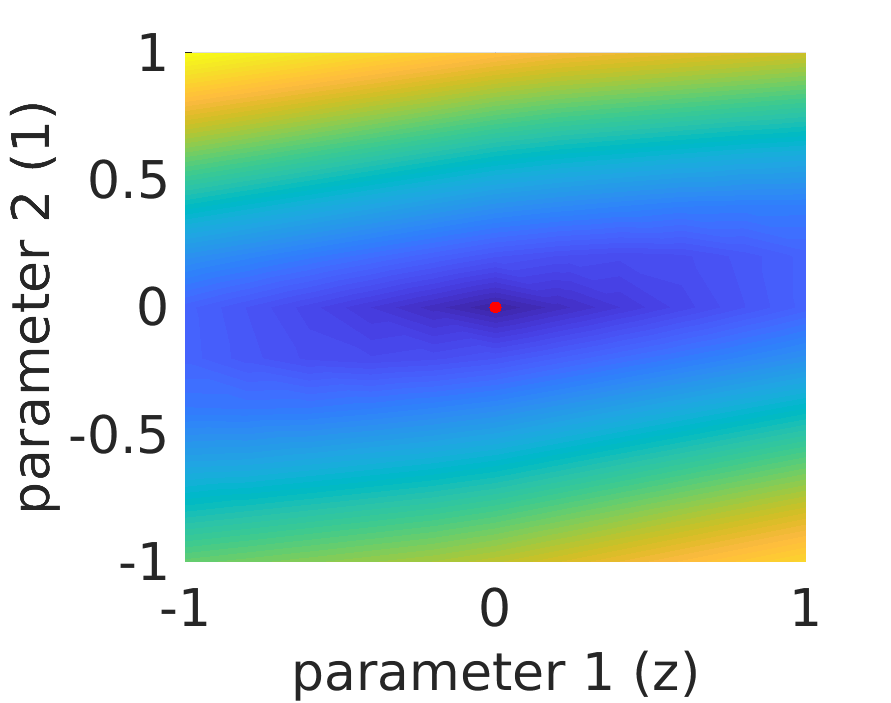}
    \includegraphics[height=0.195\textheight]{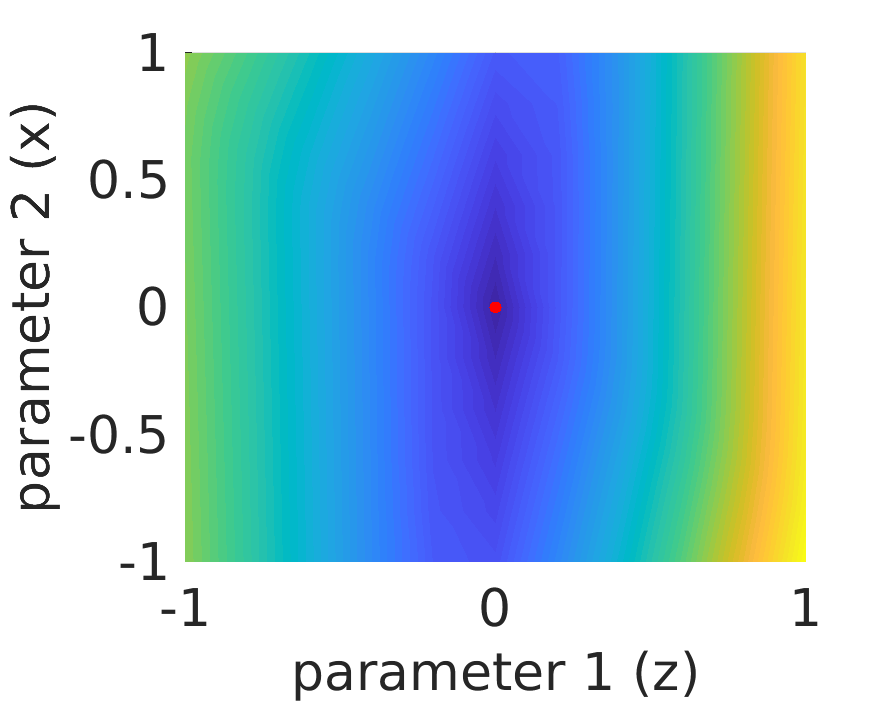}\\
    \includegraphics[height=0.195\textheight]{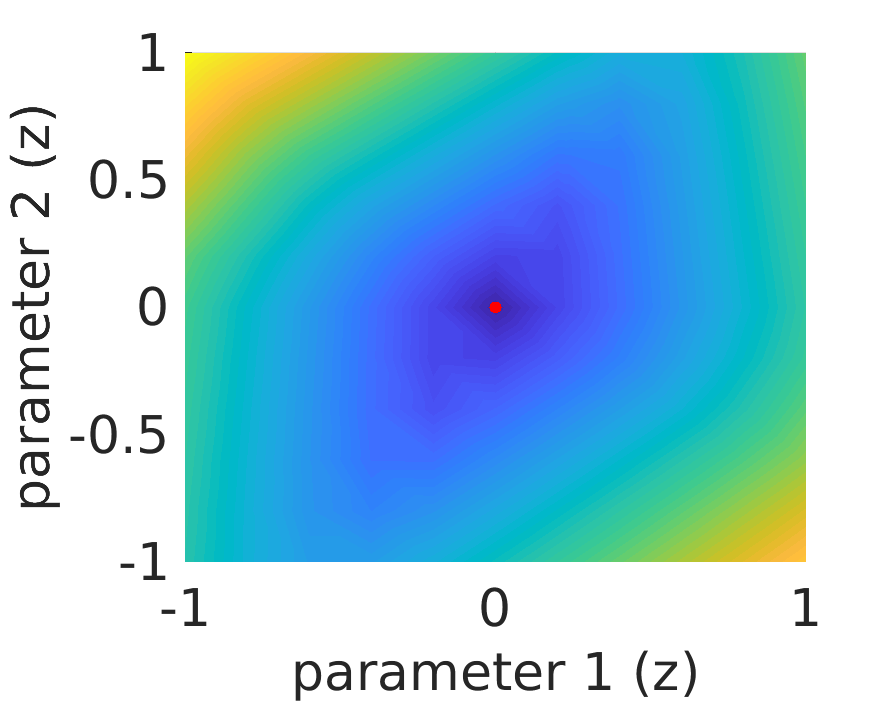}
    \includegraphics[height=0.195\textheight]{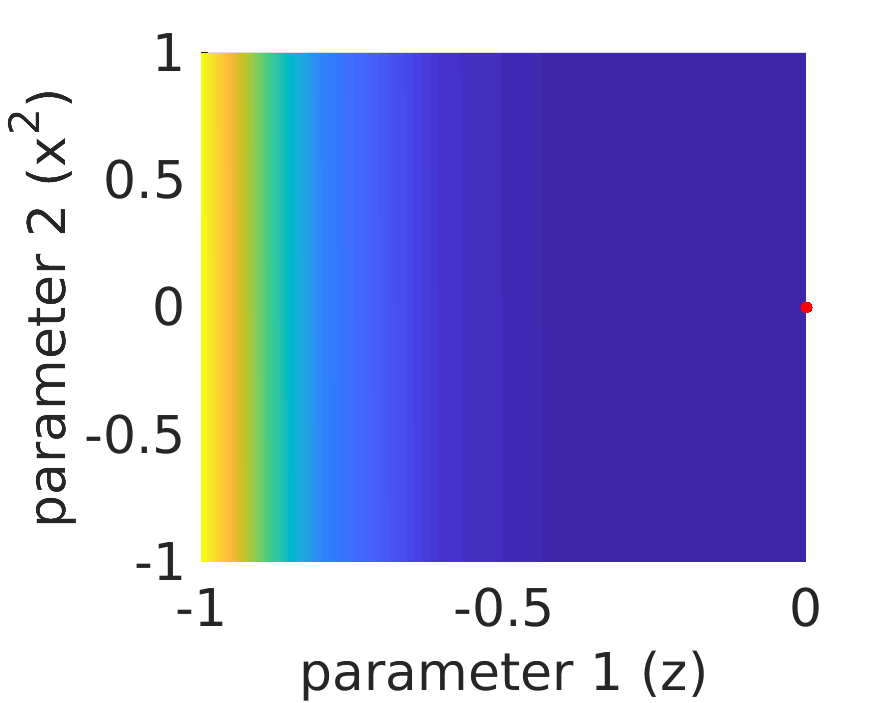}\\
    \includegraphics[height=0.195\textheight]{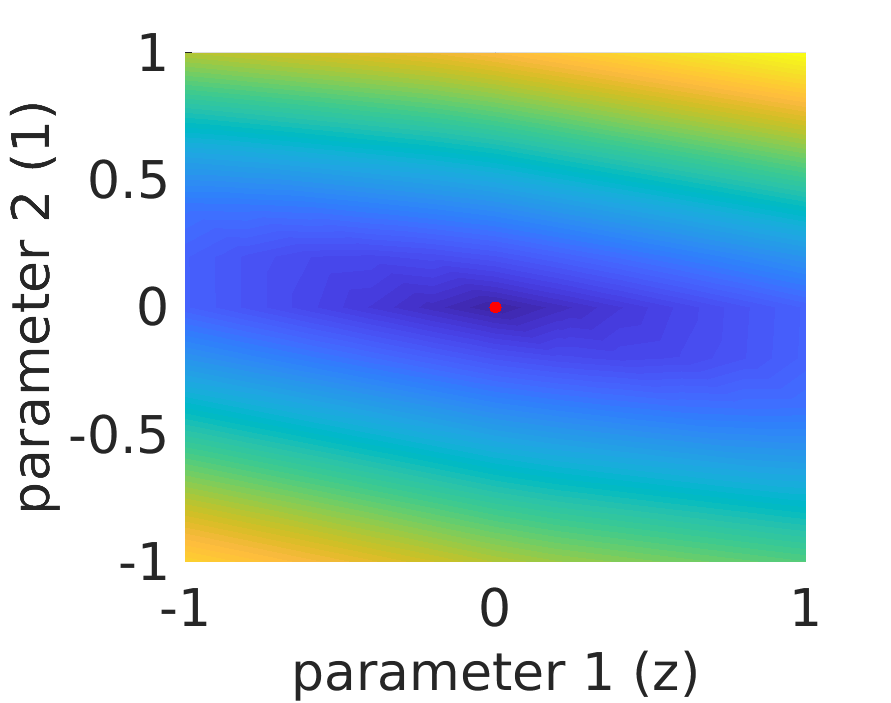}
    \includegraphics[height=0.195\textheight]{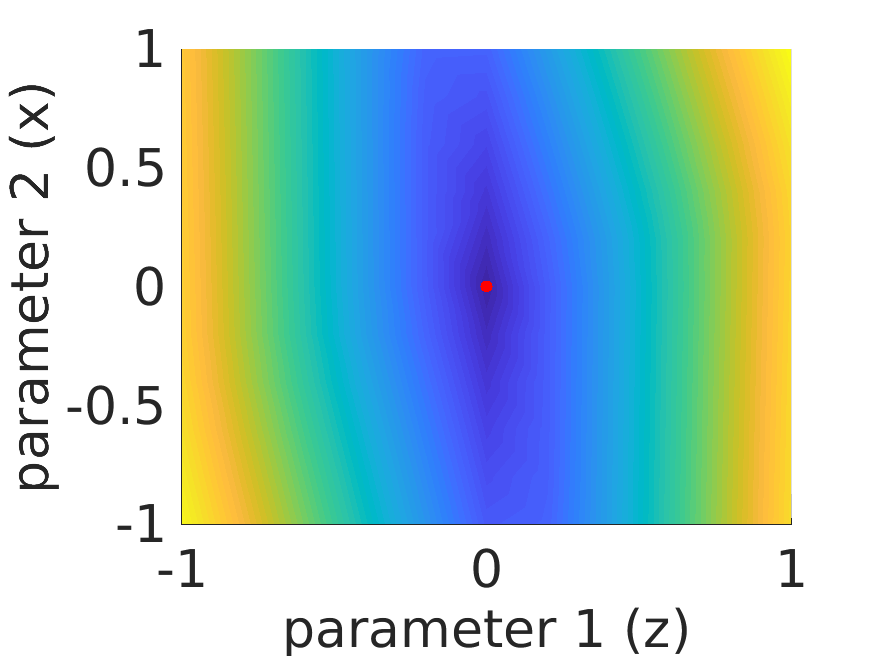}\\
    \includegraphics[height=0.195\textheight]{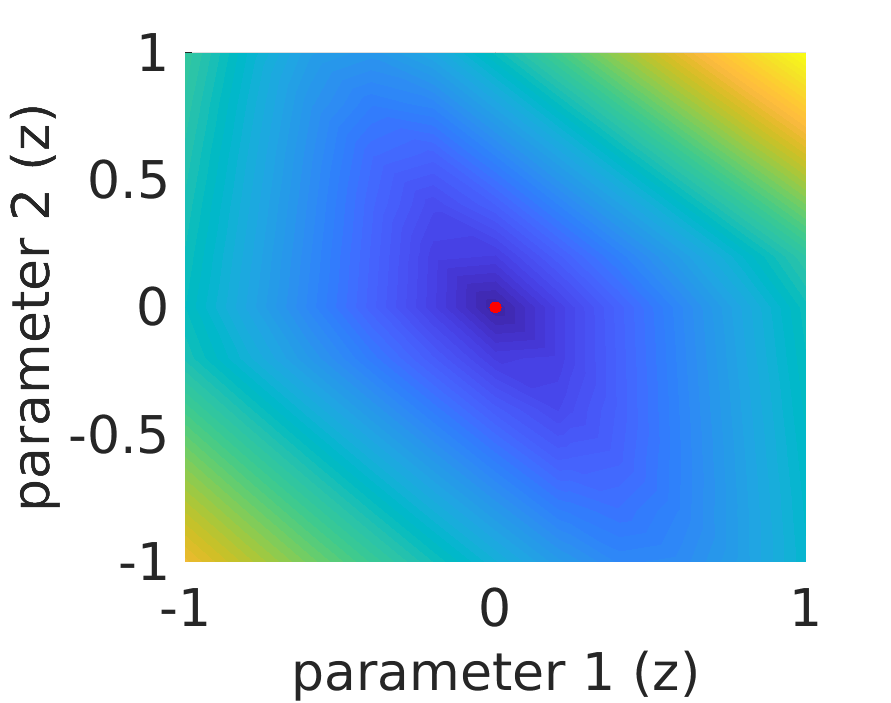}
    \includegraphics[height=0.195\textheight]{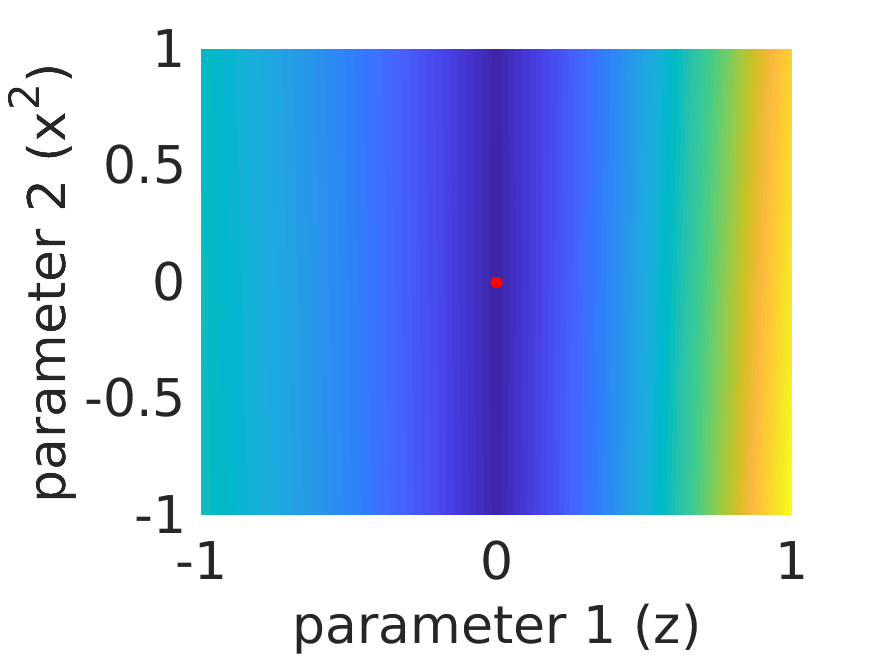}\\
    \includegraphics[height=0.195\textheight]{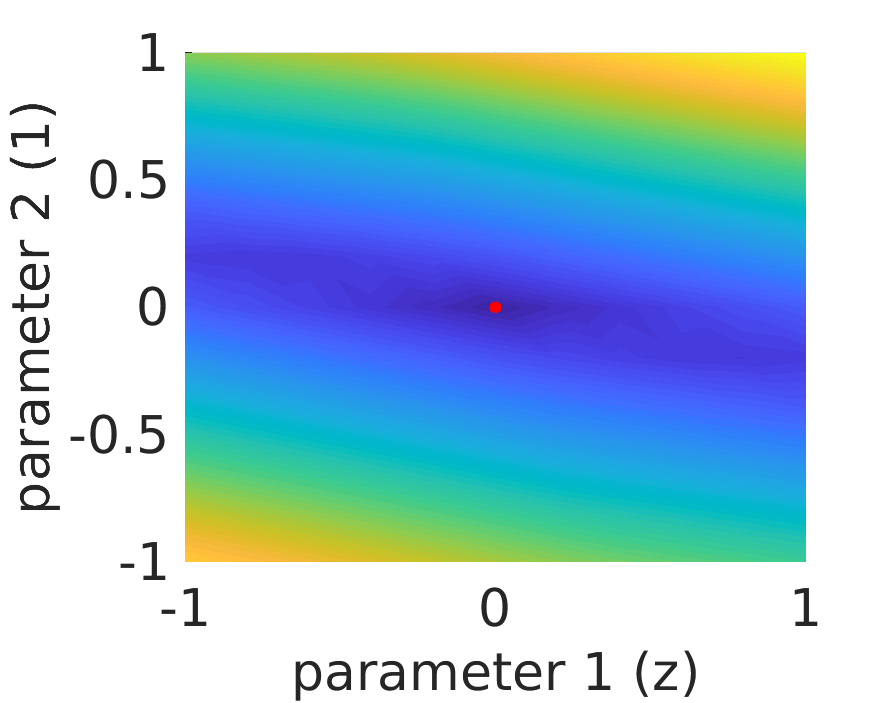}
    \includegraphics[height=0.195\textheight]{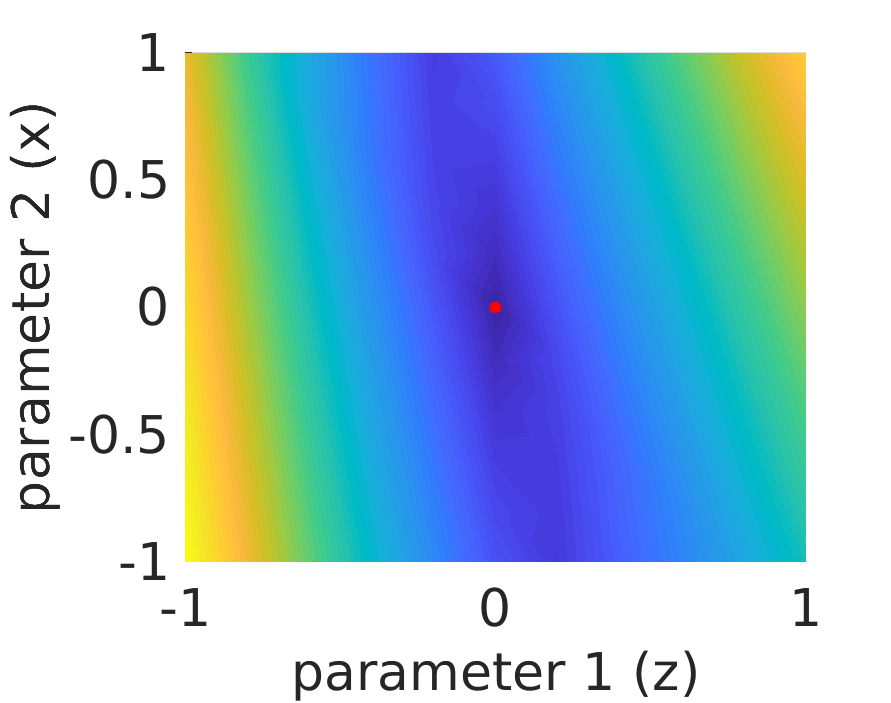}
    \caption{parameter 1 is $z$ on the $\ud z/\ud t$ equation.}
    \label{figure_identi_4}
\end{figure}

Consider the generic model (except that the right terms are fixed -- highlighted in red; $e_{n,01} = 0.5$, $R_{n,10} = 0.25$)
\begin{align}
    \frac{\ud x}{\ud t} &= p_{1,1} + p_{1,2} x \textcolor{red}{\,+\,e_{n,01}} y + p_{1,4} z + p_{1,5} x^2 \textcolor{red}{\,-\,R_{n,10}} xz,\\
    \frac{\ud y}{\ud t} &= p_{2,1} + p_{2,2} x \textcolor{red}{\,-\,e_{n,01}} y + p_{2,4} z + p_{2,5} x^2 \textcolor{red}{\,+\,R_{n,10}} xz,\\
    \frac{\ud z}{\ud t} &= p_{3,1} + p_{3,2} x \textcolor{red}{\,+\,e_{n,01}} y + p_{3,4} z \textcolor{red}{\,-\,R_{n,10}} xz.   
\end{align}
We will now add only two extra parameters (two of the black terms) at a time. Each Fig. \ref{figure_identi_1}-\ref{figure_identi_4} is obtained by picking one term (parameter 1, which is the first extra term in the system), and then study the cost function by adding another term (parameter 2, which is the second extra term in the system). We study this for all the possibles ``parameter 2''.

Take Fig. \ref{figure_identi_1}. Parameter 1 here is the term $1$ in the equation $\ud x / \ud t$, that is, $p_{1,1}$. This extra term is fixed for all subplots. Then parameter 2 (the second extra term) corresponds to (in order of subplots) $p_{1,2}, \, p_{1,4}, \, p_{1,5}, \, p_{2,1}, \, p_{2,2}, \, p_{2,4}, \, p_{2,5}, \, p_{3,1}, \, p_{3,2} \, p_{3,4}$.
Figures \ref{figure_identi_2}-\ref{figure_identi_4} follow the same logic. These four figures already show identifiability problems.

\clearpage 


\section{Predator-Prey}
\label{SI-predtorprey}

\begin{figure}[htb]
    \centering
    \includegraphics[width=.85\textwidth]{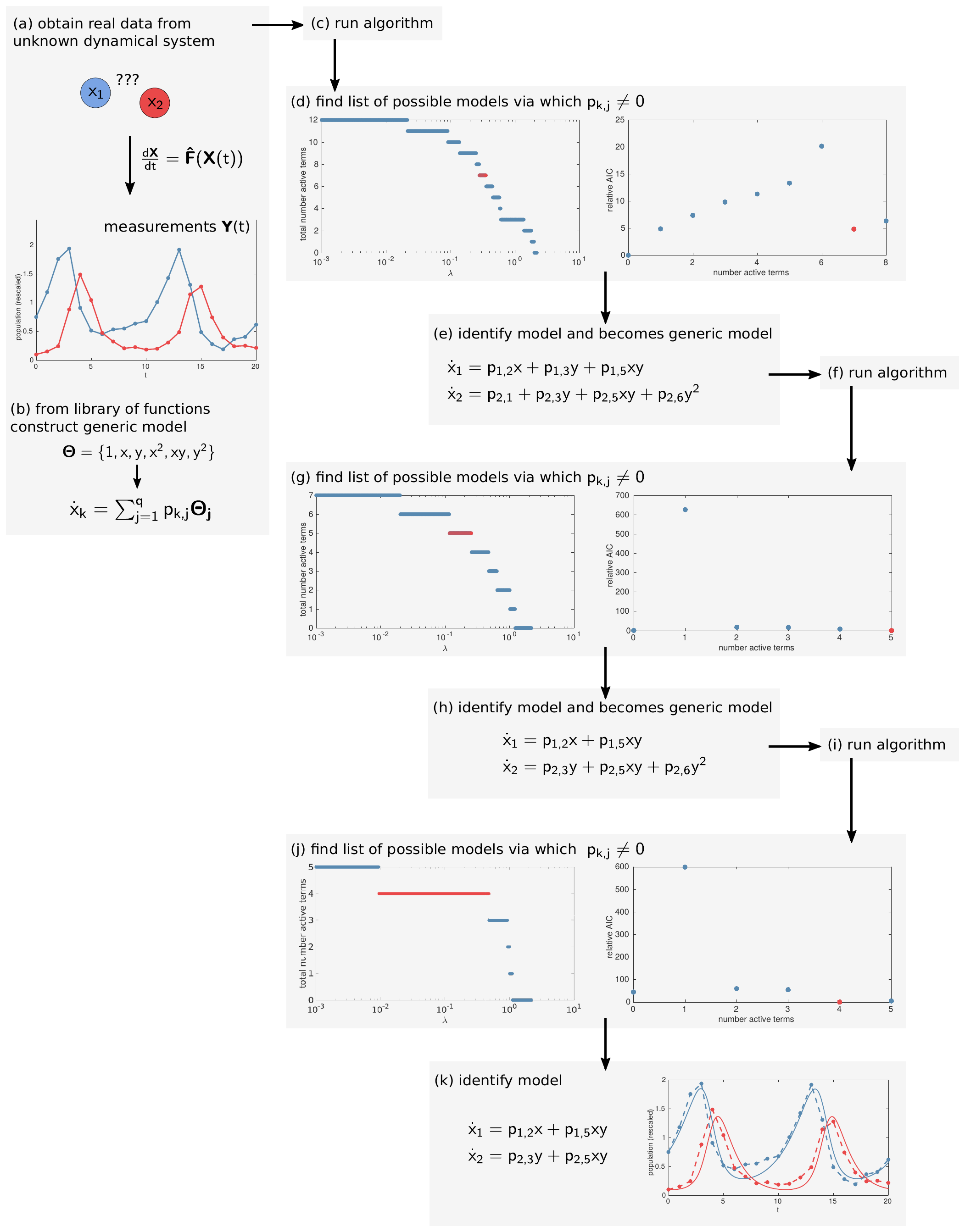}
    \caption{The recovery of the Lotka-Volterra system required an iterative formulation which consisted of down-selecting relevant monomials to describe the dynamics via AIC at the end of the variational annealing and start the algorithm again with less terms in the generic model description.}
    \label{figure_PredatorPrey}
\end{figure}

Although it is not very common to find \emph{pure} predator-prey interactions in nature, there is a classical set of data by the Hudson Bay company which corresponds the number of snowshoe hares and Canadian lynxes trapped in Canada, which in turn shows the relative population of both \cite{Odum1971}. The data is recorded yearly, so $\Delta t = 1$. We use data between 1900 and 1920, thus $N = 21$. In this particular case we really do not know the dynamics behind the system although we know that the snowshoe hare is the primary food of the lynx. Therefore, we can assume that we have a predator-prey system, and there is the classical Lotka-Volterra model to describe these type of dynamics. We consider $L = D = 2$ (Fig. \ref{figure_PredatorPrey}(a)). 
We build the library of functions with all the monomials up to degree two in two variables, and with it we construct our generic model (Fig. \ref{figure_PredatorPrey}(b)). We run our algorithm and varying $\lambda$ we obtain a list of possible models. By looking at the corresponding AIC values for each one, we find that the model with 7 active terms is the best one (Fig. \ref{figure_PredatorPrey}(d)). We now consider that our generic model is the resulting model with 7 active terms. Again, we run the algorithm to find that the best model is one containing only 5 terms (Fig. \ref{figure_PredatorPrey}(f-h)). We iterate this process, and run the algorithm considering the model with 5 active terms as the generic one. We find that the best model is the one containing 4 terms (Fig. \ref{figure_PredatorPrey}(i-k)). This identified model corresponds to the Lotka-Volterra one. Once we do only parameter estimation on it, we obtain the dynamical system shown in Fig. \ref{figure_PredatorPrey}(k). We compare the original data (dashed) with the resulting model (solid), which show an excellent match.

\clearpage


\section{$\alpha$ parameter in VA algorithm}
\label{SI-alpha}

We study how the parameter $\alpha$ used to increase the value of $R_f = R_{f,0} \alpha^{\beta}$ during the VA algorithm affects the recovery. We use the class Lorenz system,
\begin{align}
    \frac{\ud x}{\ud t} &= \sigma(y-x),\\
    \frac{\ud y}{\ud t} &= x(\rho-z) - y,\\
    \frac{\ud z}{\ud t} &= -\beta z + xy,
\end{align}
where $\sigma = 10$, $\rho = 28$, and $\beta = 8/3$. We numerically simulate the system using Runge-Kutta 4th order and a time step of $\Delta t = 0.01$, producing time-series similar to the experimental data set. 
We add some error modeled as additive Gaussian noise of mean zero and standard deviation $\omega=0.01$. Therefore, the measurement function is $\mathbf{h}(\mathbf{X}) = \mathbf{X} + \mathcal{N}(0,\omega)$.
We consider $N=501$, and $y$ to be the hidden variable.

As we increase $\alpha$ the recovery rate decreases, and for $\alpha \geq 1.3$ the recovery is 0\% (Table \ref{table_alpha}).

\begin{table}[H]
\centering
\caption{Recovery rates for varying $\alpha$.}
\label{table_alpha}
    \begin{tabular}{cc}
    $\alpha$ & recovery rate (\%) \\
    \hline
    1.1 & 93 \\ \hline
    1.2 & 87 \\ \hline
    1.25 & 20 \\ \hline
    1.3 & 0 \\ \hline
    1.4 & 0 \\ \hline
    1.5 & 0  
    \end{tabular}
\end{table} 